\documentclass[12pt]{amsart}
\usepackage{amsmath}
\usepackage{epsfig}
\usepackage{latexsym}
\usepackage{amscd}

\newcommand{\bd}{\partial}
\newcommand{\hm}{homeomorphic}
\newcommand{\R}{\ensuremath{\mathbf{R}}}
\newcommand{\RR}{\ensuremath{\mathbf{R}^2}}
\newcommand{\RRR}{\ensuremath{\mathbf{R}^3}}
\newcommand{\FF}{\ensuremath{\mathcal{F}}}
\newcommand{\sbs}{\subseteq}
\newcommand{\irr}{irreducible}
\newcommand{\inc}{incompressible}
\newcommand{\birr}{$\bd$-\irr}
\newcommand{\rirr}{\RR-\irr}
\newcommand{\pl}{parallel}
\newcommand{\bpl}{$\bd$-\pl}
\newcommand{\ei}{end \irr}
\newcommand{\eei}{eventually \ei}
\newcommand{\tm}{3-manifold}
\newcommand{\nh}{null homotopic}
\newcommand{\p}{^{\prime}}
\newcommand{\er}{end reduction}
\newcommand{\ra}{\rightarrow}
\newcommand{\ns}{\emptyset}
\newcommand{\inte}{int \,}
\newcommand{\qe}{quasi-exhaustion}
\newcommand{\dpo}{\ensuremath{D_{n+2,2p+1}}}
\newcommand{\dpt}{\ensuremath{D_{n+2,2p+2}}}
\newcommand{\hs}{\ensuremath{\RR\times[0,\infty)}}
\newcommand{\mgp}{minimal general position}
\newcommand{\ga}{\ensuremath{\gamma}}
\newcommand{\be}{\ensuremath{\beta}}

\newcommand{\al}{\ensuremath{\alpha}}
\newcommand{\ze}{\ensuremath{\zeta}}

\newcommand{\rh}{\ensuremath{\rho}}

\newcommand{\ta}{\ensuremath{\tau}}

\newcommand{\n}{^{-1}}

\newcommand{\ep}{\ensuremath{\varepsilon}}
\newcommand{\tht}{\ensuremath{\theta}}

\newcommand{\pp}{^{\prime\prime}}

\newcommand{\ka}{\ensuremath{\kappa}}

\newcommand{\de}{\ensuremath{\delta}}

\newcommand{\la}{\ensuremath{\lambda}}

\newcommand{\PP}{\ensuremath{\mathcal{P}}}
\newcommand{\QQ}{\ensuremath{\mathcal{Q}}}
\newcommand{\RRRR}{\ensuremath{\mathcal{R}}}
\newcommand{\SSSS}{\ensuremath{\mathcal{S}}}
\title[End reductions and covering translations]
{End reductions and covering translations \\ of contractible open 
3-manifolds}
\author{Robert Myers} 
\address{Department of Mathematics, Oklahoma State University, Stillwater, OK 74078}
\email{myersr@math.okstate.edu}
\thanks{This research was supported in part by NSF grant DMS-0072429.}

\newtheorem{thm}{Theorem}[section]
\newtheorem{prop}[thm]{Proposition}
\newtheorem{cor}[thm]{Corollary}
\newtheorem{lem}[thm]{Lemma}

\begin{document}

\begin{abstract} 
This paper uses Brin and Thickstun's theory of end reductions of 
non-compact 3-manifolds to study groups of covering translations 
of irreducible contractible open 3-manifolds $W$ which are not 
homeomorphic to \RRR. We associate to $W$ an object $\mathcal{S}(W)$ 
called the \textit{simplicial complex of minimal \rirr\ end reductions 
of $W$}. Whenever $W$ covers another 3-manifold the group of covering 
translations is isomorphic to a fixed point free group of automorphisms 
of $\mathcal{S}(W)$. We apply this result to give uncountably many 
examples of \rirr\ such $W$ which cover 3-manifolds with infinite 
cyclic fundamental group and non-trivially cover only 3-manifolds 
with infinite cyclic fundamental group. We also give uncountably many 
examples of \rirr\ $W$ which are not eventually end irrreducible and 
do not non-trivially cover any 3-manifold.  
\end{abstract}

\maketitle

%Section 1
\section{Introduction}

The universal covering conjecture states that whenever $M$ is a 
closed, connected, orientable, irreducible 3-manifold with 
infinite cyclic fundamental group, then its universal covering 
space must be homeomorphic to \RRR. Valentin Po\'{e}naru has 
announced a proof of this conjecture \cite{Po}. His proof is the 
outgrowth of a massive project which investigates the structure 
of simply connected 3-manifolds. He translates the problem into 
another problem in higher dimensions, solves that problem, and 
then uses the result to solve the original 3-dimensional problem. 
His methods are highly original and perhaps not that familiar to 
most 3-manifold toplogists. (See \cite{Ga} for an introduction to  
some of Po\'{e}naru's techniques.) 

The present paper continues an investigation of an alternative approach to 
the universal covering conjecture which uses more traditional 
3-manifold techniques. An irreducible contractible open 3-manifold 
which is not homeomorphic to \RRR\ is called a \textit{Whitehead manifold}. 
Thus one wants to show that Whitehead manifolds cannot cover 
compact 3-manifolds. In \cite{My genusone} I showed that genus one 
Whitehead manifolds (those which are monotone unions of solid tori) 
cannot non-trivially cover any 3-manifold. Wright then showed in 
\cite{Wr} that the same result is true for the much larger class 
of \textit{eventually end irreducible} Whitehead manifolds; his proof 
was based on the interplay of two new ideas: the ``Ratchet Lemma'' 
(which uses the assumption of eventual end irreducibility) and 
the ``Orbit Lemma'' (which does not). 
Tinsley and Wright \cite{TW} then extended this 
result to certain examples of Whitehead manifolds which are not 
eventually end irreducible; these manifolds contain a family of 
disjoint proper planes which split them into an infinite collection 
of certain genus one Whitehead manifolds. Their proof that these 
Whitehead manifolds cannot non-trivially cover any other 3-manifolds 
is based on the Orbit Lemma and an extension of the Ratchet Lemma 
called the Special Ratchet Lemma which applies to Whitehead manifolds 
which split along proper planes into eventually end irreducible 
Whitehead manifolds. They also gave other examples of manifolds of this 
type which could non-trivially cover other non-compact 3-manifolds; 
however, their methods did not rule out their covering compact 3-manifolds. 
In \cite{My free} I used the Special Ratchet Lemma and the Orbit Lemma 
to show that certain Whitehead manifolds which covered non-compact 
3-manifolds with free fundamental groups could cover only 3-manifolds 
with free fundamental groups, and so could not cover compact 3-manifolds. 
These examples are different from those of \cite{TW} but like them 
they can be split along a family of proper planes into 
eventually end irreducible Whitehead manifolds. One can associate to 
this decomposition a tree whose vertices correspond to the eventually 
end irreducible Whitehead manifolds and whose edges correspond to the 
planes. For these examples I showed that the decomposition is unique up 
to ambient isotopy. Thus the isotopy class of any self-homeomorphism 
induces an automorphism of the tree. When the self-homeomorphism is 
a non-trivial covering translation the Special Ratchet Lemma and the 
Orbit Lemma can be used to show that the automorphism does not fix 
any points of the tree. It follows that the group of covering translations 
is isomorphic to the corresponding group of automorphisms and that this 
group is free. 

This paper explores the situation when the Whitehead manifold is not 
eventually end irreducible and cannot be split into eventually end 
irreducible pieces by a collection of proper planes. It turns out that 
even in this situation the manifold still contains some very useful 
open submanifolds which are eventually end irreducible. They are called 
\textit{end reductions}. They were developed by Brin and Thickstun 
\cite{BT} as a tool in their study of the ends of non-compact 3-manifolds. 
An elementary example of an end reduction is one of the eventually 
end irreducible pieces into which the collection of planes breaks a 
Whitehead manifold of the type described above; such an example has 
topological boundary consisting of the planes adjacent to it. In general, 
however, an end reduction is embedded in a very contorted manner in the 
manifold containing it and has a very complicated topological boundary.    

In this paper an extension of the Special Ratchet Lemma called the 
``End Reduction Ratchet Lemma'' is proven. This result and the Orbit 
Lemma are used to prove that a non-trivial covering translation of 
a Whitehead manifold cannot fix the isotopy class of an end reduction. 
The isotopy in this statement need not be an ambient isotopy; it can 
refer to a non-ambient isotopy of embedding maps. 

In order to apply this theorem one needs some structure on the set 
of isotopy classes of end reductions. We define an object called the 
\textit{simplicial complex of minimal \rirr\ end reductions} of 
a Whitehead manifold $W$, denoted $\mathcal{S}(W)$. An end reduction $V$ 
of $W$ is \textit{minimal} if the only end reductions of $W$ that it 
contains are isotopic to it. It is \rirr\ if it does not contain 
any proper planes which split it in a non-trivial fashion. Two isotopy 
classes of end reductions are joined by an edge if there is an 
isotopically unique \rirr\ end reduction $E$  containing representatives 
$V_0$ and $V_1$ of the vertices such that every end reduction of $W$ 
which is contained in $E$ is isotopic to $V_0$, $V_1$, or $E$. Simplices 
of higher dimension are defined by an inductive extension of this 
definition. 

In \cite{ST} Scott and Tucker gave an example of a Whitehead manifold 
which is not eventually end irreducible, which is an infinite cyclic 
covering space of a non-compact 3-manifold, and which was purported to 
be \rirr. Unfortunately this example contained a mistake which resulted 
in its not being \rirr. In \cite{My r2p2} I showed how to modify this 
example to make it \rirr. In the present paper I determine $\mathcal{S}(W)$ 
for this modified example and uncountably many non-homeomorphic variations 
of it. In all cases $\mathcal{S}(W)$ is a triangulation of the real line. 
It follows that any 3-manifold which it non-trivially covers must have 
infinite cyclic fundamental group. A further variation produces 
uncountably many examples which cannot non-trivially cover any 3-manifold. 

The paper is organized as follows. In section 2 we review Brin and 
Thickstun's theory of end reductions and prove that any end reduction 
of a Whitehead manifold is a Whitehead manifold. 
End reductions are associated to certain compact subsets; we speak of an 
end reduction \textit{at} the compact set. In this section  
we also prove the technically 
convenient result that any end reduction can be chosen to be an end 
reduction at a knot. In section 3 we prove the End Reduction Ratchet 
Lemma. In section 4 we apply it and the Orbit Lemma to show that 
a non-trivial covering translation cannot fix the isotopy class of 
an end reduction. In the next two sections some further results about end 
reductions are proven which will be needed later. In section 5 it is 
shown that the end reductions associated to nested compact sets have 
representatives which are themselves nested; moreover, the smaller 
end reduction is an end reduction of the larger end reduction. In section 6 
it is shown that any \rirr\ end reduction can be isotoped off any 
non-trivial plane. Section 7 collects some standard facts about gluing 
together 3-manifolds in such a way as to achieve such properties as being 
incompressible, irreducible, $\bd$-irreducible, anannular, or atoroidal. 
Section 8 defines certain compact 3-manifolds which will be used in our 
construction and proves that they have certain properties we will need. 
In section 9 we construct the family $\mathcal{F}$ of Whitehead manifolds 
$W$ we will be considering. In section 10 we prove that these manifolds 
are \rirr. We also associate to each finite set of integers $\mathcal{P}$ a 
certain open subset $V^{\mathcal{P}}$ of $W$. It is a Whitehead manifold. 
We prove that it is \rirr\ if and only if the integers are consecutive. 
In section 11 we prove that $V^{\mathcal{P}}$ is a genus one Whitehead 
manifold if and only if $\mathcal{P}$ has a single element. In section 12 
we prove that each $V^{\mathcal{P}}$ is an end reduction of $W$. 
In section 13 we prove that each \rirr\ end reduction of $W$ is isotopic 
to one of the $V^{\mathcal{P}}$. In section 14 we prove that 
$V^{\mathcal{P}}$ and $V^{\mathcal{Q}}$ are isotopic if and only if 
$\mathcal{P}=\mathcal{Q}$. In section 15 we introduce a variation into 
our construction and classify the resulting uncountably many 
variations of those $V^{\mathcal{P}}$ of genus one up to homeomorphism. 
In section 16 we determine $\mathcal{S}(W)$ for the family of 
Whitehead manifolds we have constructed. The vertices are the isotopy 
classes of those $V^{\mathcal{P}}$ for which $\mathcal{P}$ has a 
single element $p$. The edges are the isotopy classes of those 
$V^{\mathcal{P}}$ for which $\mathcal{P}=\{p,p+1\}$. It follows that 
$\mathcal{S}(W)$ is a triangulation of $\mathbf{R}$. In section 17 we 
apply these results to show that uncountably many of these manifolds 
cover only 3-manifolds with infinite cyclic fundamental group and 
uncountably many of them do not non-trivially cover any 3-manifold.

%Section 2
\section{End reductions} 

In this section we review the Brin-Thickstun theory of end reductions 
of irreducible open 3-manifolds developed in \cite{BT}. We work in 
somewhat less generality by restricting attention to connected, 
orientable, irreducible open 3-manifolds $W$ which are not \hm\ to \RRR. 
Each such $W$ is contained in the class $\mathcal{Z}$ of 3-manifolds 
considered in \cite{BT}. A compact, connected 3-manifold $J\sbs W$ is 
\textit{regular in W} if $W-J$ is irreducible and has no component with 
compact closure. Note that in our situation the first condition is 
satisfied if and only if $J$ is not contained in a 3-ball in $W$; if, 
in addition,  
$W$ has one end and trivial second homology, then the second condition 
is satisfied if and only if $\bd J$ is connected. 

Let $J$ be a regular 3-manifold in $W$, and let $V$ be a connected open subset 
of $W$ which contains $J$. We say that $V$ is \textit{\ei\ rel $J$ in 
$W$} if for each compact, connected \tm\ $K$ with $J \sbs K \sbs V$ 
there is a compact, connected \tm\ $L$ with $K \sbs L \sbs V$ so that 
any loop in $V-L$ which is \nh\ in $W-J$ must be \nh\ in $V-K$. We say 
that $W$ is \textit{\ei\ rel $J$} if $W$ is \ei\ rel $J$ in $W$. We say 
that $W$ is \textit{\eei\ }if it is \ei\ rel $J$ for some $J$. (In \cite{BT} 
neither $K$ nor $L$ is required to be connected, but since $V$ is 
connected one can enlarge either of these sets to a connected 
set in $V$, so the two definitions are equivalent.)

It will be useful to have the following alternative characterization 
in terms of exhaustions. 
An \textit{exhaustion} of a connected open \tm\ $U$ is a sequence 
$\{C_n\}_{n\geq 0}$ of compact, connected \tm s in $U$ such that 
$C_n \sbs int \, C_{n+1}$ for each $n\geq 0$ and 
$U=\cup_{n=0}^{\infty} C_n$. It is a \textit{regular exhaustion} 
if each $C_n$ is regular in $U$.  

%Lemma 2.1
\begin{lem}[Brin-Thickstun] Let $J$ be a regular 3-manifold in $W$, 
and let $V$ be a connected open subset of $W$ which contains $J$. 
Then $V$ is \ei\ rel $J$ in $W$ if and only if $V$ has a regular 
exhaustion $\{C_n\}_{n\geq 0}$ such that $C_0=J$ and $\bd C_n$ is 
\inc\ in $W-J$ for $n \geq 1$. \end{lem}

\begin{proof} This is part of Lemma 2.2 of \cite{BT}. \end{proof}

Let $X$ and $Y$ be manifolds and $h:X \times [0,1] \rightarrow Y$ 
a map; we also use the notation $h_t:X\rightarrow Y$, where 
$h_t(x)=h(x,t)$, to describe this map. This is an \textit{isotopy} if 
for each $t$ we have that $h_t$ is an embedding. It is an 
\textit{ambient isotopy} if each $h_t$ is a homeomorphism. If $X$ and 
$X\p$ are subsets of $Y$, then they are \textit{isotopic} if there is 
an isotopy $h_t:X \rightarrow Y$ such that $h_0$ is the inclusion 
map of $X$ into $Y$ and $h_1(X)=X\p$; if $A$ is a subset of $X$ such 
that the restriction of each $h_t$ to $A$ is the identity of $A$, then 
we say that $X$ and $X\p$ are \textit{isotopic rel $A$}. If $h_t$ is the 
restriction of an ambient isotopy $g_t:Y \rightarrow Y$ with $g_0$ the 
identity of $Y$, then we say that they are \textit{ambient isotopic}; 
if each $g_t$ is the identity on $A$ we say that they are \textit{ambient 
isotopic rel $A$}. The closure of the set of all points $x$ such that 
for some $t$ one has $g_t(x) \neq x$ is called the \textit{support} of 
$g$. 

Note that if $X$ and $X\p$ are isotopic via the isotopy $h_t:X\ra Y$, 
then $X\p$ and $X$ are isotopic via the isotopy $p_t:X\p\ra Y$ given 
by $p_t(x\p)=h_{1-t}(h\n_1(x\p))$. 

Let $J$ be a regular 3-manifold in $W$, and let $V$ a connected open subset 
of $W$ containing $J$. We say that $V$ is an \textit{\er\ of $W$ at $J$} 
if $V-J$ is \irr, $W-V$ has no components with compact closure, 
$V$ is \ei\ rel $J$ in $W$, and whenever $N$ is a regular 3-manifold 
in $W$ such that $J\sbs \inte N$ and $\bd N$ is \inc\ in $W-J$ it must be 
the case that for some compact subset $C$ of $W-J$ one has that $V$ is  
ambient isotopic rel $W-C$ to a subset $V\p$ of $W$ such that 
$N \sbs V\p$. The last property is called the 
\textit{weak engulfing property}. 

%Theorem 2.2
\begin{thm}[Brin-Thickstun] Let $W$ be a connected, orientable, 
\irr, open \tm. Let $J$ be a regular \tm\ in $W$. Then an \er\ $V$ of 
$W$ at $J$ exists. If $V$ and $V\p$ are \er s of $W$ at $J$, then they 
are isotopic rel $J$. \end{thm}

\begin{proof} Existence and uniqueness are Theorems 2.1 and 2.3, 
respectively, of \cite{BT}. \end{proof}

We remark, as pointed out in \cite{BT}, that two end reductions 
of $W$ at $J$ need not be ambient isotopic. 

%Lemma 2.3
\begin{lem}[Brin-Thickstun] Let $V$ be an \er\ of $W$ at $J$. 
Suppose $J\p$ is a 
regular 3-manifold in $V$ which contains $J$ such that $\bd J\p$ 
is \inc\ in $W-J$. Then $V$ is an \er\ of $W$ at $J\p$. \end{lem}

\begin{proof} This is Corollary 2.2.1 of \cite{BT}. \end{proof}

Suppose \ka\ is a knot in $W$ which does not lie in a 3-ball in $W$. 
We say that $V$ is an \textit{\er\ of $W$ at \ka\ } if $V$ is an \er\ of $W$ at 
a regular neighborhood of \ka\ in $W$. 

%Lemma 2.4
\begin{lem} Let $V$ be an \er\ of $W$ at $J$. Then there is a knot \ka\ 
in $\inte J$ such that $V$ is an \er\ of $W$ at \ka. \end{lem}

\begin{proof} Since $J$ is regular in $W$ no component of $\bd J$ is 
a 2-sphere. Therefore by \cite{My homology} there is a knot \ka\ in 
$\inte J$ such that $\bd J$ is \inc\ in $J-\ka$. The result then 
follows from Lemma 2.3. \end{proof}

%The material reviewed so far is sufficient for an understanding of 
%the end reduction ratchet lemma and the main theorem, which are 
%presented in the next two sections. 
In dealing with the examples to be 
presented later we will need a means of recognizing 
an end reduction $V$ of a specific $W$ at a specific $J$. In their 
proof of existence Brin and Thickstun provide a procedure for 
constructing $V$, which we now briefly describe. 

Let $X$ and $Y$ be \tm s, with $X \sbs \inte Y$. 
Suppose $\bd X$ has a component which is compressible in 
$Y$ and is not a 2-sphere. Then there is a disk 
$D$ in $Y$ such that $D\cap \bd X$ is a non \nh\ simple closed curve in 
$\bd X$. 
We can construct a new \tm\ $X\p$ by doing surgery on $\bd X$ 
along $D$. If $D\sbs Y-\inte X$, then we obtain $X\p$ by adding a 
2-handle $H$ to $X$. If $D \sbs X$ then we obtain $Y-\inte X\p$ by adding a 
2-handle $H$ to $Y-\inte X$, i.e. by removing a 1-handle from $X$. In the two 
cases the handle $H$ is a regular neighborhood $D \times [-1,1]$ of $D$ 
in $Y-\inte X$ or $X$ respectively. We say that a \tm\ $X^*$ in $Y$ is 
\textit{obtained from $X$ by completely compressing $\bd X$ in $Y$} if 
there is a sequence $X=X_0, X_1, \ldots, X_k=X^*$ of 3-manifolds 
in $Y$ such 
that $X_{i+1}$ is obtained by adding or removing a handle $H_i$ along 
a compressing disk $D_i$ for $\bd X_i$ in $Y$, the $D_i$ are pairwise 
disjoint, the attaching/removing annuli $\bd D_i \times [-1,1]$ 
are pairwise disjoint, and each component of $\bd X^*$ which is 
not a 2-sphere is incompressible in $Y$. If, 
in addition, there is a 3-manifold $A$ in $Y$ which contains $\bd X$ and 
all the $H_i$, then we say that the compressions are \textit{confined 
to $A$}. 

For future reference we note that if $X^*$ is obtained from $X$ by 
completely compressing $\bd X$ in $Y$ with the compressions confined 
to $A$, then by general position one may arrange for the following 
property to hold. Suppose $i<j$. Then $H_i$ and $H_j$ are either disjoint 
or meet in a finite union of disjoint 3-balls each of which has the form 
$D\times [a,b]$, $-1<a<b<1$, when considered as a subset of $H_i$ and the 
form $E\times [-1,1]$, $E$ a disk in $\inte D$, when considered as a 
subset of $H_j$. 

After a 2-handle is added it may subsequently have ``sub-handles'' 
removed by the removal of 1-handles and added by addition of 2-handles, 
etc.\ , so that the net result within the 2-handle is the addition of 
several smaller parallel copies of itself.  Similarly after a 1-handle 
has been removed it may subsequently have sub-handles added by the addition of 
2-handles and sub-handles removed by the removal of 1-handles, and so 
on, so that the net result within the 1-handle is the removal of several 
smaller parallel copies of itself.    

Note also that since the removing/attaching annulus for $H_j$ is 
not null homotopic in $\bd X_j$ it cannot lie in the interior of any 
of the disks $D\times\{\pm1\}$ belonging to any of the handles $H_i$ 
with $i<j$. Since the entire set of annuli is pairwise disjoint   
this implies that it must lie in the original $\bd X$.

Now suppose that $J$ is a regular \tm\ in $W$. It follows from our 
assumptions on $W$ that it has a regular exhaustion $\{C_n\}$ with 
$C_0=J$. By passing to a subsequence we may assume that each $\bd C_n$,  
for $n\geq 1$,  
can be completely compressed in $W-C_0$ with the compressions confined 
to $C_{n+1}-C_0$. We obtain $C_1^*$ from $C_1$ by completely compressing 
$\bd C_1$ in $W-C_0$ so that $C_1^* \sbs C_2$. Now $\bd C_2$ can be 
completely compressed in $W-C_0$ with the compressions confined to 
$C_3-C_0$. The components of $\bd C_1^*$ which are not 2-spheres  
are \inc\ in $W-C_0$ and are therefore \inc\ in the 
smaller set $C_3-C_0$. Any compressing disk for $\bd C_2$ can be 
isotoped in $C_3-C_0$ so as to lie in $C_3-C_1^*$. The same is true 
for the compressing disks for the intermediate surfaces obtained in 
the process of completely compressing $\bd C_2$ in $W-C_0$. It follows 
that $\bd C_2$ can be completely compressed in $W-C_1^*$ with the 
compressions confined to $C_3-C_1^*$, so that we obtain $C_1^*\sbs 
C_2^*\sbs C_3$. We continue in this fashion, and set $C_0^*=C_0$ to 
get a sequence $\{C_n^*\}$. 

%Note that $C_n^*$ may have several components. 
%Let $Z$ be such a component. In our context there is a component $\bd^+Z$ 
%of $\bd Z$ which bounds a compact \tm\ $\widehat{Z}$ in $W$ which contains 
%$Z$ and is disjoint from any other components of $C_n^*$. Let $\bd^-Z=
%\bd Z-\bd^+Z$. We let $\widehat{C}_n^*$ be the union of the $\widehat{Z}$. 
%In a similar way we define $\bd^+C_n^*$ and $\bd^-C_n^*$. 

The union $V^*$ of the sequence $\{C_n^*\}$ is an 
open subset of $W$ which is called the \textit{constructed end 
reduction of $W$ at $J$}. 

%Theorem 2.5
\begin{thm}[Brin-Thickstun] The component $V$ of $V^*$ containing 
$J$ is an end reduction of $W$ at $J$. \end{thm}

\begin{proof} The proof is that of Theorem 2.1 of \cite{BT}. \end{proof}
  
Note that $C_n^*$ may have several components and that $W-C_n^*$ may 
have components with compact closure. Suppose $Z$ is a compact, connected 
\tm\ in $W$ of one of these two types. Since the handles used to compress 
$\bd C_{n+1}$ miss $C_n^*$ they also miss $Z$. It follows that $Z\sbs 
\inte C_{n+1}^*$ and that either $Z$ is in $V$ or $Z$ is in $W-V$. 
By choosing an appropriate union of those $Z$ which are in $V$ we can 
produce a regular \tm\ $C\p_n$ in $W$ so that $\{C\p_n\}$ is an 
exhaustion of $V$ satisfying Lemma 2.1.

We recall that a \textit{Whitehead manifold} is an \irr, contractible 
open \tm\ which is not \hm\ to \RRR. 

%Lemma 2.6
\begin{lem} Every \er\ of a Whitehead manifold is a Whitehead manifold. 
\end{lem}

\begin{proof} Let $W$ be a Whitehead manifold and $V$ an end reduction 
of $W$ at a regular \tm\ $J$. 

We first show that $V$ is simply connected. We may assume that $V$ is 
the component of a constructed end reduction $V^*$ of $W$ at $J$ 
which contains $J$. Let $\{C_n\}_{n\geq 0}$ be the exhaustion of $W$ 
used in the construction of $V^*$. Let $\al$ be a loop in $V$. Then 
$\al \sbs C^*_p$ for some $p\geq 1$, where $C^*_p$ is obtained from 
$C_p$ by completely compressing $\bd C_p$ in $W-J$, where the 
compressions are confined to $C_{p+1}-C_0$. Recall that 
there is a sequence $X_0, \ldots, X_k$ such that $X_0=C_p$, $X_k=C^*_p$, 
and $X_{i+1}$ is obtained from $X_i$ by either adding a 2-handle or 
removing a 1-handle. 

We first claim that \al\ can be homotoped in $C^*_p$ so that it lies 
in $C_p\cap C^*_p$. Suppose, by induction, that we have homotoped \al\ 
in $C^*_p$ so that it lies in $X_{i+1}\cap C^*_p$. If $X_{i+1}$ is obtained 
from $X_i$ by removing a 1-handle, then \al\ automatically lies in 
$X_i\cap C^*_p$. So assume that $X_{i+1}$ is obtained from $X_i$ by 
adding a 2-handle $H_i$. As noted above $H_i\cap C^*_p$ consists of 
either $H_i$ itself or several parallel copies of itself. In any case 
\al\ can be homotoped out of each such component by pushing it into $X_i$ 
across the intersection of that component with the attaching annulus for 
$H_i$. This homotopy takes place in $C^*_p$.    

Since $W$ is simply connected there is a map $f:B \ra W$, where $B$ is 
a disk and $f(\bd B)=\al$. There is a $q > p$ such that $f(B) \sbs C_q$. 
$C^*_q$ is obtained from $C_q$ by completely compressing $\bd C_q$ 
in $W-C^*_{q-1}$ with the compressions confined to $C_{q+1}-C^*_{q-1} 
\sbs C_{q+1}-C^*_p$. There is a sequence $Y_0,\ldots, Y_m$ such that 
$Y_0=C_q$, $Y_m=C^*_q$, and $Y_{i+1}$ is obtained from $Y_i$ by either 
adding a 2-handle or removing a 1-handle. Note that $C^*_p\sbs C^*_q$. 

We claim that there is a map $f_m:B\ra C^*_q$ such that $f_m$ agrees with 
$f$ on $\bd B$.  Let $f_0=f$. Suppose, by induction, that we have an 
$f_i:B\ra Y_i$ 
which agrees with $f$ on $\bd B$. If $Y_{i+1}$ is obtained by attaching 
a 2-handle to $Y_i$, then $f_i(B)$ is automatically in $Y_{i+1}$, so we 
let $f_{i+1}=f_i$. So assume 
that $Y_{i+1}$ is obtained by removing a 1-handle $H_i$ from $Y_i$. 
Then $H_i=D\times[-1,1]$ for a disk $D$. Since $H_i\cap C^*_p=\ns$ we may 
assume that $f_i\n(D\times\{0\})$ consists of simple closed curves in 
$\inte B$. We may then redefine $f_i$ in a neighbhorhood of the outermost 
disks bounded by these curves to obtain the desired map $f_{i+1}$.    
This finishes the proof of simple connectivity. 

We next show that $V$ is \irr. Let $S$ be a 2-sphere in $V$. Since 
$W$ is \irr, $S$ bounds a 3-ball $B$ in $W$. Since $W-V$ has no 
components with compact closure we must have that $B\sbs V$. 

Since $V$ is orientable it then follows from the sphere theorem that 
$\pi_2(V)=0$. Since $V$ is non-compact $H_3(V)=0$. It then follows 
from the Hurewicz theorem that $\pi_i(V)=0$ for all $i\geq 3$. By 
the Whitehead theorem $V$ is contractible. 

Since $J$ is a regular \tm\ in $W$ it cannot lie in a 3-ball in $W$, 
hence cannot lie in a 3-ball in $V$. Therefore $V$ is not \hm\ to \RRR. 
\end{proof}

%Section 3
\section{The End Reduction Ratchet Lemma}

In this section we prove a technical result, called the End Reduction 
Ratchet Lemma, which generalizes earlier work of Wright \cite{Wr} 
and of Tinsley and Wright \cite{TW}. The Ratchet Lemma of \cite{Wr} is  
basically the special case $W=V$ of our result. The Special Ratchet 
Lemma of \cite{TW} is the special case where $V$ is an \eei\ open 
subset of $W$ bounded by proper planes. 

%Theorem 3.1
\begin{thm}[End Reduction Ratchet Lemma] Let $W$ be a connected, 
orientable, \irr\ open \tm\ which is not \hm\ to \RRR. Let $J$ be a regular 
\tm\ in $W$. Suppose $V$ is an end reduction of $W$ at $J$. Let $g$ 
be a homeomorphism of $W$ onto itself such that each of $g(J)$ and 
$g^{-1}(J)$ can be isotoped into $V$. Then there is a compact \tm\ 
$R$ in $W$ containing $J$ such that a loop in 
$W-\cup_{i=-\infty}^{\infty}g^i(R)$ is \nh\ in $W-J$ if and only if 
it is \nh\ in $W-g^i(J)$ for each $i \in \mathbf{Z}$. \end{thm}

\begin{proof} $g(J)$ and $g^{-1}(J)$ are isotopic to subsets $J^+$ and 
$J^-$ of $V$, respectively. Since $J$ is compact 
the covering isotopy theorem \cite{Ce,EK} implies that there are 
ambient isotopies $k^{\pm}_t:W\ra W$ with $k^{\pm}_0$ the identity of $W$  
and $k^{\pm}_1(g^{\pm 1}(J))=J^{\pm}$; moreover there are compact 
sets $T^{\pm}$ such that $k^{\pm}_t(g^{\pm 1}(J))\sbs T^{\pm}$ for all 
$t\in [0,1]$ and $k_t(x)=x$ for all $x\in W-T^{\pm}$ and $t\in [0,1]$. 
Let $K$ be a compact, connected \tm\ 
in $V$ which contains $J \cup J^+ \cup J^-$. Since $V$ is \ei\ rel $J$ 
in $W$ there is a compact, connected \tm\ $L$ in $V$ such that $K\sbs L$ 
and  any loop in $V-L$ which is \nh\ in $W-J$ is \nh\ in $V-K$. Let $R$ be 
a compact, connected \tm\ in $W$ which contains $T^+ \cup T^- \cup L$. 

Suppose $\ga$ is a loop in $W-\cup_{i=-\infty}^{\infty}g^i(R)$ which 
is \nh\ in $W-J$. Let $D$ be a disk and $f:D\ra W-J$ a map which is 
transverse to $\bd L$ with $f(\bd D)=\ga$. There is a collection 
(possibly empty) of disjoint disks $E_1, \ldots, E_p$ in $\inte D$ 
whose boundaries are in $f^{-1}(\bd L)$ and whose interiors contain all  
other points of $f^{-1}(L)$. Each loop $f(\bd E_j)$ lies in $\bd L$ 
and is \nh\ in $W-J$; it follows that it is \nh\ in $V-K$ and hence is 
\nh\ in $V-J^+$. We may therefore modify $f$ so that $f(D)\sbs W-J^+$. 
Since $\ga \sbs W-R \sbs W-T^+$, we have that $\ga=(k_1^+)^{-1}(\ga)$ 
and so bounds the singular disk $(k^+_1)^{-1}(f(D))$ in $W-g(J)$. 
Thus $\ga$ is \nh\ in $W-g(J)$. 

A similar argument shows that a loop in 
$W-\cup_{i=-\infty}^{\infty}g^i(R)$ which is \nh\ in $W-J$ must be 
\nh\ in $W-g^{-1}(J)$. Translation by appropriate powers of $g$ then 
completes the proof. \end{proof}

%Section 4
\section{The main theorem}

%Theorem 4.1
\begin{thm} Let $W$ be an \irr\ contractible open \tm\ which is not 
\hm\ to \RRR. Suppose $W$ is a non-trivial covering space of a \tm\ 
$M$. Let $G\cong \pi_1(M)$ be the group of covering translations. 
Let $J$ be a regular \tm\ in $W$, and let $V$ be an \er\ of $W$ at $J$. 
If $g$ is an element of $G$ such that $g(J)$ and $g\n(J)$ 
can each be isotoped into $V$, then $g$ is the identity of $W$. \end{thm}

A group $G$ acting on a set $X$ acts \textit{without fixed points} if 
the only element of $G$ which fixes a point is the identity. 

%Corollary 4.2
\begin{cor} $G$ acts without fixed points on the set of isotopy classes 
of end reductions of $W$ and on the set of ambient isotopy classes of 
end reductions of $W$. \end{cor}

\begin{proof}[Proof of Corollary 4.2] Suppose $g\in G$ fixes the isotopy 
class of the end reduction $V$ of $W$. Thus there is an isotopy 
$h_t:V\ra W$ such that $h_0$ is the inclusion map of $V$ into $W$ and 
$h_1(V)=g(V)$. Let $J$ be a regular \tm\ in $W$ at which $V$ is the end 
reduction of $W$. We will construct isotopies 
$f^{\pm}_t:g^{\pm1}(J)\ra W$ such that 
$f^{\pm}_0$ is the inclusion map of $g^{\pm1}(J)$ into $W$ and 
$f^{\pm}_1(g^{\pm}(J))\sbs V$. 

We define $f^-_t:g\n(J)\ra W$ by $f^-_t(x)=g\n(h_t(g(x)))$. 
If $x\in g\n(J)$, then $g(x)\in J\sbs V$, so we have 
$f^-_0(x)=g\n(h_0(g(x)))=g\n(g(x))=x$. On the other hand 
$f^-_1(g\n(J))=g\n(h_1(g(g\n(J))))=g\n(h_1(J))\sbs g\n(h_1(V))=
g\n(g(V))=V$. 

Recall that $g(V)$ is isotopic to $V$ via the isotopy 
$p_t(x)=h_{1-t}(h_1\n(x))$. Thus $p_0$ is the inclusion map of 
$g(V)$ into $W$, and $p_1(g(V))=V$. Let $f^+_t$ be the restriction 
of $p_t$ to $g(J)$. If $x\in g(J)$, then $f^+_0(x)=p_0(x)=x$. On the 
other hand $f^+_1(g(J))=p^+_1(g(J))\sbs p^+_1(g(V))=V$. 

The case of an ambient isotopy clearly follows from 
the case just done. \end{proof}

\begin{proof}[Proof of Theorem 4.1] By Lemma 2.5 $V$ is a Whitehead 
manifold. We will show that if $g$ is not the identity, then $V$ must 
be \hm\ to \RRR, thereby contradicting this statement. 
The argument follows the tradition 
of \cite{Wr}, \cite{TW}, and \cite{My free}, in combining the Orbit Lemma 
of Wright with an appropriate generalization of his Ratchet Lemma, in 
this case the End Reduction Ratchet Lemma. 

A group $G$ acting on an  $n$-manifold $W$ acts \textit{totally 
discontinuously} if for every compact subset $C$ of $W$ one has that 
$g(C)\cap C = \ns$ for all but finitely many elements of $G$. 
$G$ acts without fixed points 
and totally discontinuously on $W$ if and only if the projection to 
the orbit space of the action is a regular covering map with group 
of covering translations $G$ and the orbit space is an $n$-manifold 
\cite{Ma}. In this case if $W$ is contractible, then $G$ is torsion 
free. 

%Proposition 4.3
\begin{prop}[Orbit Lemma (Wright)] Let $W$ be a con\-tract\-i\-ble, open 
$n$-man\-i\-fold, $n \geq 3$. Let $g$ be a non-trivial homeomorphism of $W$ 
onto itself such that the group $<g>$ of homeomorphisms generated by 
$g$ acts without fixed points and totally discontinuously on $W$. 
Given compact subsets $B$ and $Q$ of $W$, there is a compact subset 
$C$ of $W$ containing $B$ such that every loop in $W-C$ is homotopic in $W-B$ 
to a loop in $W-\cup_{i=-\infty}^{\infty} g^i(Q)$. \end{prop}

\begin{proof} Except for the statement that $C$ contains 
$B$ this is Lemma 4.1 of \cite{Wr}; we can clearly enlarge the $C$ of that 
result to satisfy this requirement. A somewhat shorter alternative 
proof for the special case in which $W$ is an \irr\ contractible 
open \tm\ is given in \cite{My free}. \end{proof} 

Returning to the proof of Theorem 4.1, we assume that $g$ is not the 
identity. It then satisfies the hypotheses of the Orbit Lemma. 
We shall prove that $V$ is 
\textit{$\pi_1$-trivial at $\infty$}, i.e. for every compact subset 
$A$ of $V$ there is a compact subset $A\p$ of $V$ containing 
$A$ such that every loop in $V-A\p$ is null-homotopic in $V-A$. 
By a result of C. H. Edwards \cite{Ed} and 
C. T. C. Wall \cite{Wl} every irreducible, contractible, open 
3-manifold which is $\pi_1$-trivial at $\infty$ must be homeomorphic 
to $\mathbf{R}^3$.  

We may assume that $V$ is the component of the constructed end reduction 
$V^*$ containing $J$. Thus we have a regular exhaustion $\{C_n\}$ of $W$ 
with $C_0=J$ from which we obtain $\{C_n^*\}$ with union $V^*$ and 
$\{C\p_n\}$ with union $V$. 

Let $A$ be a compact subset of $V$. Let $K$ be a compact \tm\ in 
$V$ which contains $J \cup A$. 
Since $V$ is \ei\ rel $J$ in $W$, there 
is a compact \tm\ $L$ in $V$ such that every loop in $V-L$ which is 
\nh\ in $W-J$ is \nh\ in $V-K$. 
%By Lemma 2.1 there is a regular 
%exhaustion $\{C_n\}_{n\geq 0}$ of $V$ such that $C_0=J$ and 
%$\bd C_n$ is \inc\ in $W-J$ for $n \geq 1$. 
By replacing $L$ by a 
$C\p_n$ containing $L$ we may assume that $\bd L$ is \inc\ in $W-J$. 
Then by Lemma 2.3 $V$ is an end reduction of $W$ at $L$. 

By the End Reduction 
Ratchet Lemma there is a compact \tm\ $R$ in $W$ containing $J$ such 
that a loop in $W-\cup_{i=-\infty}^{\infty} g^i(R)$ is null-homotopic in 
$W-J$ if and only if it is null-homotopic in $W-g^i(J)$ for all $i \in 
\mathbf{Z}$. 
Apply the Orbit Lemma with $B=L$ and $Q=R$ to get a compact subset $C$ 
of $W$ containing $L$ such that every loop in $W-C$ is homotopic in $W-L$ 
to a loop in $W-\cup_{i=-\infty}^{\infty} g^i(R)$. By enlarging $C$ 
we may assume that it is some $C_p$, $p>n$.   

Now $C_p^*$ is the \tm\ obtained by completely compressing $\bd C_p$ 
in $W-C^*_n$. Let $A\p=C^*_p\cap V$.  
We must show that every loop $\gamma$ in $V-A\p$ is 
null-homotopic in $V-A\p$. 

We first show that $\gamma$ is homotopic in $W-L$ to a loop 
$\gamma\p$ in $W-C_p$. 
Since $\gamma$ is in $V$ it cannot meet any component of 
$C^*_p$ which is not in $V$. Therefore $\gamma \cap C_p$  must 
lie in the union of the 1-handles which were removed in 
the process of completely compressing $\bd C_p$. 
These handles miss $C_n^*$ and hence miss $C_n\p=L$, so 
we can homotop $\gamma$ in $W-L$ to a $\gamma\p$ in $W-C_p$.

Now $\gamma^{\prime}$ is 
homotopic in $W-L$ to a loop $\gamma^{\prime\prime}$ in 
$W-\cup_{i=-\infty}^{\infty} g^i(R)$. Since $W$ is contractible 
$\gamma^{\prime\prime}$ is null-homotopic in $W$. Since $<g>$ is totally 
discontinuous $\gamma^{\prime\prime}$ is null-homotopic in 
$W-g^i(J)$ for some $i$. Since $\gamma^{\prime\prime}$ lies in $W-L$ 
the End Reduction Ratchet Lemma implies that $\gamma^{\prime\prime}$ is 
null-homotopic in $W-J$. Since $J \subseteq L \subseteq A\p$ we have 
that $\gamma$ is null-homotopic in $W-J$. Since $\gamma$ lies in $V-L$ 
and $V$ is \ei\ rel $J$ in $W$ we have that 
$\gamma$ is null-homotopic in $V-J$. Thus $\gamma$ is null-homotopic 
in $V-K \subseteq V-A$, as required. \end{proof}

%We now completely compress $\bd C$ in $W-L$ to obtain a \tm\ $C^*$ 
%whose boundary components are 2-spheres and surfaces which are \inc\ 
%in $W-L$. We let $C\p$ be the component of $C^*$ containing $L$. 
%Since $W$ and $W-L$ are \irr\ and $L$ is regular, we may assume that 
%$\bd C\p$ is \inc\ in $W-L$ and has no 2-sphere components. 

%By the weak engulfing property of $V$ there is an ambient isotopy $h_t$ 
%of $V$ rel $L$ to an open subset $V\p$ of $W$ which contains $C\p$. 
%Note that $A\sbs L \sbs C\p$; thus $A\sbs L \sbs h_1^{-1}(C\p)$. 
%We claim that we may take $A^*=h_1^{-1}(C\p)$. 

%It suffices to show that every loop $\gamma$ in $V\p-C\p$ is  
%null-homotopic in $V\p-A$. 
%First note that $\gamma \cap C$ is contained in the union of the 
%1-handles which were removed from $C$ in the process of obtaining $C^*$ 
%and possibly some of the 3-balls which were added to the boundary 
%of the original $C^+$. We first homotop $\gamma$ out of any such 
%3-balls so that the intersection with $C$ lies in the original $C\p$. 
%We then homotop $\gamma$ in $W-L$ to a loop $\gamma\p$ which misses 
%the 1-handles that were removed from $C$. So $\gamma\p$ lies in $W-C$.   

%Section 5
\section{Nested end reductions}

%Theorem 5.1
\begin{thm} Let $W$ be a Whitehead manifold. Let $J$ and $K$ be 
regular 3-manifolds in $W$ such that $J\sbs \inte K$. Let $U$ and 
$V$ be end reductions of $W$ at $J$ and $K$, respectively. Then $U$ is  
isotopic rel $J$ to $U\p$ such that $U\p\sbs V$. Moreover $U\p$ is an 
end reduction of $V$ at $J$. \end{thm}

\begin{proof} Let $\{C_n\}_{n\geq0}$ be a regular exhaustion for $W$ 
with $C_0=J$ and $C_1=K$. We may assume that $V$ is the component of the 
constructed end reduction $V^*$ of $W$ at $K$ associated to the exhaustion 
$\{C_n\}_{n\geq1}$ which contains $K$. Let $C_n^*$ be the \tm\ arising 
in the construction of $V^*$. Recall that $C_n^*$ is obtained from $C_n$ by 
completely compressing $\bd C_n$ in $W-K$ with the 
compressions confined to $C_{n+1}-C_{n-1}^*$. 

It suffices to show that after passing to a subsequence 
the constructed end reduction $U^{\#}$ of $W$ at $J$ associated 
to the exhaustion $\{C_n\}_{n\geq0}$ can be built in such a 
way that the \tm\ $C_n^{\#}$ obtained by completely 
compressing $C_n$ in $W-J$ with the compressions confined to 
$C_{n+1}-C_{n-1}^{\#}$ can be obtained by first constructing $C_n^{\#}$ as 
above and then completely compressing $\bd C_n^*$ in $W-J$ with 
the compressions confined to $C_{n+1}^*-C_{n-1}^{\#}$. 

We first completely compress $\bd C_1$ in $W-J$. Let $E_1,\ldots, E_p$ be 
the sequence of compressing disks and $H_1, \ldots, H_p$ the associated 
sequence of 1-handles and 2-handles. We require the handles to intersect 
as described in Section 2. Let $C_1=Q_0, \ldots, Q_p$ be the associated 
sequence of 3-manifolds. We may assume that the compressions 
are confined to $C_2-J$. 

Let $D_1, \ldots, D_q$ be the sequence of compressing disks 
involved in the complete compression of $\bd C_2$ in $W-K$. 
Let $C_2=X_0, X_1, \ldots, X_q=C_2^*$ be the associated sequence 
of 3-manifolds. We may assume that the compressions are 
confined to $C_3-K$. 

We will show that by replacing $E_1, \ldots, E_p$ by a new set of disks 
we can completely compress $\bd C_1$ in $W-J$ with the 
compressions confined to $C_2^*-J$. The resulting 3-manifold 
will be our new $C_1^{\#}$. This will be done in $q$ steps. After 
each step we will rename our new disks and 3-manifolds 
with the names of their predecessors. 

Assume as the inductive hypothesis 
that $\bd C_1$ can be completely compressed in $W-J$ 
with the compressions confined to $X_{i-1}-J$. 
If $X_i$ is obtained by adding a 2-handle to $X_{i-1}$, 
then there is nothing to do. So assume that it is 
obtained by removing a 1-handle from $X_{i-1}$ associated to 
the compressing disk $D_i$ for $\bd C_2$ in $C_2-K$. 
Let $E=E_1\cup\cdots\cup E_p$. Put $E$ in general postion 
with respect to $D_i$. 

Let $\{\al_1,\ldots, \al_m\}$ be the set of components of $E\cap D_i$ 
which are outermost on $E$. So $\al_s=\bd\Delta_s$ for a 
disk $\Delta_s$ in $E$, the $\Delta_s$ are disjoint, and this union 
contains $E\cap D_i$. Let $F=E-(\Delta_1\cup\ldots\cup\Delta_m)$. 
It is a disk with $m$ holes. 
Now $\al_s=\bd \Delta_s\p$ for a disk $\Delta_s\p$ in $D_i$. 
The $\Delta_s\p$ need not be disjoint. If two $\Delta_s\p$ 
intersect, then one is contained in the interior of the other. 

We may assume that $\Delta_1\p$ is innermost on $D_i$ among the $\Delta_s\p$. 
Let $F_1=F\cup\Delta_1\p$. 
Isotop $F_1$ by moving $\Delta_1\p$ to a \pl\ disk which misses $D_i$ 
in such a way that now $F_1\cap D_i=\al_2\cup\cdots\cup\al_m$. 
We may now assume that $\Delta_2\p$ is innermost on $D_i$ among the 
remaining $\Delta_s\p$. Let $F_2=F_1\cup\Delta_2\p$. Isotop $F_2$ by moving 
$\Delta_2\p$ to a \pl\ disk which misses $D_i$ in such a way that 
now $F_2\cap D_i=\al_3\cup\dots\cup\al_m$. We continue in this fashion 
until we get a surface $E\p=F_m$ which is a union of $q$ disks 
which are disjoint from each other and from $D_i$. Let 
$E_j\p$ be the component of $E\p$ with $\bd E_j\p=\bd E_j$. 

Let $A_j$ be the removing/attaching annulus for the handle $H_j$. 
Let $G_j=\bd C_1-\inte(A_1\cup\cdots\cup A_j)$. By construction 
$E_j\p\cap G_{j-1}=\bd E_j\p=\bd E_j$, which is a core of $A_j$. 
Since $E_j\p\cap E_k\p=\ns$ for $k<j$ and $\bd Q_{j-1}\p$ is the 
union of $G_j$ and two \pl\ copies of each $E_k\p$ with 
$k<j$ we have that $E_j\p\cap\bd Q_{j-1}\p=\bd E_j\p$. It follows 
that $Q_j\p$ is obtained from $Q_{j-1}\p$ by removing a 
1-handle or adding a 2-handle along $E_j\p$. Since 
$\bd Q_{j-1}\p$ and $\bd Q_{j-1}$ are homeomorphic by a homeomorphism 
which is the identity on $G_{j-1}$ and $\bd E_j$ does not bound 
a disk in $\bd Q_{j-1}$ we have that $\bd E_j\p$ does not bound a 
disk in $\bd Q_{j-1}\p$. Thus $E_j\p$ is a compressing disk. 

Now suppose that $D\p$ is a compressing disk for $\bd Q_p\p$ 
in $W-J$. Since $\bd Q_p\p-\inte G_p$ consists of disks 
we can isotop $D\p$ so that $\bd D\p$ is in $G_p$. If $D\p\cap \bd Q_p 
\neq \bd D\p$, then the excess intersections must lie in the interior 
of the disks comprising $\bd Q_p-\inte G$. 
We do surgery on $D\p$ along these intersections to get 
a disk $D$ with $\bd D=\bd D\p$ and $D\cap \bd Q_p=\bd D$. 
Since $Q_p=C_1^{\#}$ and $\bd C_1^{\#}$ is \inc\ in $W-J$ 
we have that $\bd D\p$ bounds a disk in $\bd Q_p$ and 
hence bounds a disk in $\bd Q_p\p$. 

Thus $Q_p\p$ has 
been obtained from $C_1$ by completely compressing 
$\bd C_1$ in $W$ with the compressions confined to $X_i-J$. 
We then rename $Q_p\p$ as $C_1^{\#}$ and continue the 
induction on $i$ to finally get a 
$C_1^{\#}$ that is obtained from $C_1$ by completely 
compressing $\bd C_1$ in $W-J$ with the compressions 
confined to $C_2^*-J$. Note that $\bd C_1$ has also 
been completely compressed in $U-J$. 

We next completely compress $\bd C_2^*$ in $W-J$ to obtain 
$C_2^{\#}$. We may assume that the compressions are confined 
to $C_3-C_1^{\#}$. Recycling our notation we let $E_1,\ldots, E_p$ 
be the sequence of compressing disks and $H_1,\ldots, H_p$ the 
sequence of handles. $C_2^*=Q_0, Q_1, \ldots,Q_p=C_2^{\#}$ is 
the sequence of 3-manifolds. 

We now let $D_1,\ldots,D_q$ be the sequence of compressing 
disks involved in the complete compression of $\bd C_3$ in 
$W-K$ and let $C_3=X_0, X_1, \ldots, X_q=C_3^*$ be the 
associated sequence of 3-manifolds. We may assume that 
the compressions are confined to $C_4-C_2^*$. 

As in the previous step we replace the $E_j$ by disks 
$E_j\p$ with $\bd E_j\p=\bd E_j$by doing surgery on the $E_j$ along 
their intersections with the $D_i$. This again results 
in the complete compression of $\bd C_2^*$ in $W-J$ 
with the compressions confined to $C_3^*-C_1^{\#}$. 
It follows that we have completely compressed 
$\bd C_2$ in $W-J$ and $\bd C_2^*$ in $U-J$. 

We then continue this process to complete the proof. \end{proof}

%Section 6
\section{End reductions and planes}

Let $W$ be a Whitehead manifold. An embedded proper plane $\Pi$ in $W$ 
is \textit{trivial} if some component of $W-\Pi$ has closure \hs. 
$W$ is \textit{\rirr}\ if every proper plane in $W$ is trivial. 

%Theorem 6.1
\begin{thm} Let $W$ be a Whitehead manifold. Suppose $\Pi$ is a non-trivial 
plane in $W$ and $J$ is a regular submanifold of $W$ such that 
$J\cap\Pi=\ns$. Then any open subset $V$ of $W$ which is end irreducible 
rel $J$ in $W$ can be isotoped so that $V\cap \Pi=\ns$. In particular, 
any end reduction of $W$ at $J$ can be isotoped off $\Pi$. \end{thm}

\begin{proof} By Lemma 2.1 there is a regular exhaustion $\{C_n\}$ of $V$ 
such that $C_0=J$ and $\bd C_n$ is \inc\ in $W-C_0$ for all $n>0$. 
Let $W_n=W-\inte C_n$. 

We first isotop $\bd C_1$ in $W_0$ so that it is in minimal general position 
with respect to $\Pi$. This is to be done via an ambient isotopy which is 
fixed on $\bd W_0$. Since $\bd C_1$ is \inc\ in $W_0$ we must now have that 
$\Pi\cap\bd C_1=\ns$. Note that our isotopy may have moved all the $C_n$ and 
$W_n$ for $n>0$. 

We next isotop our new $\bd C_2$ in our new $W_1$ so that it is in 
minimal general position with respect to $\Pi$. This is to be done by 
an ambient isotopy which is fixed on $\bd W_1$. Since $\bd C_2$ is \inc\ 
in $W_1$ we must now have that $\Pi\cap\bd C_2=\ns$. Again our isotopy 
has possibly moved all the $C_n$ and $W_n$ for $n>1$. 

We continue in this fashion to move each $\bd C_n$ in turn off of $\Pi$. 
The points of each $C_n$ are moved a finite number of times and thereafter 
remain fixed. Thus the restriction of this series of ambient isotopies 
to $V$ converges to an isotopy of $V$ in $W$ which moves $V$ off $\Pi$. 
(Note that since we have established no control over the points of $W-V$ 
we cannot guarantee that this isotopy is ambient.) \end{proof}

%Theorem 6.2
\begin{thm} Let $W$ be a Whitehead manifold. Suppose $\Pi$ is a non-trivial 
plane in $W$. Then any \rirr\ open subset $V$ of $W$ which is end \irr\ rel 
$J$ in $W$ for some regular submanifold $J$ of $W$ can be isotoped so that 
$V\cap \Pi=\ns$. In particular, any \rirr\ \er\ of $W$ can be 
isotoped off $\Pi$. \end{thm}

\begin{proof} By Lemma 2.4 $V$ is an \er\ of $W$ at a knot $\kappa$. 
Put $\kappa$ in \mgp\ with respect to $\Pi$. If $\kappa\cap\Pi=\ns$, then 
we are done by Theorem 5.1. So assume that $\kappa\cap\Pi$ has $k>0$ 
components. We will show that $V$ contains a non-trivial plane. 

By Lemma 2.1 $V$ has a regular exhaustion $\{C_n\}$ such that $C_0$ is 
a regular neighborhood of $\kappa$ and $\bd C_n$ is \inc\ in $W-C_0$ 
for all $n>0$. Let $W_n=W-\inte C_n$. 

$\Pi\cap C_0$ consists of $k$ disks. Let $\Pi_0=\Pi\cap W_0$. Isotop $\Pi_0$ 
in $W_0$ rel $\bd \Pi_0$ so that $\Pi_0$ is in \mgp\ with respect to 
$\bd C_1$. Then no component of $\Pi_0\cap\bd C_1$ bounds a disk in $\Pi_0$. 
Every component of $\Pi_0\cap\bd C_1$ which is innermost on $\Pi$ bounds 
a disk in $\Pi$ which lies in $C_1$ and contains at least one component of 
$\Pi\cap C_0$. Thus there are at most $k$ such components. 

Let $\Pi_1=\Pi_0\cap W_1$. Isotop $\Pi_1$ in $W_1$ rel $\bd\Pi_1$ so that 
$\Pi_1$ is in \mgp\ with respect to $\bd C_2$. Then no component of 
$\Pi_1\cap\bd C_2$ bounds a disk in $\Pi_1$. Every component of $\Pi_1
\cap\bd C_2$ which is innermost on $\Pi$ bounds a disk in $\Pi$ which 
lies in $C_2$ and contains at least one of the disks bounded by an 
innermost component of $\Pi_0\cap\bd C_1$. Thus there are at most the 
same number of innermost components. 

We continue in this fashion by suitably isotoping 
$\Pi_{n+1}=\Pi_n\cap W_{n+1}$. There is an $n_0$ such that for all 
$n\geq n_0$ we have that $\Pi_n\cap\bd C_{n+1}$ and 
$\Pi_{n+1}\cap \bd C_{n+2}$ have the same number of components which are 
innermost on $\Pi$. Then the disk bounded by an innermost component of 
$\Pi_{n+1}\cap \bd C_{n+2}$ contains exactly one disk bounded by an 
innermost component of $\Pi_n\cap\bd C_{n+1}$. The union of each such 
family of concentric disks is a plane which by construction is proper 
in $V$. It is non-trivial in $V$ since otherwise one could use the copy 
of \hs\ which it splits off in $V$ to isotop $\kappa$ so as to reduce 
the number of components of $\kappa\cap\Pi$, a contradiction. \end{proof}

%\section{Nested end reductions}

%Section 7
\section{Some gluing lemmas}

Let $Q$ be a compact \tm. Let $F$ be a compact surface in $\bd Q$, and let $G$ be 
a union of components of $\bd Q-\inte F$. The triple $(Q,F,G)$ has the 
\textit{halfdisk property} if whenever $D$ is a proper disk in $Q$ such that 
$D \cap F$ and $D \cap G$ are each arcs and the union of these arcs is $\bd D$, then $\bd D=\bd D\p$ 
for some disk $D\p$ in $\bd Q$. The ordered triple $(Q,F,G)$ has the \textit{band property} 
if whenever $D$ is a proper disk in $Q$ such that $D \cap F$ consists of two disjoint 
arcs which lie in the same component of $F$ and the remainder of $\bd D$ lies in $G$, then 
$\bd D=\bd D\p$ for some disk $D\p$ in $\bd Q$. Note that if $Q$ is \birr, then 
$(Q,F,G)$ automatically has both properties. 

%Lemma 7.1
\begin{lem} Let $Y$ be a compact \tm. Let $S$ be a compact, 2-sided, proper 
surface in $Y$. Let $Y\p$ be the \tm\ obtained by splitting $Y$ along $S$. 
Let $S\p$ be the surface in $\bd Y\p$ \hm\ to two copies of $S$ which are 
identified to obtain $Y$. Let $Z$ be a union of components of $\bd Y$ and 
$Z\p$ the surface in $\bd Y\p$ obtained by splitting $Z$ along $Z \cap S$. 
Suppose that 

(1) $Y\p$ is \irr, 

(2) $S\p$ and $Z\p$ are \inc\ in $Y\p$, and 

(3) $(Y\p,S\p,Z\p)$ has the halfdisk property. 

Then 

(a) $Y$ is \irr, and 

(b) $Z$ in \inc\ in $Y$. \end{lem}

\begin{proof} This is a standard argument which will be left to the reader. \end{proof}

%Lemma 7.2
\begin{lem} Suppose $Y$, $S$, and $Z$ satisfy all the hypotheses of Lemma 5.1.  
Let $Z_0$ be a union of components of $Z$ and $Z_0\p$ the surface in $\bd Y\p$ 
obtained by splitting $Z_0\p$ along $Z_0\cap S$. Suppose that  

(4) every proper \inc\ annulus in $Y\p$ whose boundary lies in $(\inte S\p) \cup 
(\inte Z\p)$ either is parallel to an annulus in $\bd Y\p$ or cobounds a 
compact submanifold of $Y\p$ with an annulus in $\bd Y\p$ which meets $Z_0\p$. 

(5) $(Y\p,S\p,Z\p)$ has the band property,

(6) each component of $S$ separates $Y$, and 

(7) no component of $S$ is an annulus.  

Then 

(c) every proper \inc\ annulus in $Y$ whose boundary lies in $Z$ either is parallel to an 
annulus in $\bd Y$ or cobounds a compact submanifold of $Y$ with an annulus in $Z_0$. \end{lem}

\begin{proof} Let $A$ be a proper \inc\ annulus in $Y$ such that $\bd A$ lies in $Z$. We 
assume that $A$ is in general position with respect to $S$ and that $A \cap S$ has a 
minimal number of components. We may assume that no component is a simple closed curve 
which bounds a disk in $A$ or is an arc whose boundary lies in one component of $\bd A$. 

Suppose $A\cap S=\ns$. Then $A$ lies in $Y\p$. If $A$ is parallel in $Y\p$ to an 
annulus $A\p$ in $\bd A\p$, then $A\p$ must lie in $\bd Y$ since otherwise it would 
contain an annulus component of $S\p$, and thus $S$ would have an annulus component. 
If $A\cup A\p=\bd Q$, where $Q$ is a compact submanifold of $Y\p$ and $A\p$ is 
an annulus in $\bd Y\p$ which meets $Z_0$, then $A\p$ must lie in $Z_0\p$ since 
otherwise it would contain an annulus component of $S\p$. Thus $Q$ must lie in $Y$.  

Suppose $A \cap S$ contains a simple closed curve. Then there is such a curve \al\ 
which cobounds an annulus $A_0$ in $A$ with a component $\be$ of $\bd A$. 
No component of $A \cap S$ is an arc. We have that $\al$ lies in $\inte S\p$ and \be\ 
lies in $\inte Z\p$, so $A_0$ cobounds a compact submanifold $Q$ of $Y\p$ with an 
annulus $A_0\p$ in $\bd Y\p$.  Since $S$ has no annulus components $A_0\p$ must be 
the union of an annulus $A_0\pp$ in $Z\p$ and an annulus $A_1$ in $S\p$ which meet in 
a component $\beta\p$ of $Z\p\cap S\p$. If $A_0$ is \pl\ to $A_0\p$ across $Q$, then 
an isotopy of $A$ in $Y$ which moves $A_0$ it across $Q$ and then past $A_1$ 
removes at least $\alpha$ from the intersection, thereby contradicting minimality. 
If $A_0$ is not \pl\ to $A_0\p$ across $Q$, then $A_0\pp$ must lie in $Z_0\p$. 
$A_1$ meets $A$ in $\alpha$ and possibly other components. 

Suppose first that 
$A_1\cap A=\alpha$. Now $\alpha$ splits $A$ into two subannuli, one of which is $A_0$. 
Let $A\p$ be the union of the other subannulus and $A_1$. $A\p$ is a proper 
\inc\ annulus in $Y$ which can be isotoped so that it meets $S$ in one fewer 
component than $A$ did. Apply induction to the number of components of the 
intersection. Then $A\p$ cobounds a compact submanifold $Q\p$ of $Y$ with an annulus 
$A\pp$ in $\bd Y$. If $A\pp$ does not contain $A_0\pp$, then $Q\cap Q\p=A_1$ and 
$A\pp$ lies in $Z_0$. In this case $Q\cup Q\p$ is the desired compact submanifold 
of $Y$ cobounded by $A$ and the annulus $A_0\pp\cup A\pp$. If $A\pp$ contains $A_0\pp$, 
then $Q\p$ contains $Q$. In this case the desired compact submanifold is the closure 
in $Y$ of $Q\p-Q$. It is cobounded by $A$ and the annulus in $Z_0$ which is the 
closure of $A\pp-A_0\p$. 

Suppose now that $A_1\cap A$ has some components other than $\alpha$. Then there must be 
some annulus components of $A\cap Q$ other than $A_0$. Suppose there is such a 
component $A_2$ with $\bd A_2$ in $\inte A_1$. Then $A_2$ is a proper \inc\ annulus 
in $Y\p$ with $\bd A_2$ in $\inte S\p$. So $A_2$ cobounds a compact submanifold $Q\p$ 
of $Y\p$ with an annulus $A_2\p$ in $\bd Y$. $A_2\p$ must lie in $A_1$, for otherwise 
$S$ must have an annulus component. Thus $A_2\p$ cannot meet $Z_0\p$, and so it 
is \pl\ to $A_2$ across $Q\p$. Hence an isotopy of $A$ could be performed which 
would reduce the number of intersection curves, contradicting minimality. 
It follows that there is only one component $A_2$ of $A\cap Q$.  
It joins a curve $\alpha\p$ in $\inte A_1$ to a curve $\beta\pp$ in $\inte A_0\p$. 
There is a subannulus $A_2\p$ of $A_1$ with $\bd A_2\p=\beta\p\cup\alpha\p$. 
$A_2\p\cap A=\al\p$. Now $\alpha\p$ splits $A$ into two subannuli, one of which is 
$A_2$. Let $A\p$ be the union of $A_2\p$ and the other subannulus. $A\p$ is a proper 
\inc\ annulus in $Y$ with $\bd A\p=\beta\cup\beta\p$, which lies in $Z_0$. 
$A\p\cap S$ has one fewer component than did $A\cap S$. We apply induction to 
conclude that $A\p$ satisfies (?). Let $Q\p$ be the compact submanifold of $Y$ 
cobounded by $A\p$ and an annulus $A\pp$ in $\bd Y$. Now $A\pp$ either equals $A_0\p$ 
or the component of $Z_0$ containing it is a torus which is the union of these two 
annuli. In the first case $Q\p$ contains $Q$, and $A_2$ is a proper \inc\ annulus 
in $Q\p$ which splits $Q\p$ into two submanifolds, one of which is the desired 
submanifold of $Y$ bounded by $A\cup A_0\p$. In the second case $Q\p$ meets $Q$ in 
$A_0\cup A_2\p$. $A_2$ splits $Q$ into two components. The union of $Q\p$ with 
the component which meets it in $A_2\p$ is the desired compact submanifold.  

Suppose some component \al\ of $A\cap S$ is an arc. Then it must meet both components 
of $\bd A$. By (6) $A\cap S$ must hve an even number of components, and there must be a 
component $\Delta$ of $A \cap Y\p$ which meets some component of $S\p$ in two disjoint 
arcs and has the remainder of its boundary in $Z\p$. By (5) $\bd \Delta = \bd \Delta\p$ 
fo some disk $\Delta\p$ in $\bd Y$. $\Delta \cup \Delta\p$ bounds a 3-ball in $Y\p$, 
and an isotopy of $\Delta$ across this 3-ball reduces the number of components, again 
contradicting minimality. So this case cannot occur.  \end{proof}

%Lemma 7.3
\begin{lem} Suppose $Y$, $S$, $Z$, and $Z_0$ satisfy all the hypotheses of 
Lemmas 7.1 and 7.2. Suppose that 

(8) every \inc\ torus in $Y\p$ bounds a compact submanifold of $Y\p$, and 

(9) no component of $\bd Y$ is a torus.  

Then 

(d) every \inc\ torus in $Y$ is bounds a compact 
submanifold of $Y$. \end{lem}

\begin{proof} Let $T$ be an \inc\ torus in $Y$. Put $T$ in general position 
with respect to $S$ so that $T\cap S$ has a minimal number of components. 
We may assume that no component bounds a disk in either surface. 

Suppose $T\cap S=\ns$. Then $T$ lies in $Y\p$, so $T$ bounds a compact 
submanifold $Q$ of $Y\p$. Since $Q$ lies in $Y$ we are done. 

Suppose $T\cap S\neq\ns$. Then $S$ splits $T$ into annuli which lie 
in $Y\p$. Since each component of $S$ separates $Y$ there is such 
an annulus $A$ with $\bd A$ in a single component of $S\p$. 
Then there is an annulus $A\p$ 
in $\bd Y\p$ such that $A\cup A\p=\bd Q$ for a compact submanifold 
$Q$ of $Y\p$. 

Assume that $A\p$ lies in $S\p$. Then $A\p$ does not meet $Z_0$, and 
so $A$ is \pl\ across $Q$ to $A\p$. Thus there is an isotopy of $T$ in $Y$ 
which removes at least $\bd A$ from the intersection, contradicting minimality. 
Hence $A\p$ does not lie in $S\p$. Since $S$ has no annulus components 
$A\p$ must consist of an annulus $A_1$ in $\bd Y$ together with a collar 
$C$ in $S\p$ on two components of $\bd S\p$. Let $A_2$ be the union of 
$C$ and $T-\inte A$. Note that $\bd A_2=\bd A_1$. Now $A_2$  is a proper 
\inc\ annulus in $Y$ which can be 
isotoped to have boundary disjoint from $\bd S$. Therefore $A_2$ cobounds 
a submanifold $Q\p$ of $Y$ with an annulus $A_2\p$ in $\bd Y$. If $A_2\p=A_1$, 
then $Q\p$ contains $Q$ and $A$ is a proper annulus in $Q\p$ which splits it 
into two components, 
one of which is the desired submanifold having boundary $T$. If $A_2\p\neq A_1$, 
then the union of these two annuli is a torus component of $\bd Y$ contrary 
to the hypotheses on $Y$. \end{proof}

%Section 8
\section{Some building blocks}

In this section we define certain compact \tm s which will be used to construct 
our examples and establish some of their properties which will be used in 
conjunction with the gluing lemmas of the previous section. 

We first introduce some general notation. Let $\Gamma$ be a graph in a 
\tm\ $Q$. We assume that $\Gamma\cup\bd Q$ is either empty or consists of 
vertices of order one. We denote a regular neighborhood of $\Gamma$ in $Q$ by 
$N(\Gamma,Q)$. The \textit{exterior} 
of $\Gamma$ in $Q$ is the closure of the complement of $N(\Gamma,Q)$ in $Q$ and is denoted 
by $X(\Gamma,Q)$. The \textit{lateral surface} determined by $\Gamma$ is the set 
$S(\Gamma,Q)=N(\Gamma,Q) \cap X(\Gamma,Q)$. The $Q$ will be suppressed when its identity 
is clear from the context. If $\gamma$ is an edge of $\Gamma$, then a simple 
closed curve in $S(\Gamma,Q)$ which bounds a disk in $N(\Gamma,Q)$ which meets 
$\gamma$ transversely in a single point is called a \textit{meridian} of $\gamma$.  

Let $P=D\times [0,2]$, where $D$ is a closed disk. Set $L=D\times [0,1]$, $R=D\times [1,2]$, 
and $D_j=D\times \{j\}$ for $j=0,1,2$. Attach a 1-handle $H$ to $L$ along $\bd D\times(0,1)$ 
to obtain a solid torus $L\cup H$. . Set $G=\bd(L\cup H)-\inte(D_0\cup D_1)$ and 
$O=\bd R-\inte(D_1\cup D_2)$. 
Let $J=P\cup H$. Let $J^\#$ be the genus two handlebody obtained from $J$ by identifying 
$D_0$ and $D_2$. Let $P^\#$ be the solid torus in $J^\#$ which is the image of $P$ under 
the identification. With the exception of $J$ and $J^\#$ and of $P$ and $P^\#$ we will 
usually use the same notation for subsets of $J$ and their images in $J^\#$, relying on 
the context for which is meant. Thus we regard $J^\#$ as $P^\# \cup H$. 

We next let \be\ and \ep\ be the two disjoint arcs in the solid torus $L\cup H$ shown in 
Figure 1. Each of them runs from $\inte D_0$ to $\inte D_1$. We have that \be\ 
meets $H$ in a Whitehead clasp while \ep\ is a product arc in $L$. 
This configuration is called an \textit{eyebolt}.

%INSERT FIGURE 1 HERE.  

\begin{figure}
\epsfig{file=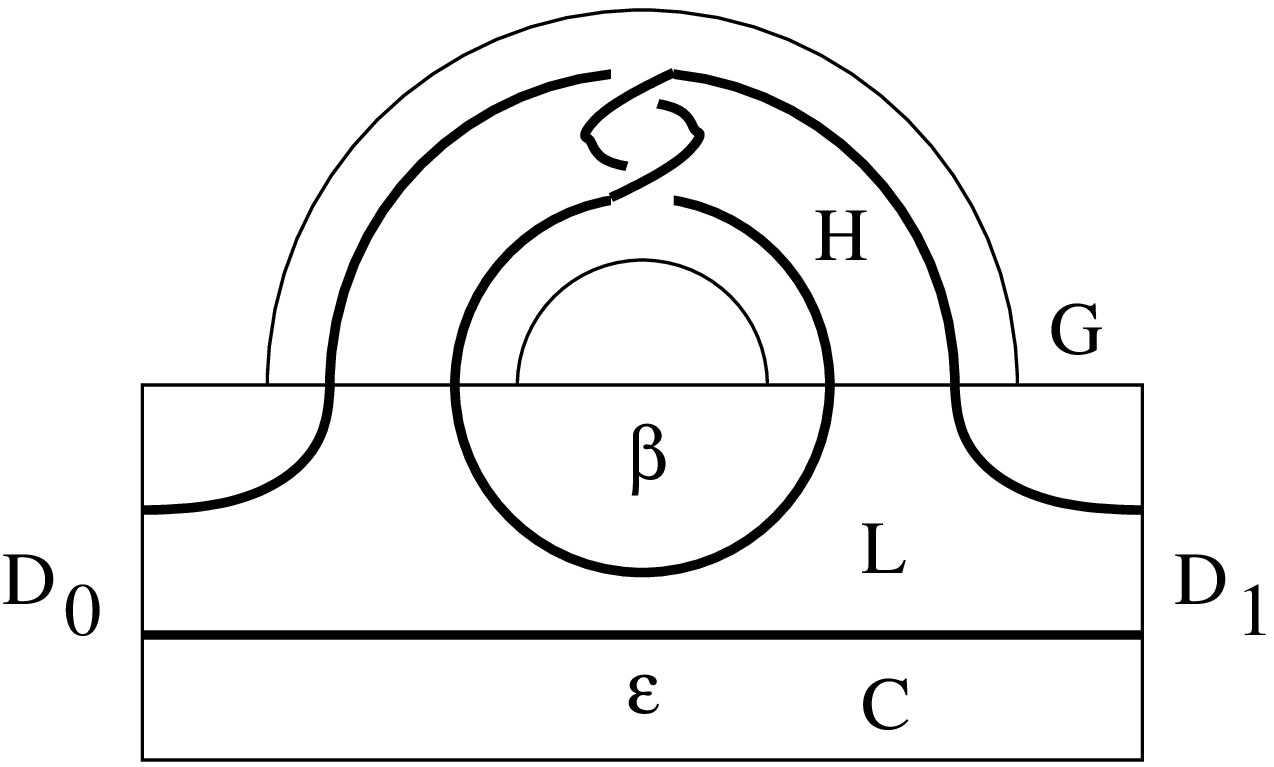, width=4in}
\caption{An eyebolt $(L\cup H,\beta\cup\ep)$}
\end{figure}

We then let \al, \ga, \de, \ze, and \rh\ be the arcs in the 3-ball $R$ shown in Figure 2. 
We require that $\ga\cap D_1=\be\cap D_1$, that $\de \cap D_1=\ep\cap D_1$, 
and that $\al \cap D_2$ and $\ze \cap D_2$, respectively, are identified with 
$\be \cap D_0$ and $\ep \cap D_0$ in $J^\#$. We have that \al, \ze, and \rh\ meet in 
a common endpoint while \ga, \de, and \rh\ also meet in a common endpoint. Except as 
indicated all the arcs are disjoint. We have that $\al \cup\ze$ and $\ga\cup\de$ are 
linked trefoil knotted arcs whose union is a tangle \ta\ in $R$. We have that 
$\de\cup\rh\cup\ze$ is a product arc in $R$. Let $\rh\p=\rh\cap X(\tau,R)$. 
We may assume that \rh\ is the union of $\rh\p$ with two disjoint arcs in 
$N(\tau,R)$. This configuration is called a textit{junction}

%INSERT FIGURE 2 HERE.  

\begin{figure}
\epsfig{file=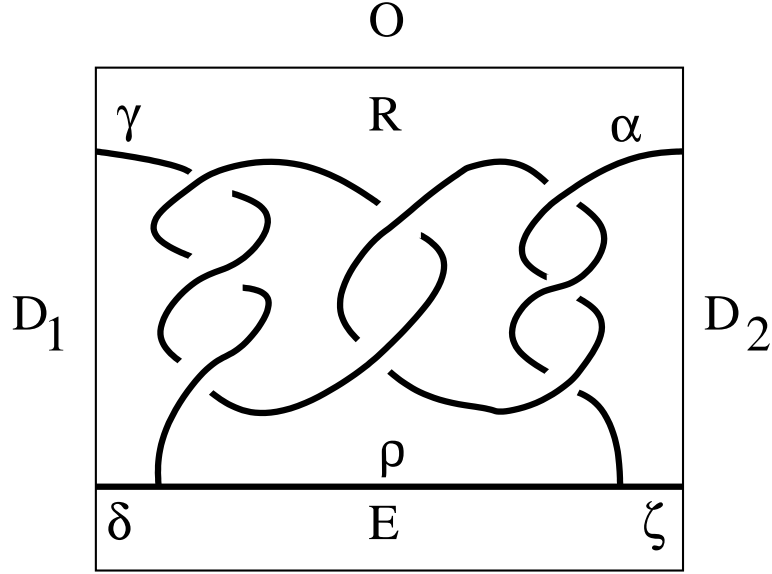, width=4in}
\caption{A junction $(R,\alpha\cup\gamma\cup\delta\cup\zeta\cup\rho)$}
\end{figure}

In $J^\#$ set $\lambda=\de\cup\ep\cup\ze$, $\eta=\al\cup\be\cup\ga$, 
$\ka=\la\cup\eta$, and $\mu=\ka\cup\rh$ as in Figure 3. 

%INSERT FIGURE 3 HERE. 

\begin{figure}
\epsfig{file=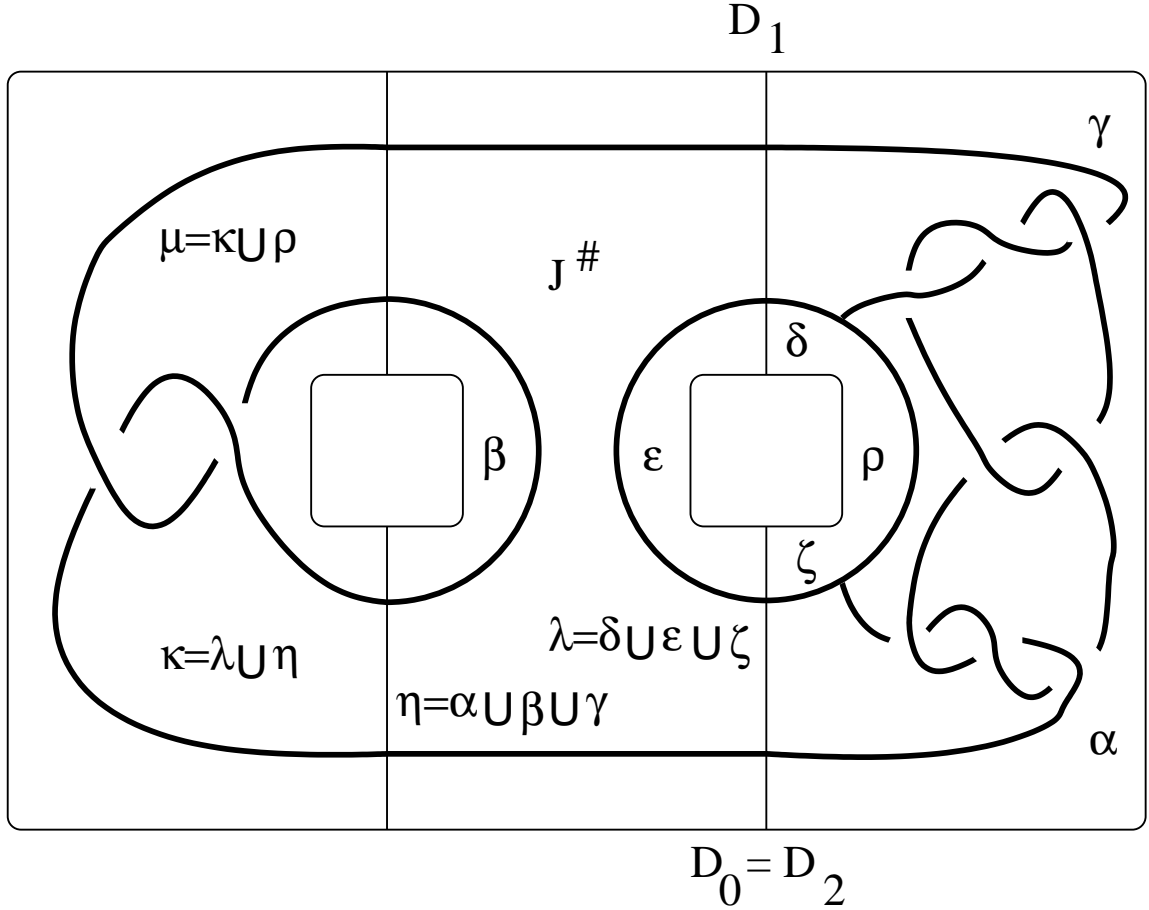, width=6in}
\caption{The arc $\lambda$, arc $\eta$, knot $\kappa$, and graph $\mu$ 
in the genus two handlebody $J^{\#}$} 
\end{figure}

%Lemma 8.1
\begin{lem} Let $Y=X(\be \cup \ep, L \cup H)$, and let $F_i=Y \cap D_i$ for 
$i=0,1$. Then 

(0) $Y$ is \irr, 

(1) $F_i$  is \inc\ in $Y$, 

(2) $G$ is \inc\ in $Y$, 

(3) $S(\ep,L\cup H)$ is \inc\ in $Y$, 

(4) $S(\be,L\cup H)$ is \inc\ in $Y$, 

(5) $(Y,F_0\cup F_1, \bd Y-\inte(F_0\cup F_1))$ has the halfdisk property, 

(6) $(Y,F_0\cup F_1, \bd Y-\inte(F_0\cup F_1))$ has the band property. 

(7) Every proper \inc\ annulus in $Y$ whose boundary lies in 
$(\inte(F_0\cup F_1)\cup(\bd Y-(F_0\cup F_1))$ is \pl\ to an annulus 
in $\bd Y$, and 

(8) $Y$ contains no \inc\ tori. 
\end{lem}

\begin{proof} For $i=0,1$ attach 
3-balls $B_i$ to $L\cup H$ so that $B_i\cap(L\cup H)=D_i$; let 
$D\p_i=\bd B_i-\inte D_i$.  Let $\omega_i$ be 
a proper unknotted arc in $B_i$ joining the points of $(\be\cup\ep)\cap D_i$. 
Then $\be\cup\ep\cup\omega_0\cup\omega_1$ is a Whitehead curve in the solid 
torus $K=L\cup H\cup B_0\cup B_1$. Let $X$ be the exterior of this curve in 
this solid torus. $X$ is \irr, \birr, anannluar, and atoroidal 
\cite{My simple}.   

(0) and (8) hold because $Y$ is a handlebody. 

(1) holds for homological reasons. 

(2) Suppose $D$ is a compressing disk for $G$ in $Y$.  
$\bd D=\bd D\p$ for a disk $D\p$ on $\bd K$.   
If $D\p$ contains, say, $D_0\p$ but not $D_1\p$, then 
$\bd D$ is \pl\ in $G$ to $\bd D_0$, and so $F_0$ is compressible 
in $Y$, contradicting (1). Assume $D\p$ contains $D_0\p$ and $D_1\p$. 
Then $D$ splits $K$ into a solid torus and a 3-ball 
$B\p$ such that $B_0$ and $B_1$ lie in $B\p$. Thus 
$\be\cup\ep\cup\omega_0\cup\omega_1$ 
lies in $B\p$, contradicting the irreducibility of $X$. 
Thus we have that $D\p$ lies in $G$. 

(3) $S(\ep,L\cup H)$ is \inc\ for homological reasons. 

(4) $S(\be,L\cup H)$ is \inc\ for homological reasons. 

(5) Suppose $D$ is a proper disk in $Y$ such that $D\cap(F_0\cup F_1)$ is an arc 
$\theta$ in, say, $F_0$, and $\bd D-\inte \theta=D\cap(\bd Y-\inte(F_0\cup F_1))$ 
is an arc $\varphi$. If $\varphi$ lies in one of the annuli $S(\be,L\cup H)$ or 
$S(\ep,L\cup H)$, then $D$ can be isotoped so that $\bd D$ lies in $F_0$, and so 
the result follows from (1). So assume that $\varphi$ lies in $G$. By (2) we may 
assume that $\theta$ separates $\be\cap D_0$ from $\ep\cap D_0$. 
There is a disk $E$ in $B_0$ which meets $\omega_0$ 
transversely in a single point and meets $D_0$ in $\theta$. Then $D\cup E$ is a 
proper disk in the solid torus $K$ which meets the 
Whitehead curve transversely in a single 
point, thereby contradicting the anannularity of $X$.  

(6) Suppose $D$ is a proper disk in $Y$ such that $D\cap(F_0 \cup F_1)$ consists 
of two arcs $\theta\p$, $\theta\pp$ in, say, $F_0$ and such that 
$\bd D-\inte(\theta\p\cup\theta\pp)=D\cap(\bd Y-\inte(F_0\cup F_1))$ consists 
of the two arcs $\varphi\p$ and $\varphi\pp$. If $\theta\p$ or $\theta\pp$ 
is \bpl\ in $F_0$ or $\varphi\p$ or $\varphi\pp$ is \bpl\ in $\bd Y-\inte(F_0\cup F_1)$, 
then there is an isotopy of $D$ after which it meets $F_0$ and $\bd Y-\inte(
F_0\cup F_1)$ each in a single arc, and the union of these arcs is $\bd D$. 
Then by (5) we are done. So we may assume that none of $\tht\p$, $\tht\pp$, 
$\varphi\p$, or $\varphi\pp$ are \bpl. In particular this implies that 
$\varphi\p\cup\varphi\pp$ lies in $G$ and that $\tht\p$ and $\tht\pp$ are \pl\ 
arcs in $F_0$ which join $\bd D_0$ to itself and separate $\be\cap D_0$ from 
$\ep\cap D_0$. 

Put $D$ in general position with respect to $L\cap H\cap Y$ so that the 
intersection has a minimal number of components. $L\cap H \cap Y$ is \inc\ 
in $X(\be\cap H, H)$ by \cite{My simple}  
and is \inc\ in $X((\be\cup\ep)\cap L,L)$ for 
homological reasons. These two spaces are \irr\ since they are handlebodies. 
Thus $D\cap L\cap H\cap Y$ has no simple closed curve components. 

If $\nu$ 
is an outermost arc on $D$ of this intersection such that  $\bd\nu$ lies in, 
say, $\varphi\p$, 
then $\bd \nu=\bd\nu\p$ for an arc $\nu\p$ in $\inte\varphi\p$ such that 
$\nu\cup\nu\p=\bd\Delta$ for a disk $\Delta$ in $D$. 

Suppose $\nu\p$ lies in the annulus 
$G\cap H$. Then it is \bpl\ in $G\cap H$, and so by the incompressibility 
of $L\cap H\cap Y$ in $X(\be\cap H,H)$ we have that $\nu$ is \bpl\ in 
$L\cap H\cap Y$. It follows that $\nu$ could have been removed by an isotopy 
of $D$, thereby contradicting minimality. 

Suppose $\nu\p$ lies in the disk with three holes $G\cap L$. Since $\Delta$ 
misses $\ep$ we must have that $\bd \Delta=\bd\Delta\p$ for a disk 
$\Delta\p$ in the annulus $\bd L-\inte(D_0\cup D_1)$. Let $\Delta\pp$ 
be the component of $H\cap L$ which does not contain $\nu$. 
Then $\Delta\p$ either contains $\Delta\pp$ or is disjoint from it. 
In the first case $\Delta$ would separate $\Delta\pp$ from $D_0\cup D_1$, which 
is impossible since they are joined by a subarc of $\beta$. Thus $\Delta\p$ 
is disjoint from $\Delta\pp$. It follows that $\nu\p$ is \bpl\ in $G\cap L$. 
By the incompressibility of $L\cap H\cap Y$ in $X((\beta\cup\ep)\cap L, L)$ 
we have that $\nu$ is \bpl\ in $L\cap H\cap Y$. Again $\nu$ could have 
been removed by an isotopy of $D$, contradicting minimality.

So every component of 
$D\cap L\cap H\cap Y$ must join $\varphi\p$ to $\varphi\pp$ 
\cite{My simple}. Since $L\cap H\cap Y$ 
separates $Y$ there must be an even number of such components. Let $\nu_0$ 
and $\nu_1$ be adjacent such components which are \pl\ in $D$ across a disk 
$\Delta$ which meets $\varphi\p$ and $\varphi\pp$ in arcs $\nu\p$ and $\nu\pp$, 
respectively. If $D\cap H\neq\ns$, then there must be such a $\Delta$ lying 
in $X(\be\cap H, H)$. Suppose that 
$\nu_0$ and 
$\nu_1$ lie in the same component of $L\cap H\cap Y$, and so $\nu\p$ 
and $\nu\pp$ are \bpl\ in $G\cup H$. It follows from the incompressibility 
of $L\cap H\cap Y$ and the irreducibility of $X(\be\cap H,H)$ that there is 
an isotopy of $D$ which removes $\nu_0\cup\nu_1$ from the intersection, again 
contradicting minimality. Therefore for every such $\Delta$ we have  that $\nu_0$ 
and $\nu_1$ lie in different components of $L\cap H\cap Y$. By 
\cite{My simple} $\Delta$ 
is \bpl\ in this space. Let $\Delta_0$ be the component of $D\cap L$ containing 
$\theta\p$. There is a $\Delta$ as above which meets $\Delta_0$ in the arc $\nu_0$. 
Since $\nu_0$ is \bpl\ in $L\cap H\cap Y$ we can isotop $\Delta_0$ in $L\cap Y$ 
to a disk $\Delta_0\p$ which meets $D_0$ in $\theta\p$ and is disjoint from 
$L\cap H$ and $D_1$. $\Delta_0\p$ splits $L$ into two 3-balls $B\p$ and $B\pp$ 
with $\beta\cap D_0$ contained in, say, $B\p$ and $\ep$ contained in $B\pp$. 
The component of $\beta\cap L$ 
meeting $D_0$ must lie in $B\p$. Thus the component of $L\cap H$ meeting this 
arc must lie in $B\p$. Hence so must the component of $\beta\cap L$ joining 
the two components of $L\cap H$. Hence so must the other component of $L\cap H$. 
Hence so must the third component of $\beta\cap L$. Hence so must $D_1$. 
Hence so must $\ep$, a contradiction.    

Thus we now have that $D\cap H=\ns$. Now $\bd D=\bd D\p=\bd D\pp$, where 
$D\p$ and $D\pp$ are disks on $\bd L$ whose union is $\bd L$. We may assume 
that $D\pp$ contains $D_0\cap(\be\cup\ep)$. Then for homological reasons we 
must have that $D\pp$ also contains $\be\cap L\cap H$ and $D_1$. Thus $D\p$ 
lies in $\bd Y$, and we are done. 

(7) Let $A$ be an annulus with $\bd A$ in 
$(\inte(F_0\cup F_1)\cup(\bd Y-(F_0\cup F_1))$. 

Note that $H_1(G)$ has basis consisting of classes $d$, $\ell$, and $m$ 
represented by, respectively, $\bd D_0$, and a longitude and a 
meridian of the solid torus $L\cup H$. $H_1(Y)$ has basis consisting 
of classes $b$, $e$, and $\ell_0$ represented by, respectively, by 
meridians of $\beta$ and $\ep$ and a longitude of $L\cup H$. 
These curves and their orientations can be chosen so that under 
the inclusion induced map $d\rightarrow b+e$, $\ell\ra\ell_0$, and 
$m\ra0$. It follows that if $\theta\p$ is a meridian of $\beta$ or 
$\epsilon$, respectively, then $\theta\pp$ must be also. Hence in 
this case we may assume that $\bd A$ lies in a single component 
of $S(\beta\cup\epsilon,L\cup H)$. Since $\bd A$ consists of 
meridians of the Whitehead curve $A$ is \inc\ in $X$ and thus is 
\pl\ to an annulus $A\p$ in $\bd X$. Thus $A$ splits $X$ into a solid 
torus containing $A\p$ and another 3-manifold containing $\bd K$. 
It follows that $A\p$ and the solid torus lie in $Y$, so we 
are done in this case. 

Now suppose that $\bd A$ lies in $G$. 

Let $F$ be the disk in Figure 1. Let $F\p=F\cap Y$. Put $A$ in 
general position with respect to $F\p$ so that $A\cap F\p$ has 
a minimal number of components. Then this intersection has no 
simple closed curve components. 

Suppose some component $\nu$ of $A\cap F\p$ is \bpl\ in $A$. 
We may assume that $\nu$ cuts off an outermost disk $D$ on $A$ 
and that $\bd D=\nu\cup\nu\p$ for an arc $\nu\p$ in $\theta\p$. 
Then $\nu$ cuts off a disk $D\p$ in $F\p$. Since $G$ is \inc\ 
in $Y$ we have that $\bd(D\cup D\p)=\bd D\pp$ for a disk $D\pp$ 
in $G$. An isotopy of $A$ which moves $D$ across the 3-ball 
bounded by $D\cup D\p\cup D\pp$ and then past $D\p$ removes at 
least $\nu$ from $A\cap F\p$, thereby contradicting minimality. 

Thus every component $\nu$ of $A\cap F\p$ is a spanning arc 
on $A$. Suppose $\nu$ is outermost on $F\p$, cutting off an 
outermost disk $D$ on $F\p$. Performing a boundary compression 
on $A$ along $D$ yields a disk $D\p$ with $\bd D\p$ in $G$. 
Then $\bd D\p=\bd D\pp$ for a disk $D\pp$ in $G$. Let $B$ be the 
3-ball in $Y$ bounded by $D\p\cup D\pp$. $B$ cannot contain $D$ 
since this would make $A$ comopressible in $Y$. It follows that 
the union of $B$ and an appropriate regular neighborhood of $D$ 
is a solid torus across which $A$ is \pl\ to an annulus in $\bd Y$. 

Thus we may assume that $A\cap F\p=\ns$. 

Suppose $A$ is compressible in $X$. Since $X$ is \birr\ each 
component of $\bd A$ must bound a disk in $\bd X$, say 
$\theta\p=\bd\Delta\p$ and $\theta\pp=\bd\Delta\pp$. Since 
$G$ is \inc\ in $Y$ we must have that $\Delta\p$ and $\Delta\pp$ 
each contains at least one of the disks $D_0\p$ and $D_1\p$. 
We may assume that $\Delta\p$ contains $D_0\p$ and that $\Delta\pp$ 
either contains $\Delta\p$ or is disjoint from it. If $\Delta\p$ 
does not contain $D_1\p$, then $\Delta\p\cap F\p\neq\ns$, hence 
$A\cap F\p\neq\ns$, a contradiction. So 
$(D_0\p\cup D_1\p)\subseteq\Delta\p\subseteq\Delta\pp$. 
Therefore $\Delta\p$ and $\Delta\pp$ are isotopic to concentric 
regular neighborhoods of $D_0\p\cup D_1\p\cup(F\cap G)$ in $\bd K$. 

Put $A$ in general position with respect to $L\cap H\cap Y$ so that 
the intersection has a minimal number of components among annuli 
in $Y$ with $\bd A=\bd\Delta\p\cup\Delta\pp$. Then each component 
is a simple closed curve which bounds a disk on neither $A$ nor 
$L\cap H\cap Y$. So there is an annulus component $A\p$ of $A\cap 
Y\cap L$ with $\bd A\p=\theta\p\cup\theta$, 
where $\theta$ is a curve in $L\cap H\cap Y$. It is easily checked 
that this is homologically impossible. 

Thus we have that $A$ lies 
in $Y\cap L$. Let $C$ be the disk in Figure 1. Let $C\p=C\cap Y\cap L$. 
$A$ must meet $C\p$. Assume, as usual, general position and minimality. 
Each component is an arc. Suppose $\nu$ is an outermost \bpl\ arc on 
$A$ cutting off a disk $D$ on $A$. Then $\nu$ cuts off a disk $D\p$ on 
$C\p$. Since $G$ is \inc\ $\bd(D\cup D\p)=\bd D\pp$ for a disk $D\pp$ 
on $G$. An isotopy of $D$ across the 3-ball bounded by $D\cup D\p\cup D\pp$ 
and past $D\p$ removes at least $\nu$ from the intersection. 
Hence we may assume that all components are spanning arcs of $A$. 
Let $\nu$ be a component which is outermost on $C\p$, cutting off a 
disk $D$ on $C\p$ with $\bd D=\nu\cup\nu\p$ for an arc $\nu\p$ 
in $C\p\cap G$. Performing a boundary compression on $A$ along $D$ 
yields a disk $D\p$ with $\bd D\p$ in $G$. $\bd D\p=\bd D\pp$ for a 
disk $D\pp$ in $G$. The 3-ball $B$ bounded by $D\p\cup D\pp$ cannot 
contain $D$ as this would make $A$ compressible in $Y$. So the union 
of $B$ and an appropriate regular neighborhood of $D$ is a solid 
torus in $Y$ across which $A$ is \pl\ to an annulus in $\bd Y$. 

Finally, suppose that $A$ is \inc\ in $X$. Let $A\p$ be the 
annulus in $\bd X$ to which $A$ is \pl. Suppose, say,  $D_0\p$ lies 
in $A\p$. Since $A\cap D_0=\ns$ we have that $B_0$ lies in the 
solid torus bounded by $A\cup A\p$. It follows that $\bd X-\bd K$ 
lies in this solid torus, which is impossible. 
\end{proof}

The configuration of $\beta\cup\gamma\cup\delta\cup\epsilon$ in 
$J$ is called a \textit{right hitch}. There is an obvious twin 
configuration called a \textit{left hitch}. 

%Lemma 8.2
\begin{lem} Let $Y=X(\beta\cup\gamma\cup\delta\cup\epsilon,J)$. Then 

(0) $Y$ is \irr, 

(1) $F_0$ is \inc\ in $Y$, 

(2) $\bd J-\inte D_0$ is \inc\ in $Y$, 

(3) $S(\beta\cup\gamma\cup\delta\cup\epsilon,J)$ is \inc\ in $Y$, 

(4) $(Y,F_0,\bd Y-\inte F_0)$ has the halfdisk property, 

(5) $(Y,F_0,\bd Y-\inte F_0)$ has the band property, 

(6) every proper \inc\ annulus in $Y$ whose boundary misses $\bd F_0$ is 
either \bpl\ or is isotopic to $F_1\cup O\cup D_2$, and 

(7) every \inc\ torus in $Y$ bounds a compact submanifold of $Y$. 

\end{lem}

\begin{proof} (0) Attach a 3-ball $B_0$ and unknotted arc $\omega_0$ to 
$P\cup H$ along $D_0$ as in 
the proof of Lemma 6.1. The result is a solid torus $K$ containing a Whitehead 
curve in which one has locally tied a trefoil knot. The exterior $X$ of 
this curve is homeomorphic to the union of the Whitehead link exterior 
and a trefoil knot exterior along an annulus which is \inc\ in both. 
This implies that $X$ is \irr\ and \birr. So any 2-sphere in $Y$ must bound a 
3-ball in $X$. Since $X(\omega_0, B_0)$ has torus boundary this 3-ball 
must lie in $Y$. 

(1) Let $D$ be a compressing disk for $F_0$. For homological 
reasons $\bd D$ must be isotopic in $F_0$ to $\bd D_0$. Since an arc 
of $\beta\cap L$ joins $D_0$ to $L\cap H$ we must have that $D$ intersects 
$H$. Suppose $\theta$ is a component of $D\cap L\cap H\cap Y$ which bounds 
an innermost disk $D\p$ on $D$. For homological reasons $D\p$ cannot lie 
in $P\cap Y$. So $D\p$ lies in $H\cap Y$. By \cite{} it is \pl\ to a disk 
in $L\cap H\cap Y$, so can be removed by an isotopy. Continuing in this 
fashion we can remove all intersections, contradicting the fact that $D$ 
must meet $H$. 

(2) Let $D$ be a compressing disk for $\bd J-\inte D_0$. 
Since $X$ is \birr\ $\bd D=\bd D\p$ for a disk $D\p$ in $\bd X$. 
Let $D_0\p=\bd B_0-\inte D_0$. If $D\p$ contains $D_0\p$, then $\bd D$ 
is isotopic in $\bd Y$ to $\bd D_0$, contradicting (1). Thus $D\p$ 
lies in $\bd Y$, and we are done. 

(3) This holds for homological reasons. 

(4) Suppose $D$ is a proper disk in $Y$ which lmeets $F_0$ in an arc 
$\theta$ and $\bd Y-\inte F_0$ in an arc $\varphi$. By (1), (2), and (3) 
we may assume that neither arc is \bpl. It follows that either $\theta$ 
joins $\bd D_0$ to itself and separates $\beta\cap D_0$ from 
$\epsilon\cap D_0$ or $\theta$ joins the two components of $\bd F_0-\bd D_0$. 

In the first case attach $B_0$ and $\omega_0$ as in (2). As in the proof 
of Lemma 7.1 (5) we choose a disk $E$ in $B_0$ which meets $\omega_0$ 
transversely in a single point and meets $D_0$ in $\theta$. Then $D\cup E$ 
is a proper disk in the solid torus $B_0\cup J$ which meets a (locally 
knotted) Whitehead curve transversely in a single point, contradicting 
the fact that the corresponding link has linking number zero. 

In the second case we note that $\bd D$ can be regarded as a knot in 
$\RRR$ which is both a trivial knot and a trefoil knot, so this cannot 
occur. 

(5) Suppose $D$ is a proper disk in $Y$ such that $D\cap F_0$ consists 
of two arcs $\theta\p$ and $\theta\pp$. Let $\varphi\p$ and $\varphi\pp$ 
be the components of $\bd D-\inte(\theta\p\cup\theta\pp)$. By (4) we 
may assume that none of these arcs are \bpl. Then $\theta\p$ and $\theta\pp$ 
msut be parallel, and they either join $\bd D_0$ to itself, separating 
$\beta\cap D_)$ from $\epsilon\cap D_0$, or they join the two components of 
$\bd F_0-\bd D_0$. 

In the first case put $D$ in general position with respect to 
$L\cap H\cap Y$ so that the intersection has a minimal number of 
components. $L\cap H\cap Y$ is \inc\ in $H\cap Y$ by \cite{My simple} and in 
$P\cap Y$ for homological reasons. $H\cap Y$ is \irr\ since it is a 
handlebody. $P\cap Y$ is \irr\ since it is the union of a handlebody 
and a knot exterior along an annulus which is \inc\ in both. Thus 
the intersection has no simple closed curve components. 

Let $\nu$ be an outermost arc on $D$ cutting off a disk $\Delta$ 
with $\bd\Delta=\nu\cup\nu\p$, where $\nu\p$ is an arc in, say, $\varphi\p$. 
If $\nu\p$ lies in $G\cap H$, then $D$ can be isotoped to remove $\nu$ as 
in the proof of Lemma 7.1 (6), contradicting minimality. So $\nu\p$ must lie 
in the disk with two holes $\bd P-\inte(D_0\cup(L\cap H))$. Suppose 
$\nu\p$ separates $\bd D_0$ from the boundary component of this surface which 
does not meet $\nu\p$. Then $\Delta$ separates $D_0$ in $P$ from the 
corresponding component of $L\cap H$. But this is impossible since $D_0$ 
is connected to each component of $L\cap H$ by arcs in $P$ which miss 
$\Delta$. Thus $\nu\p$ must be \bpl\ in our disk with two holes. $D$ can 
therefore be isotoped to remove $\nu$, a contradiction. 

So, each component of the intersection must join $\varphi\p$ to $\varphi\pp$. 
Let $\nu_0$ and $\nu_1$ be adjacent such components which are \pl\ in $D$ 
across a disk $\Delta$ which meets $\varphi\p$ and $\varphi\pp$ in arcs 
$\nu\p$ and $\nu\pp$, respectively. If $D\cap H\neq\ns$, then there must 
be such a disk $\Delta$ lying in $H\cap Y$. If $\nu_0$ and $\nu_1$ lie in 
the same component of $L\cap H\cap Y$, then as in the proof of Lemma 7.1 (6) w
we can isotop $D$ to remove $\nu_0\cup \nu_1$, contradicting minimality. 
Thus $\nu_0$ and $\nu_1$ must lie in different components, and so by 
\cite{My simple} $\Delta$ is \bpl\ in $H\cap Y$. We let $\Delta_0$ be the component 
of $D\cap P$ containing $\theta\p$ and $\Delta$ the component of $D\cap H$ 
meeting $\Delta_0$ in the arc $\nu_0$. Since $\nu_0$ is \bpl\ in 
$L\cap H\cap Y$ we can isotop $\Delta_0$ to a disk $\Delta_0\p$ which 
meets $F_0$ in $\theta\p$ and is disjoint from $L\cap H$. $\Delta_0\p$ 
splits $P$ into two 3-balls $B\p$ and $B\pp$ with $\beta\cap D_0$ in $B\p$ 
and $\epsilon\cap D_0$ in $B\pp$. The component of $L\cap H$ joined to 
$\beta\cap D_0$ by a subarc of $\beta$ must lie in $B\p$. The other 
component of $L\cap H$ is joined to this one by another subarc of $\beta$, 
so it must lie in $B\p$. This is joined to $\epsilon\cap D_0$ by the 
union of a third subarc of $\beta$ with $\gamma$, $\delta$, and $\epsilon$, 
so $\epsilon\cap D_0$ must lie in $B\p$, a contradiction. 

Thus we have that $D\cap H=\ns$. Now $\bd D-\bd D\p=\bd D\pp$ for disks 
$D\p$ and $D\pp$ on $\bd P$ whose union is $\bd P$. We may assume that $D\pp$ 
contains $D_0\cap(\beta\cup\epsilon)$. By following subarcs of $\beta$ as 
above we get that $D\pp$ also contains $L\cap H$. Thus $D\p$ lies in $\bd Y$, 
and we are done in this case. 

Recall that in the second case $\theta\p$ and $\theta\pp$ join the two 
components of $\bd F_0-\bd D_0$. Then $\bd D$ can be regarded as a knot 
in \RRR\ which is both a trivial knot and a 2-strand cable of a trefoil 
knot, so this case cannot occur. 

(6) Let $A$ be a proper \inc\ annulus in $Y$ such that 
$(\bd A)\cap (\bd F_0)=\ns$. Denote the components of $\bd A$ by $\theta\p$ 
and $\theta\pp$. Isotop $A$ so that $(\bd A) \cap F_0=\ns$. Attach $B_0$ 
and $\omega_0$ as in (2) to obtain a locally unknotted Whitehead curve in 
a solid torus $K$. Let $X$ be the exterior of this curve in $K$. 

Suppose $A$ is compressible in $X$. Since $X$ is \birr\ each component of 
$\bd A$ bounds a disk on $\bd X$. Since $A$ is \inc\ in $Y$ these disks 
must contain $D_0\p$. Thus each component of $\bd A$ is isotopic in $\bd Y$ 
to $\bd D_0$. Isotop $A$ so that $\bd A$ lies in $F_0$. 

Put $A$ in general position with respect to $L\cap H\cap Y$ so that the 
intersection has a minimal number of components. By the incompressibility 
of $L\cap H\cap Y$ in $Y$ each component is a simple closed curve which 
does not bound a disk on either surface. If $A\cap L\cap H\cap Y\neq\ns$, 
then there is a component $\theta$ of the intersection such that $\theta\cup 
\theta\p=\bd A_0$ for a component $A_0$ of $A\cap P\cap Y$. However, it is 
easily checked that this is homologically impossible. Thus $A\cap H=\ns$. 

Recall the disk $C$ in Figure 1. It intersects $Y$ in a disk $C\p$. 
Put $A$ in general position with respect to $C\p$ so that $A\cap C\p$ 
has a minimal number of components. We must have that $A\cap C\p$ is 
non-empty and contains no simple closed curves. Suppose $\nu$ is a 
component of this intersection which is an outermost \bpl\ arc on $A$, 
cutting off a disk $D$ with $\bd D=\nu\cup\nu\p$, where $\nu\p$ is an 
arc in, say, $\theta\p$. There is an arc $\nu\pp$ in $F_0$ such that 
$\nu\cup\nu\pp$ bounds a disk $D\p$ in $C\p$. Since $F_0$ is \inc\ and 
$P\cap Y$ is \irr\ $\bd(D\cup D\p)=\bd D\pp$ for a disk $D\pp$ in $F_0$, 
and $D\cup D\p\cup D\pp$ bounds a 3-ball $B$ in $P\cap Y$. An isotopy of 
$A$ which moves $D$ across $B$ and then past $D\p$ removes at least 
$\nu$ from the intersection, thereby contradicting minimality. So, every 
component of $A\cap C\p$ is a spanning arc on $A$. Assume now that $\nu$ 
is such a component which is outermost on $C\p$, cutting off a disk $D\p$ 
with $\bd D\p=\nu\cup\nu\p$ for an arc $\nu\p$ in $F_0$. The result of 
performing a boundary compression on $A$ along $D\p$ is a disk $A\p$ in 
$F_0$. $\bd A\p=\bd D\pp$ for a disk $D\pp$ in $F_0$. The 2-sphere 
$A\p\cup D\pp$ bounds a 3-ball $B$ in $P\cap Y$. $B$ cannot contain $D\p$ 
as that would imply that $A$ is compressible in $Y$. It follows that the 
union of $B$ an an appropriate regular neighborhood of $D\p$ is a solid 
torus in $Y$ across which $A$ is \pl\ to an annulus in $\bd Y$. So we are 
done in this case. 

Now suppose that $A$ is \inc\ in $X$. Let $(\bd R)\times[0,1]$ be a 
collar on $\bd R$ in $R$ with $\bd R=(\bd R)\times\{0\}$. We may assume 
that $\gamma$ and $\delta$ meet this collar in product arcs. Let 
$Q=Y\cap(R-((\bd R)\times[0,1)))$. Then $Q$ is the exterior of a trefoil 
knot, and $Y\cap\bd Q$ is a proper annulus $A_Q$ in $Y$ which is isotopic 
to $F_1\cup O\cup D_2$. $A_Q$ splits $X$ into $Q$ and a space $X\p$ \hm\ 
to the Whitehead link exterior. $A_Q$ is \inc\ in both spaces. We may 
assume that $(\bd A)\cap(\bd Q)=\ns$. 

First suppose that $A\cap Q=\ns$. Then $A$ is an \inc\ annulus in $X\p$ 
and hence must be \pl\ in $X\p$ to an annulus $A\p$ in $\bd X\p$. 
Assume $A\p$ lies in $\bd K$. Since $A\cap(\bd X\p-\bd K)=\ns$ the solid 
torus in $X\p$ bounded by $A\cup A\p$ must lie in $Y$, and we are done. 
Assume $A\p$ lies in $\bd X\p-\bd K$. If $A\p$ contains the annulus 
$S(\omega_0,B_0)$, then the solid torus in $X\p$ bounded by $A\cup A\p$ 
must contain $X(\omega_0,B_0)$, which is impossible since this space contains 
$D_0\p$, which lies in $\bd K$. Thus $A\p$ does not contain $S(\omega_0,B_0)$. 
It follows that $A\p$ either must lie in $\bd Y\cap L$, in which case 
it $A$ is \bpl\ in $Y$, or must contain $A_Q$, in which case $A$ is 
\pl\ in $Y$ to $A_Q$. 

Now suppose that $A\cap Q\neq\ns$. Subject to the requirement that 
$(\bd A)\cap(\bd Q)=\ns$ put $A$ in general position with respect to 
$A_Q$ so that the intersection has a minimal number of components. Then 
each component is a simple closed curve which does not bound a disk 
on either surface. Suppose $A_1$ is a component of $A\cap Q$. Then 
$A_1$ is an \inc\ annulus in $Q$ with meridian boundary components. 
Since a trefoil knot is prime this implies that $A_1$ is \pl\ in $Q$ 
to an annulus in $\bd Q$ \cite{}. If that annulus lies in $A_Q$, 
then we can isotop $A$ to remove at least $\bd A_1$ from the intersection, 
thereby contradicting minimality. Thus each component $A_1$ of $A\cap Q$ 
is \pl\ in $Q$ to $\bd Q-\inte A_Q$. More precisely there is an embedding 
of $U_1=S^1\times[0,1]\times[0,1]$ in $Q$ with $S^1\times[0,1]\times\{0\}=
\bd Q-\inte A_Q$, $S^1\times[0,1]\times\{1\}=A_1$, and 
$S^1\times\{0,1\}\times[0,1]$ a collar on $\bd A_Q$ in $A_Q$. Let $A_0$ be 
the component of $A\cap X\p$ which contains $\theta\p$. The other 
component $\theta$ of $\bd A_0$ lies in $\inte A_Q$. We may assume that $A_1$ 
is the component of $A\cap Q$ containing $\theta$ and that 
$S^1\times\{0\}\times\{1\}=\theta$. Now $A_0$ is \bpl\ in $X\p$. 
More precisely there is a solid torus $U_0$ in $X\p$ with $\bd U_0=
A_0\cup A_0\p$, where $A_0\p$ lies in $\bd X\p$, $A_0\cap A_0\p=\bd A_0=
\bd A_0\p$, and these annuli are longitudinal in $U_0$. Let $A_Q\p=
A_0\p\cap A_Q$ and $A_Y\p=A_0\p\cap\bd Y$. Thus $A_0\p=A_Q\p\cup A_Y\p$. 

There are two possible configurations of $U_0$ and $U_1$ with respect to 
each other. In the first case $U_0\cap U_1=A_Q\p=S^1\times\{0\}\times[0,1]$. 
Then $U_0\cup U_1$ is a solid torus across which $A_0\cup A_1$ is \pl\ to 
$A_Y\p\cup(\bd Q-\inte A_Q)\cup(S^1\times\{1\}\times[0,1])$. We can then 
isotop $A$ so as to move $A_0\cup A_1$ across $U_0\cup U_1$ to 
$S^1\times\{1\}\times[0,1]$ and then past this annulus into $X\p$. This 
removes at least $\bd A_1$ from $A\cap A_Q$, thereby contradicting 
minimality. 

In the second case $U_0\cap U_1=\theta\cup(S^1\times\{1\}\times[0,1])$ 
and $A_Q\p$ properly contains $S^1\times\{1\}\times[0,1]$, i.e. $A_0$ 
is ``folded back over $A_Q$.'' Let $A_2$ be the component of $A\cap X\p$ 
other than $A_0$ which meets $A_1$. Let $\psi=A_1\cap A_2$ and 
$\psi\p=\bd A_2-\psi$. Now $A_2$ lies in $U_0$ and is \pl\ to both 
of the annuli into which $\bd A_2$ splits $\bd U_0$. If $\psi\p$ lies 
in $A_Y\p$, then as in the first case we can isotopy $A$ so as to move 
$A_1\cup A_2$ to $S^1\times\{0\}\times[0,1]$ and then past this annulus 
into $X\p$ to remove $\bd A_1$ from $A\cap A_Q$ (and in fact to conclude 
that $A$ is isotopic to $A_Q$.) If $\psi\p$ lies in $A_Q\p$, then $\psi\p=
A_2\cap A_3$ for a component $A_3$ of $A\cap Q$ which is parallel in $Q$ 
to $A_1$. We can then isotop $A$ to move $A_1\cup A_2\cup A_3$ to an annulus 
in $A_Q$ and then past that annulus into $X\p$ to remove at least 
$\bd A_1\cup\bd A_3$ from $A\cap A_Q$, again contradicting minimality. 

(7) Let $T$ be an \inc\ torus in $Y$. Let $Q$ and $A_Q$ be as in (6). 
Let $Y\p$ be the closure of $Y-Q$. $Y\p$ is a handlebody of genus two. 
There is a proper disk $D$ in $Y\p$ which meets $A_Q$ in a spanning arc. 
Put $T$ in general position with respect to $A_Q$ so that the intersection 
has a minimal number of components. As usual no component is a simple 
closed curve which bounds a disk on both surfaces.

Suppose $T\cap A_Q\neq\ns$. Let $A$ be a component of $T\cap Y\p$. 
Put $A$ in general position with respect to $D$ so that $A\cap D$ has 
a minimal number of components. This intersection must be non-empty. 
If there is a component which is \bpl\ on $A$, then by the irreducibility 
of $Y\p$ and the incompressibility of $A_Q$ there is an isotopy of 
an outermost disk on $A$ past the corresponding disk on $D$ which reduces 
the intersection, contradicting minimality. So every arc of the 
intersection is a spanning arc on $A$. There is such an arc $\nu$ of 
$A\cap D$ which cobounds an outermost disk $\Delta$ on $D$ with an arc 
$\nu\p$ in $D\cap A_Q$. Performing a boundary compression on $A$ along 
$\Delta$ gives a proper disk $A\p$ in $Y\p$ which cobounds a 3-ball $B$ 
with a disk in $A_Q$. $\Delta$ cannot lie in $B$ since this would make 
$A$ compressible in $Y\p$. It follows that the union of $B$ and an 
appropriate regular neighborhood of $\Delta$ is a solid torus across 
which $A$ is parallel to an annulus in $A_Q$. There is then an isotopy 
of $T$ in $Y$ which removes at least $\bd A$ from the intersection, 
contradicting minimality.  
Thus $T\cap A_Q=\ns$. $T$ cannot lie in the handlebody $Y\p$ since then 
it would be compressible. So $T$ must lie in $Q$. Since this is a trefoil 
knot exterior it must be \pl\ to $\bd Q$ \cite{Ja}, and we are done. \end{proof}

%Lemma 8.3
\begin{lem} Let $Y=X(\tau \cup \rho, R)$, and let $F_i=Y \cap D_i$ for 
$i=1,2$. Then 

(0) $Y$ is \irr, 

(1) $F_i$  is \inc\ in $Y$, 

(2) $O$ is \inc\ in $Y$, 

(3) $S(\tau\cup\rho, R)$ is \inc\ in $Y$, 

(4) $(Y,F_1\cup F_2, \bd Y-\inte(F_1\cup F_2))$ has the halfdisk property, 

(5) $(Y,F_1\cup F_2, \bd Y-\inte(F_1\cup F_2))$ has the band property,  

(6) every proper \inc\ annulus in $Y$ whose boundary lies in 
$(\inte(F_1\cup F_2)\cup(\bd Y-(F_1\cup F_2))$ is \pl\ to an annulus 
in $\bd Y$, and  

(7) $Y$ contains no \inc\ tori. 

\end{lem}

\begin{proof} The tangle $\tau$ is \textit{excellent} in the sense that 
$X(\tau,R)$ is \irr, \birr, anannular, and atoroidal \cite{My simple}. 
(This is not the same tangle as in \cite{My simple}, but the same proof 
works.) 

(0) $Y$ is \hm\ to the result of attaching a 1-handle to $X(\tau,R)$ and 
is therefore \irr. 

$F_i$ and $O$ are each \inc\ in $X(\tau,R)$ so they each must be 
\inc\ in the smaller space $X(\tau\cup\rho, R)$. This establishes 
(1) and (2).  

To prove (3), (4), and (5) it suffices to prove that 
$\bd Y-\inte(F_1\cup O)$ is \inc\ in $Y$. 
Suppose $D$ is a compressing disk for this surface. 
Let $E$ be the disk shown in Figure 2. Then $E\cap Y$ is a proper 
disk $E\p$ in $Y$ which splits $Y$ into a space $Y\p$ homeomorphic to 
$X(\tau, R)$. Let $E_1$ and $E_2$ be the two copies of $E\p$ in $\bd Y\p$ 
which are identified to obtain $E\p$ in $Y$. Let $F_1\p$ and $F_2\p$ be 
the annuli in $\bd Y\p$ into which $\bd E\p$ splits $F_1$ and $F_2$, 
respectively. Put $D$ in general position with respect to $E\p$.   
We assume that $D\cap E\p$ has a minimal number of components. 

First assume that $D\cap E\p=\ns$. Then $D$ lies in $Y\p$, 
$D\cap(E_1\cup E_2)=\ns$, and $\bd D=\bd D\p$ for a disk $D\p$ in $\bd Y\p$. 
If $D\p$ lies in $\bd Y$, then we are done. If $D\p$ does not lie in $\bd Y$, 
then $\inte D\p$ must contain $E_1$ or $E_2$ (or both). We may assume that 
it contains $E_1$. Since $E_1\cap F_1\p\neq\ns$ and $F_1\p$ is connected 
and disjoint from $\bd D\p$ we must have that $F_1\p$ lies in $\inte D\p$. 
It follows that a meridian of $\gamma$ bounds a disk in $\bd Y\p$, which 
is clearly not the case. 

Now assume that $D\cap E\p\neq\ns$. By minimality and irreducibility each 
component of the intersection is an arc. Let $\xi$ be an outermost such 
arc on $D$, and let $D_0$ be the outermost disk which it cuts off on $D$. 
We may assume that $\bd\xi$ lies in $\inte F_2$. We may regard $D_0$ as 
a disk in $Y\p$ with $\xi$ lying in, say, $E_1$. Then $\bd D_0=\bd D_0\p$ 
for a disk $D_0\p$ in $\bd Y\p$. Note that $F_1\p\cap \bd D_0\p=\ns$. 
If $F_1\p$ lies in $\inte D_0\p$, then a meridian of $\gamma$ again bounds a 
disk in $\bd Y\p$, which is not the case. Thus $D_0\p$ misses $F_1\p$. 
Hence the disk $E_0$ on $E$ bounded by the the union of $\xi$ and an 
arc in $\bd F_2\p$ must lie in $D_0\p$. Let $B$ be the 3-ball in $Y\p$ 
bounded by $D_0\cup D_0\p$. We may regard $B$ and $E_0$ as lying in $Y$. 
An isotopy of $D$ which moves $D_0$ across $B$ and then past $E_0$ 
removes at least $\xi$ from $D\cap E$ while keeping $\bd D$ in our 
surface, thereby contradicting minimality. 

Next we prove (6). Suppose $A$ is a proper, \inc\ annulus in $Y$ 
whose boundary misses $\bd(F_1\cup F_2)$. Since every simple closed 
curve in $F_1 \cup F_2$ is \bpl, we may isotop $A$ so that it misses 
$F_1\cup F_2$. Thus $\bd A$ lies in $S(\tau\cup\rho,R)\cup O$, which 
is \inc\ in $Y$. Put $A$ in general position with respect to $E\p$ 
so that $A\cap E\p$ has a minimal number of components. 

Assume that $A\cap E\p=\ns$. Then we may regard $A$ as lying in $Y\p$, 
and so it is \pl\ in $Y\p$ to an annulus $A\p$ in $\bd Y\p$. If $A\p$ 
lies in $\bd Y$, then we are done. If this is not the case, then $A\p$ 
must contain $E_1$ or $E_2$ (or both). We may assume that it contains 
$E_1$. Since $O-E_1$ is simply connected $\bd A\p$ must lie in the 
disk with two holes $S\p$ obtained by splitting $S(\tau\cup\rho,R)$ 
along its intersection with $\bd E\p$. Inspection of $\bd Y\p$ shows 
that there is no annulus $A\p$ satisfying all these requirements.  

Now assume that $A\cap E\p\neq\ns$. Since $A$ is \inc\ the components 
of the intersection are 
either \bpl\ arcs in $A$ or spanning arcs of $A$. 

Suppose there is a \bpl\ arc. Let $\xi$ be an outermost such arc on $A$ 
which cuts off an outermost disk $D_0$ from $A$. We may regard $D_0$ as 
a disk in $Y\p$ which meets $E_1$. Let $\xi\p=(\bd D_0)-\inte \xi$. 
Now $\xi\p$ must lie in either $S\p$ or $O\p$, where $O\p$ is the disk 
obtained by splitting $O$ along its intersection with $\bd E\p$. 
In either case $\bd D_0=\bd D_0\p$ for a disk $D_0\p$ in $\bd Y\p$ which 
lies in the union of $E_1$ and either $S\p$ or $O\p$. 
Let $E_0=D_0\p\cap E_1$. Let $B$ be the 3-ball in $Y\p$ bounded by 
$D_0\cup d_0\p$. We may regard $B$ as being embedded in $Y$. An 
isotopy of $D$ which moves $D_0$ across $B$ and then past $E_0$ removes 
at least $\xi$ from the intersection thereby contradicting minimality. 
Thus there are no arcs which are \bpl\ in $A$. 
 
Now suppose that there is a spanning arc $\xi$. We may assume that it 
is outermost on $E$. Let $E_0$ be an outermost disk on $E$ cut off by 
$\xi$. We may assume that $(\bd E_0)-\inte \xi$ lies in the surface 
$S(\tau\cup\rho)\cup F_2\cup O$, which is \inc\ in $Y$. 
Let $D$ be the disk obtained 
by boundary compressing $A$ along $E_0$. Then $\bd D=\bd D\p$ for a disk 
$D_0\p$ in $S(\tau\cup\rho)\cup F_2\cup O$. 
Let $B$ be the 3-ball in $Y$ bounded by 
$D_0\cup D_0\p$. Then the union of $B$ and an appropriate regular 
neighborhood of $E_0$ is a solid torus across which $A$ is \pl\ to 
an annulus in $\bd Y$. 

Finally we prove (7). Any \inc\ torus $T$ in $Y$ can be isotoped off the 
disk $E\p$, and so can be moved into a space \hm\ to $X(\tau,R)$. 
Since this space is atoroidal $T$ must be boundary parallel in it. But 
this is impossible since the space has no torus boundary component. 
\end{proof}

%Lemma 8.4
\begin{lem} Let $Y=X(\beta\cup\gamma\cup\delta\cup\epsilon\cup
\rho\cup\zeta,J)$ and $A=Y\cap D_2$. Then 

(1) $Y$ is \irr, 

(2) $F_0\cup A$ is \inc\ in $Y$, 

(3) $\bd Y-\inte(F_0\cup A)$ is \inc\ in $Y$, and 

(4) $(Y,F_0,\bd Y-\inte F_0)$ has the halfdisk property. 

\end{lem}

\begin{proof} (1) $Y$ is obtained from the space of Lemma 8.2, call it 
$\widehat{Y}$, by drilling out a regular neighborhood of the arc 
$(\rho\cup\zeta)\cap\widehat{Y}$. Any 2-sphere in $Y$ must bound a 
ball in $\widehat{Y}$; since $\bd Y$ is connected this ball must lie 
in $Y$. 

(2) $F_0$ is \inc\ in $\widehat{Y}$ and so must be \inc\ in the smaller 
space $Y$. $A$ is \inc\ in $Y$ for homological reasons. 

(3) Suppose $D$ is a compressing disk for $G\cup O$ in $Y$. Since 
$G\cup O\cup D_2$ is \inc\ in $\widehat{Y}$ $\bd D$ must bound a disk 
in $G\cup O\cup D_2$ which contains $D_2$. But this implies that $A$ 
is compressible in $Y$, contradicting (2). 

Suppose $D$ is a compressing disk for the other component $S$ of 
$\bd Y-\inte(F_0\cup A)$. It follows from Lemma 6.2 that $\bd D$ must 
be parallel in $S$ to $S\cap A$, which again contradicts (2). 

(4) Now suppose that $D$ is a proper disk in $Y$ which meets $F_0\cup A$ 
in an arc $\theta$ and $\bd Y-\inte(F_0\cup A)$ in an arc $\varphi$. 
We may assume that neither arc is \bpl. This implies that $\theta$ must 
lie in $F_0$. 

Assume that $\varphi$ lies in $G\cup O$. Then by Lemma 6.2 $\bd D=\bd D\p$ 
for a disk $D\p$ in $\bd\widehat{Y}$. If $D\p$ does not lie in $\bd Y$, 
then it must contain $D_2$. It follows that $A$ is compressible in $Y$, 
again contradicting (2). 

Assume that $\varphi$ lies in $S$. Let $\widehat{S}$ be the surface 
obtained by removing from $S$ a collar on $S\cap A$ in $S$ and replacing 
it by a meridinal disk for $\zeta$. We may assume that $\varphi$ lies in 
$\widehat{S}$ and misses the meridinal disk. Hence we may assume that it 
lies in $\bd\widehat{S}$. By Lemma 6.2 $\bd D=\bd D\p$ for a disk $D\p$ 
in $\bd\widehat{S}$. If $D\p$ does not lie in $\bd Y$, then it must contain 
the meridinal disk. It follows that $A$ is compressible in $Y$, contradicting 
(2) yet a final time. \end{proof}

Let $\bar{R}$ be a copy of $R$. Given any subset $\Sigma$ of $R$ denote 
the corresponding copy by $\bar{\Sigma}$. In the following two lemmas 
$\bar{R}$ is glued to $J$ by identifying $\bar{D}_2$ with $D_0$ so that 
$\bar{\alpha}\cap\bar{D}_2=\beta\cap D_0$ and 
$\bar{\zeta}\cap\bar{D}_2=\varepsilon\cap D_0$.  

%Lemma 8.5
\begin{lem} Let $Y=X(\bar{\delta}\cup\bar{\rho}\cup\bar{\zeta}\cup
\bar{\alpha}\cup\beta\cup\epsilon\cup\gamma\cup\delta\cup\rho\cup\zeta,
\bar{R}\cup J)$. Let $A_1=\bar{D_1}\cap Y$ and $A_2=D_2\cap Y$. Then 

(1) $Y$ is \irr, 

(2) $A_1\cup A_2$ is \inc\ in $Y$, 

(3) $\bd Y-\inte(A_1\cup A_2)$ is \inc\ in $Y$, and 

(4) $(Y,A_1\cup A_2,\bd Y-\inte(A_1\cup A_2))$ has the halfdisk property. 

\end{lem}

\begin{proof}  Let $\widehat{Y}=X(\bar{\zeta}\cup
\bar{\alpha}\cup\beta\cup\epsilon\cup\gamma\cup\delta, 
\bar{R}\cup J)$. This space is the exterior of a Whitehead 
link in which one component has acquired a local granny knot. 
It is therefore \irr\ and \birr\ \cite{My simple}. 
$Y$ is \hm\ to the space obtained from $\widehat{Y}$ by drilling 
out regular neighborhoods of the arcs 
$(\bar{\delta}\cup\bar{\rho})\cap Y$ and $(\rho\cup\zeta)\cap Y$. 

(1) Any two sphere in $Y$ must bound a 3-ball in $\widehat{Y}$. Since 
$\bd\widehat{Y}$ is connected this 3-ball must lie in $Y$. 

(2) $A_1\cup A_2$ is \inc\ for homological reasons. 

(3) Suppose $D$ is a compressing disk for $\bar{O}\cup G\cup O$. 
We may assume that $\bd D$ lies in $G$. 
Since $G$ is \inc\ in $X(\beta\cup\epsilon,L\cup H)$ 
we have that $D$ must intersect $F_0\cup F_1$. We may assume that 
$D$ is in minimal general position with respect to this surface. 
The boundary of an innermost disk on $D$ would be \bpl\ in $F_0$. 
This is homologically impossible. 

Now suppose that $D$ is a compressing disk for the other component $S$ 
of $\bd Y-\inte(A_1\cup A_2)$. We may assume that $\bd D$ lies in 
$\bd \widehat{Y}$. Then $\bd D=\bd D\p$ for a disk $D\p$ in 
$\bd \widehat{Y}$ which must contain one or both of the disks in 
$\bd\widehat{Y}$ whose boundaries are meridians of $\bar{\rho}$ and $\rho$. 
But this is homologically impossible. 

(4) Suppose $D$ is a proper disk in $Y$ which meets $A_1\cup A_2$ 
in an arc $\theta$ and $\bd Y-\inte(A_1\cup A_2)$ in an arc $\varphi$ 
such that $\theta\cup\varphi=\bd D$. The only way this can occur is 
for $\theta$ to be \bpl. The result then follows from (3). \end{proof}

%Lemma 8.6
\begin{lem} Let $Y=X(\bar{\delta}\cup\bar{\rho}\cup\bar{\zeta}\cup
\bar{\alpha}\cup\beta\cup\epsilon\cup\gamma,
\bar{R}\cup J)$. Let $A_1=\bar{D_1}\cap Y$. Then 

(1) $Y$ is \irr, 

(2) $A_1$ is \inc\ in $Y$, 

(3) $\bd Y-\inte A_1$ is \inc\ in $Y$, and 

(4) $(Y,A_1,\bd Y-\inte A_1)$ has the halfdisk property. 

\end{lem}

\begin{proof}  Let $\widehat{Y}=X(\bar{\zeta}\cup
\bar{\alpha}\cup\beta\cup\epsilon\cup\gamma\cup\delta, 
\bar{R}\cup J)$. This space is the exterior of a Whitehead 
link in which one component has acquired a local granny knot. 
It is therefore \irr\ and \birr\ \cite{My simple}. 
$Y$ is \hm\ to the space obtained from $\widehat{Y}$ by drilling 
out a regular neighborhood of the arc  
$(\bar{\delta}\cup\bar{\rho})\cap Y$. 

(1) Any 2-sphere in $Y$ must bound a 3-ball in $\widehat{Y}$. Since 
$\bd\widehat{Y}$ is connected this 3-ball must lie in $Y$. 

(2) Suppose $D$ is a compressing disk for $A_1$ in $Y$. 
Put $D$ in minimal general position with respect to $F_0$. 
An innnermost disk on $D$ cannot lie in $J$ by Lemma 8.4 and 
cannot lie in $\bar{R}$ for homological reasons. So $A_1$ is 
\inc\ in $Y$.  

(3) Suppose $D$ is a compressing disk for $\bar{O}\cup G\cup O\cup D_2$. 
We may assume that $\bd D$ lies in $G$. 
Put $D$ in minimal general position with respect to $F_0$. 
An innermost disk on $D$ with boundary in $D\cap F_0$ cannot lie in 
$J$ by Lemma 8.4 and cannot lie in $\bar{R}$ for homlogical reasons. 
Thus $D$ must lie in $J$ and so by Lemma 8.4 $\bd D=\bd D\p$ for 
a disk $D\p$ in $\bd(Y\cap J)$. For homological reasons $D\p$ cannot 
contain $F_0$, and so must lie in $G\cup D_2$. Thus $\bar{O}\cup G
\cup D_2$ is \inc\ in $Y$.  

Now suppose that $D$ is a compressing disk for the other component $S$ 
of $\bd Y-\inte A_1$. We may assume that $\bd D$ lies in 
$\bd \widehat{Y}$. Then $\bd D=\bd D\p$ for a disk $D\p$ in 
$\bd \widehat{Y}$ which must contain the disk in 
$\bd\widehat{Y}$ whose boundary is a meridian of $\bar{\rho}$. 
This contradicts (2), so $S$ is \inc\ in $Y$. 

(4) Suppose $D$ is a proper disk in $Y$ which meets $A_1$ 
in an arc $\theta$ and $\bd Y-\inte A_1$ in an arc $\varphi$ 
such that $\theta\cup\varphi=\bd D$. The only way this can occur is 
for $\theta$ to be \bpl. The result then follows from (3). \end{proof}

%Section 9
\section{The construction of $W$}

In this section we construct a Whitehead manifold $W$ which is an 
infinite cyclic covering space of another \tm\ $W^\#$. In the sections that 
follow it will be shown that $W$ is \rirr\ and that whenever it 
non-trivially covers a \tm\ the group of covering translations must be 
infinite cyclic. Later we will modify the construction to obtain 
uncountably many such examples as well as examples which cannot 
non-trivially cover any \tm. 

We will also define certain open subsets of $W$ which we will later 
show to be a complete set of representatives for the isotopy classes of 
end reductions of $W$. 

For each integer $n\geq 0$ take a copy of each of the objects defined 
in section 8. Denote the $n^{th}$ copy of $D_j$ by $D_{n,j}$, that of 
each of the other objects by a subscript $n$. We regard the arcs with 
subscript $n$ as embedded in the \tm s with subscript $n+1$. 

We embed $J^\#_n$ in $\inte J^\#_{n+1}$ as $N(\mu_n,J^\#_{n+1})$ in the 
following manner. $L_n$ is sent to $N(\la_n,J^\#_{n+1})$. We assume that 
this set is the image in $J^\#_{n+1}$ of $N(\ep_n,L_{n+1})\cup
N(\de_n\cup\ze_n,R_{n+1})$, that $\rho_n$ meets $N(\la_n,J^\#_{n+1})$ in 
$\rho_n-\inte \rho\p_n$, and that each of $\al_n$ and $\ga_n$ meet 
$N(\la_n,J^\#_{n+1})$ in an arc. Let $\al\p_n$ and $\ga\p_n$, respectively, 
be the intersections of $\al_n$ and $\ga_n$ with $J^\#_{n+1}-
\inte N(\la_n,J^\#_{n+1})$, and let $\eta\p_n=\al\p_n\cup\be_n\cup\ga\p_n$. 
Send $H_n$ to $N(\eta\p_n,J^\#_{n+1}-\inte(\la_n,J^\#_{n+1}))$. We assume 
that this set is the image in $J^\#_{n+1}$ of 
$N(\be_n,L_{n+1})\cup N(\al\p_n\cup\ga\p_n,R_{n+1}-
\inte N(\de_n\cup\ze_n))$. Thus $L_n\cup H_n$ is sent to 
$N(\ka_n,J^\#_{n+1})$. We assume that the intersection of this set with 
$R_{n+1}$ is $N(\tau_n,R_{n+1})$. Finally we send $R_n$ to 
$N(\rho\p_n,X(\tau_n,R_{n+1}))$. The result is shown in Figure 4. 

%INSERT FIGURE 4 HERE. 

\begin{figure}
\epsfig{file=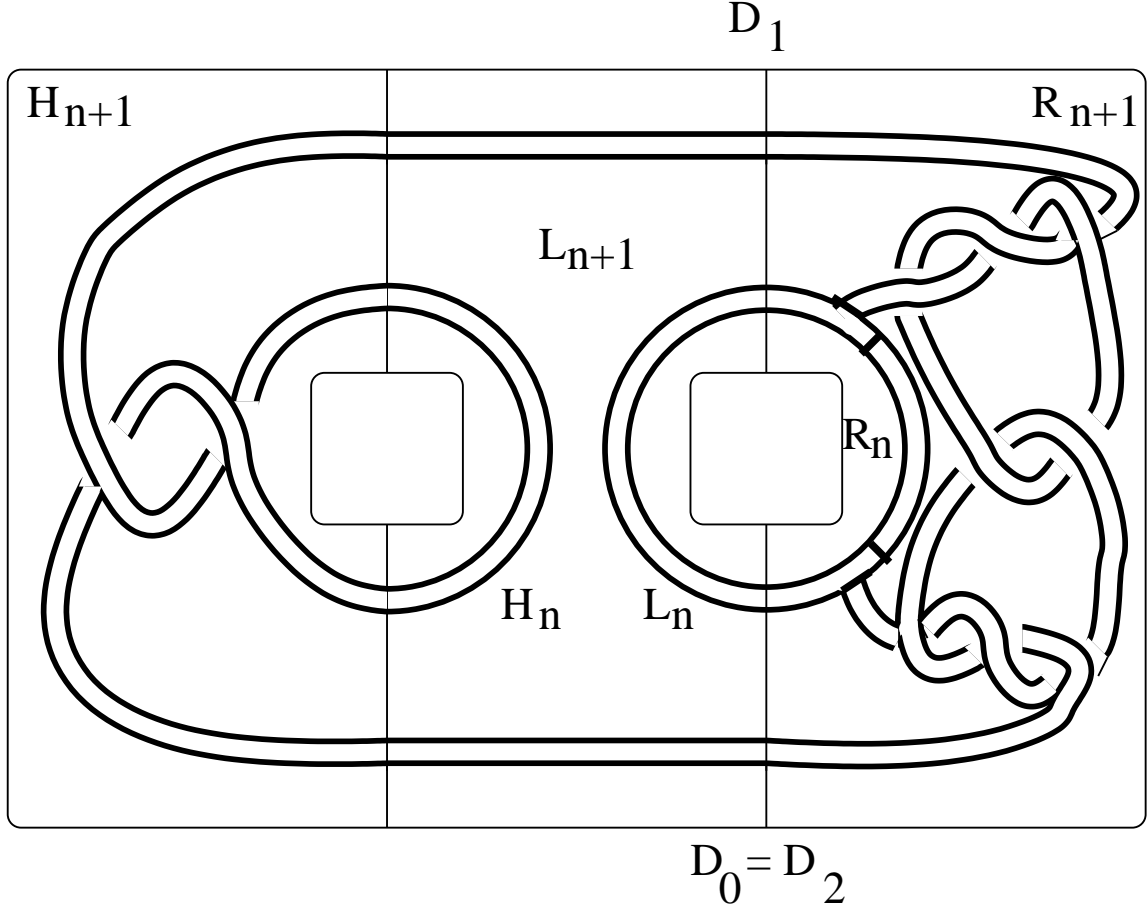,width=6in}
\caption{The embedding of $J_n^{\#}$ in $J_{n+1}^{\#}$}
\end{figure}

Now let $W^\#$ be the direct limit of the $J^\#_n$, and let 
$p:W\rightarrow W^\#$ be the universal covering map. Then $\pi_1(W^\#)$ 
is infinite cyclic. Let $h:W\rightarrow W$ be a generator of the group 
of covering translations. We regard $p\n(P^\#_n)$ as $D_n\times\R$ 
with $P_{n,j}=D_n\times [2j,2j+2]$, $L_{n,j}=D_n\times [2j,2j+1]$, 
and $R_{n,j}=D_n\times [2j+1,2j+2]$. We have that $p\n(H_n)$ is a 
disjoint union of 1-handles $H_{n,j}$, where $H_{n,j}$ is attached to 
$\bd D_n\times (2j,2j+1)$, thereby yielding a copy $J_{n,j}=
P_{n,j}\cup H_{n,j}$ of $J_n$. Set $D_{n,k}=D_n\times\{k\}$ for 
$k\in\mathbf{Z}$. For all the objects with subscript $n$ contained in 
$L_{n+1}\cup H_{n+1}$ or contained in $R_{n+1}$ denote the preimage 
contained in $L_{n+1,j}\cup H_{n+1,j}$ or in $R_{n+1,j}$ by using the 
subscripts $n,j$. We assume that $h$ is chosen so that 
$h(D_{n,k})=D_{n,k+2}$ and the image under $h$ of any object with subscripts 
$n,j$ has subscripts $n,j+1$. The embedding of $p\n(J^\#_n)$ in 
$p\n(J^\#_{n+1})$ is shown in Figure 5.  

%INSERT FIGURE 5 HERE. 

\begin{figure}
\epsfig{file=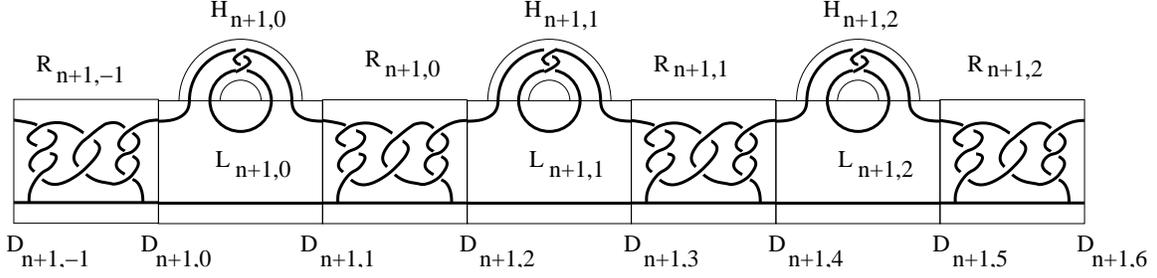, width=6in}
\caption{The embedding of $p^{-1}(J_n^{\#})$ in $p^{-1}(J_{n+1}^{\#})$}
\end{figure}

We now describe the open sets which, up to isotopy, will constitute all 
the end reductions of $W$. 

Let $\mathcal{P}=\{p_1, p_2, \cdots, p_m\}$ be a finite non-empty set 
of distinct integers ordered so that $p_1<p_2<\cdots<p_m$. We say that 
\PP\ is \textit{good} if $p_{i+1}=p_i+1$ for $1\leq i\leq m-1$. Otherwise 
\PP\ is \textit{bad}. Note that if $m=1$, then \PP\ is automatically good. 

For $n\geq0$ we let $C^{\PP}_n$ be the union of those $R_{n,j}$ with 
$p_1-1\leq j\leq p_m$, those $L_{n,j}$ with $p_1\leq j\leq p_m$, and those 
$H_{n,p}$ with $p\in\PP$. Each $C^{\PP}_n$ is a cube with $m$ handles. 
We have that $C^{\PP}_n\sbs \inte C^{\PP}_{n+1}$. This embedding is illustrated 
in Figure 6 for $\PP=\{0,1,3\}$. 

%INSERT FIGURE 6 HERE. 

\begin{figure}
\epsfig{file=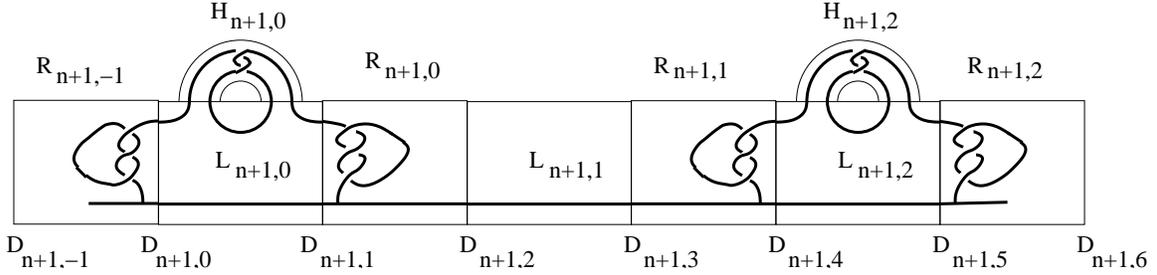, width=6in}
\caption{The embedding of $C^{\mathcal{P}}_n$ in $C^{\mathcal{P}}_{n+1}$ 
for $\mathcal{P}=\{0,2\}$}
\end{figure}

The sequence $\{C^{\PP}_n\}_{n\geq0}$ is denoted by $C^{\PP}$; 
its union is denoted by $V^{\PP}$. 
%We say that $V^{\PP}$ is 
%\textit{good} or \textit{bad} if \PP\ is good or bad, respectively. 
Whenever \PP\ is good and $m>1$ we denote $V^{\PP}$ by $V^{p,q}$, 
where $p=p_1$ and $q=p_m$. When $\PP=\{p\}$ we use the notation $V^p$. 
The expressions $C^{p,q}_n$, $C^p_n$, $C^{p,q}$, and $C^p$ are 
defined similarly.   

%Section 10
\section{$\RR$-irreducibility}

In this section we show that the manifold $W$ constructed in section 9 is 
\rirr. In the process we will also show that certain important open 
submanifolds of $W$ are \rirr. 

We begin by recalling a criterion due to Scott and Tucker \cite{ST}  
for a $\mathbf{P}^2$-\irr\ open \tm\ to be \rirr. A reformulation of this 
criterion was used in \cite{My r2p2}; we will need a slightly more general 
reformulation here. 

A proper plane $\Pi$ in an open \tm\ $U$ is \textit{homotopically trivial} 
if for every compact subset $C$ of $U$ the inclusion map of $\Pi$ into 
$U$ is end-properly homotopic to a map whose image is disjoint from $C$. 

%Lemma 10.1
\begin{lem} Let $U$ be an \irr, open \tm, and let $\Pi$ be a proper plane 
in $U$. If $\Pi$ is homotopically trivial, then $\Pi$ is trivial. \end{lem} 

\begin{proof} This is Lemma 4.1 of \cite{ST}. \end{proof}

%Lemma 10.2
\begin{lem} Let $U$ be an \irr, open \tm, and let $\{C_n\}_{n\geq 1}$ be 
a sequence of compact 3-dimensional submanifolds of $U$ such that 
$C_n\subseteq \inte C_{n+1}$ and 

(1) each $C_n$ is \irr,

(2) each $\bd C_n$ is \inc\ in $U-\inte C_n$, and 

(3) if $D$ is a proper disk in $C_{n+1}$ which is in general position 
with respect to $\bd C_n$ such that $\bd D$ is not null-homotopic in 
$\bd C_{n+1}$, then $D\cap \bd C_n$ has at least two components which 
are not null-homotopic in $\bd C_n$ and bound disjoint disks in $D$. 

Then any proper plane in $U$ can be end-properly homotoped off $C_n$ 
for any $n$. \end{lem}

\begin{proof} This is Lemma 4.2 of \cite{ST}. \end{proof}

The precise criterion we shall use is the following slightly stronger 
version of Lemma 3.3 of \cite{My r2p2}. Define a \tm\ $Y$ to be \textit{weakly 
anannular} if every proper \inc\ annulus in $Y$ has both its boundary 
components in the same component of $\bd Y$. 

%Lemma 10.3
\begin{lem} Let $U$ be a connected, \irr, open \tm. Suppose that for each 
compact subset $K$ of $U$ there is a sequence $\{C_n\}_{n\geq 1}$ of 
compact 3-dimensional submanifolds of $U$ such that $C_n \subseteq 
\inte C_{n+1}$ and 

(1) each $C_n$ is \irr, 

(2) each $\bd C_n$ is \inc\ in $U-\inte C_n$, 

(3) each $C_{n+1}-\inte C_n$ is \irr, \birr, and weakly anannular, and 

(4) $K\subseteq C_1$. 

Then $U$ is \rirr. \end{lem}

\begin{proof} Let $D$ be a proper disk in $C_{n+1}$ which is 
in general position with respect to $\bd C_n$ such that 
$\bd D$ is not null-homotopic in $\bd C_{n+1}$. If some component of 
$D\cap \bd C_n$ bounds a disk on $\bd C_n$ then there is such a disk whose 
interior misses $D$. Surgery on $D$ along this disk gives a disk $D\p$ 
whose intersection with $\bd C_n$ is contained in the intersection of 
$D$ with $\bd C_n$ but has fewer components. By repeating this procedure, 
if necessary, we may assume that no component of $D\cap \bd C_n$ is 
null-homotopic in $\bd C_n$. If $D\cap \bd C_n=\ns$, then $\bd C_{n+1}$ is 
compressible in $C_{n+1}-\inte C_n$, thereby contradicting the 
$\bd$-irreduciblity of $C_{n+1}-\inte C_n$. If $D\cap \bd C_n$ has only 
one component, then it and $\bd D$ together bound an \inc\ annulus in 
$C_{n+1}-\inte C_n$ joining $\bd C_{n+1}$ and $\bd C_n$, thereby 
contradicting the weak anannularity of $C_{n+1}-\inte C_n$. 
If $D\cap\bd C_n$ has more than one component and no two of them bound 
disjoint disks in $D$, then the components must be nested on $D$. 
Then $\bd D$ and an outermost such component on $D$ again together bound 
an \inc\ annulus in $C_{n+1}-\inte C_n$ which joins $\bd C_{n+1}$ and 
$\bd C_n$. Thus there must be two components which bound disjoint disks 
in $D$. Now apply Lemma 10.2 and then Lemma 10.1. \end{proof}

Let $C=\{C_n\}$ be a sequence of compact, connected, 3-dimensional 
submanifolds of an \irr, open \tm\ $U$ such that $C_n\sbs\inte C_{n+1}$ and 
$U-\inte C_n$ has no compact components. This will be called a 
\textit{quasi-exhaustion for $U$}. Note that a quasi-exhaustion for $U$ 
whose union is $U$ is an exhaustion for $U$. A quasi-exhaustion is 
\textit{good} if it satisfies condtions (1)--(3) of Lemma 10.3. Thus 
that lemma can be rephrased by saying that if every compact subset of $U$ 
is contained in the first term of a good \qe, then $U$ is \rirr. 

Recall the \qe\ $C^{\PP}$ defined in section 9, where $\PP=\{p_1, p_2, 
\cdots, p_m\}$. The embedding of $C^{\PP}_n$ in 
$C^{\PP}_{n+1}$ is shown in Figure 7 for the case $m=1$, $p=p_1$ and in 
Figure 8 for the case $\PP=\{0,1,2\}$.  

%INSERT FIGURE 7 HERE. 

\begin{figure}
\epsfig{file=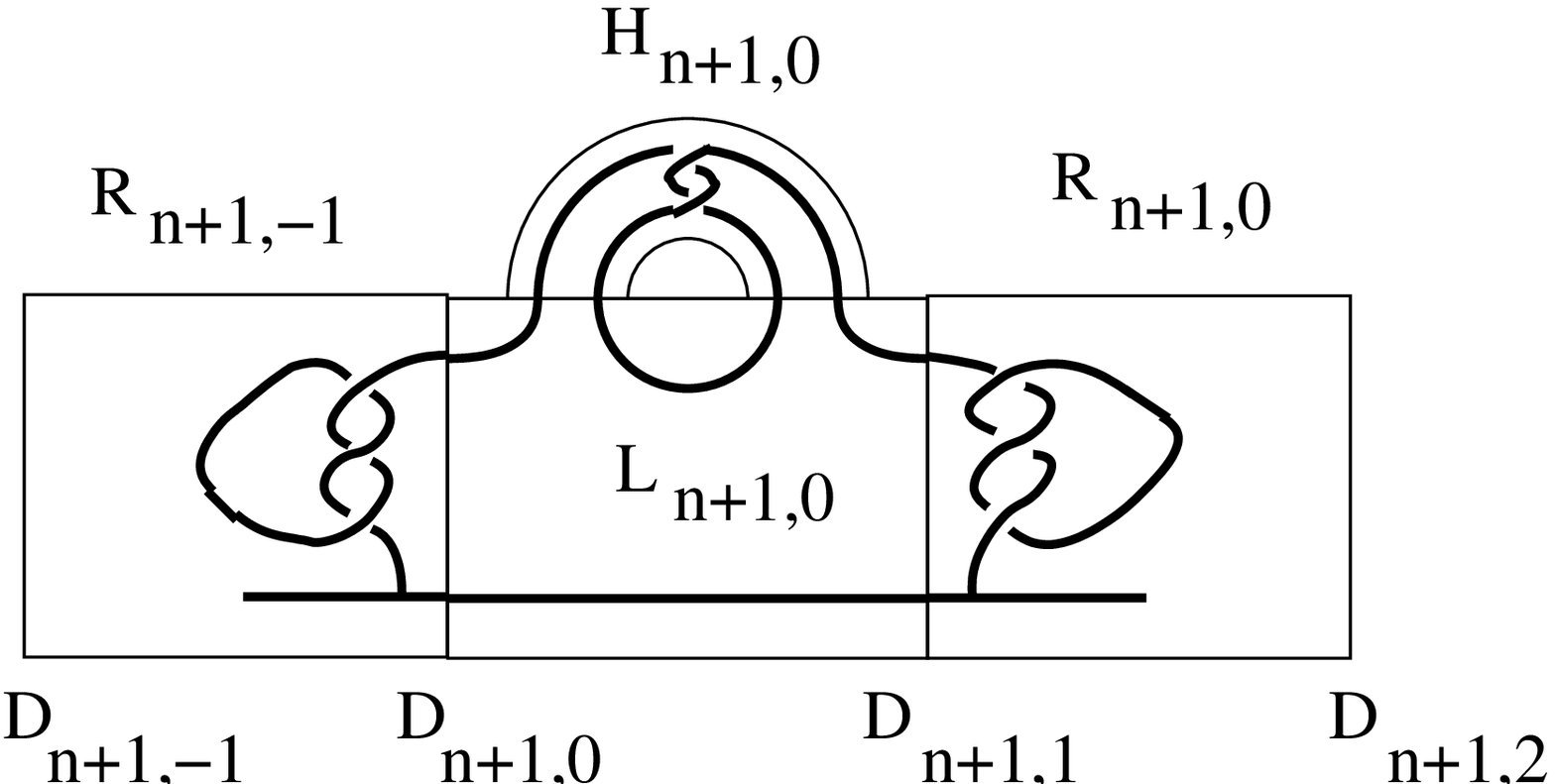, width=2.7in}
\caption{The embedding of $C^0_n$ in $C^0_{n+1}$} 
\end{figure}

%INSERT FIGURE 8 HERE. 

\begin{figure}
\epsfig{file=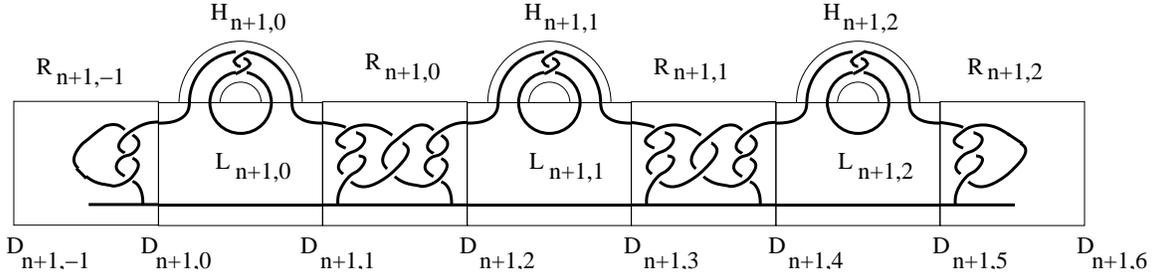, width=6in}
\caption{The embedding of $C^{[0,2]}_n$ in $C^{[0,2]}_{n+1}$} 
\end{figure}

%Proposition 10.4
\begin{prop} If \PP\ is good, then $C^{\PP}$ is good.  
\end{prop}

%Corollary 10.5
\begin{cor} $W$ is \rirr. Each $V^{\PP}$ with \PP\ good is \rirr. 
\end{cor}

\begin{proof}[Proof of Corollary 10.5] A given compact subset $K$ of $W$ 
lies in some $p\n(J^\#_n)$ and thus in a finite union of $J_{n,j}$ and hence 
in some $C^{p,q}_n$, which can be renumbered as the first term of the 
good \qe\ $C^{\PP}$. Thus by Lemma 10.3 $W$ is \rirr. 

$V^{\PP}$ is \rirr\ because it has the good exhaustion $C^{\PP}$. \end{proof} 

\begin{proof}[Proof of Proposition 10.4] 
We have that $C^m_n\sbs\inte C^m_{n+1}$, and 
$C^m_n\sbs C^{m+1}_n$. A given compact set $K$ of $W$ lies in some 
$p\n(J^\#_n)$ and thus in a finite union of $J_{n,j}$ and hence in some 
$C^m_n\sbs C^q_q$, where $q=\max\{m,n\}$. Thus $\{C^m_m\}$ is an 
exhaustion for $W$. 

$C^m_n$ is a cube with $2m+1$ handles and so is \irr. 
Let $Y=C^m_{n+1}-\inte C^m_n$. 

By Lemma 8.3 it suffices to prove the following. \end{proof}

%Lemma 8.6
\begin{lem} $Y$ is \irr, \birr, and weakly anannular. \end{lem}

\begin{proof} To simplify the notation we may assume that $p=0$. 

We first consider the case $m=1$. Then the embedding of $C^0_n$ into 
$C^0_{n+1}$ is shown in Figure 7. 

$Y$ is homeomorphic to the union of 
the Whitehead link exterior $Y_0$ and two trefoil knot exteriors $Y_1$ 
and $Y_2$. Thus each $Y_i$ is \irr\ and \birr. $Y_1\cap Y_2=\ns$. 
For $i=1$, $2$ we have that $Y_i\cap Y_0=\bd Y_i\cap \bd Y_0=A_i$, an 
annulus in the component $T$ of $\bd Y_0$ which is not equal to 
$\bd C_{n+1}$. Each $A_i$ is a regular neigborhood of a simple closed 
curve which is a meridian in both $Y_0$ and $Y_i$. Thus $A_i$ incompressible 
in both $Y_0$ and $Y_i$ and $T-\inte(A_1\cup A_2)$ is incompressible in 
$Y_0$. We apply Lemma 7.1 with $S=A_1\cup A_2$ and $Z=\bd Y$ to conclude 
that $Y$ is \irr\ and \birr. 

Every proper incompressible annulus in $Y_0$ is \bpl\ \cite{My simple}. 
For $i=1$, $2$ each proper \inc\ annulus in $Y_i$ with meridian 
boundary components is \bpl\ \cite{Ja}. We apply Lemma 7.2 with 
$S=A_1\cup A_2$ and $Z=\bd C_{n+1}$ to conclude that $Y$ is weakly 
anannular. 

We now consider the case $m>1$. The embedding of $C_n^{0,q}$ in  
$C_{n+1}^{0,q}$ is shown in Figure 8 for $q=2$. 

We split $Y$ along the surface $F=\cup_{i=1}^{2q}F_{n+1,i}$ into 
a \tm\ $Y\p$. Let $F\p$ be the surface in $\bd Y\p$ whose image in $Y$ 
is $F$. Let $Z=\bd Y$ and $Z_0=\bd C^m_n$. Let $Z\p$ and $Z_0\p$ be 
the surfaces into which $Z$ and $Z_0$, respectively, are split by their 
intersections with $F$. 

Each component of $Y\p$ is \hm\ to $X(\beta\cup\epsilon,L\cup H)$, 
$X(\beta\cup\gamma\cup\delta\cup\epsilon,J)$, or $X(\tau\cup\rho,R)$. 
It follows from Lemmas 8.1, 8.2, and 8.3 that $Y\p$ is \irr, $F\p$ and 
$Z\p$ are \inc\ in $Y\p$, and $(Y\p,F\p,Z\p)$ has the halfdisk property. 
So by Lemma 7.1 $Y$ is \irr\ and \birr. 

By Lemmas 8.1, 8.2, and 8.3 $(Y\p,F\p,Z\p)$ has the band property and 
every proper \inc\ annulus in $Y\p$ which misses $\bd F\p$ is either 
\bpl\ in $Y\p$ or has boundary which is \pl\ in $\bd Y\p$ to a pair of 
simple closed curves in $Z_0\p$. Therefore by Lemma 7.2 every proper 
\inc\ annulus in $Y$ is either \bpl\ in $Y$ or has its boundary in 
$Z_0$. In particular $Y$ is weakly annular. \end{proof} 

%Proposition 10.7
\begin{prop} If $\mathcal{P}$ is bad, then $V^{\mathcal{P}}$ is not 
\rirr. \end{prop}

\begin{proof} Suppose $\mathcal{P}=\{p_1,p_2, \ldots, p_m\}$. There 
is an $i$, $1\leq i\leq m-1$, such that $p_{i+1}>p_1+1$. Set 
$d=p_i+1$. Let $\Delta_0$ be a proper disk in $C^{\mathcal{P}}_0$ 
which lies in the interior of $L_{1,d}$ and whose boundary is a meridian 
of $\varepsilon_{1,d}$. There is a proper annulus $\Omega_1$ in 
$C^{\PP}_1-\inte C^{\PP}_0$ which lies in the interior of $L_{2,d}$ and 
whose boundary is the union of $\bd \Delta_0$ and a meridian of 
$\varepsilon_{1,d}$. Let $\Delta_1=\Delta_0\cup \Omega_1$. Continuing 
in this fashion we build for each $n>1$ a proper annulus $\Omega_n$ in 
$C^{\PP}_n-\inte C^{\PP}_{n-1}$ which lies in the interior of $L_{n+1,d}$ 
and whose boundary is the union of $\bd \Delta_{n-1}$ and a meridian of 
$\varepsilon_{n,d}$. Then $\Delta_n=\Delta_{n-1}\cup \Omega_n$  
is a proper disk in $C^{\PP}_n$ such that $\bd \Delta_n$ 
is a meridian of $\varepsilon_{n,d}$. The union $\Pi$ of the $\Delta_n$ is 
a plane which is proper in $V^{\PP}$. Its complement in $V^{\PP}$ has 
two components. One contains $V^{p_1}$; the other contains $V^{p_m}$. 
By \cite{Ha} none of the $V^p$ embeds in \RRR. It follows that $\Pi$ is 
non-trivial in $V^{\PP}$. \end{proof}

%11
\section{Genus}

Suppose $\{C_n\}$ is an exhaustion for a Whitehead manifold $V$. 
The \textit{genus of} $\{C_n\}$ is defined to be the maximum of the 
genera of the surfaces $\bd C_n$ if this set of numbers is bounded 
and to be $\infty$ if it is unbounded. The \textit{genus of} $V$ 
is defined to be the minimum of the genera of those exhaustions for $V$ 
which have finite genus if such exhaustions exist and to be $\infty$ 
if all exhaustions for $V$ have infinite genus. 

%Proposition 11.1
\begin{prop} $V^{\PP}$ has finite genus. It has genus one if and only 
if $\PP$ has exactly one element. \end{prop}

\begin{proof} Suppose $\PP=\{p_1,\ldots, p_m\}$. Then $V^{\PP}$ has 
an exhaustion by cubes with $m$ handles, so it has genus at most $m$. 
Since $V^{\PP}$ contains the manifold $V^{p_1}$ which does not embed 
in \RRR\ \cite{Ha} we have that $V^{\PP}$ is not \hm\ to \RRR\ and so has 
genus at least one. So if $m=1$, then $V^{\PP}$ has genus one. 

Suppose $m>1$. If $\PP$ is bad, then by Proposition 10.7 $V^{\PP}$ is 
not \rirr. Since every genus one Whitehead manifold is \rirr\ \cite{Kn} 
the genus of $V^{\PP}$ must be greater than one. 

Now suppose that $\PP$ is good. As in the $m>1$ case of the proof of 
Lemma 10.6 we let $Y=C^{[1,m]}_{n+1}-\inte C^{[1,m]}_n$, etc. By 
Lemmas 8.1, 8.2, and 8.3 all the hypotheses of Lemmas 7.1, 7.2, and 
7.3 are satisfied. Thus every \inc\ torus in $Y$ bounds a compact 
submanifold of $Y$ and every proper \inc\ annulus in $Y$ is either 
\bpl\ in $Y$ or cobounds a compact submanifold of $Y$ with an annulus 
in $\bd C^{[1,m]}_n$. 

Suppose $V^{\PP}$ has an exhaustion $\{K_n\}$ of genus one. We may 
assume that each $K_n$ has genus one, that each $\bd K_n$ is \inc\ in 
$V^{\PP}-\inte K_0$, and that $C_0\sbs \inte K_0$. Choose $q>0$ 
such that $K_1\sbs \inte C_q$. Let $T=\bd K_1$. Then $T$ is an 
\inc\ torus in $C_q-\inte C_0$ which separates $\bd C_q$ from $\bd C_0$. 
Hence $T$ cannot bound a compact submanifold of this space. 

On the other hand $C_q-\inte C_0$ is the union of spaces 
$Y_j=C_j-\inte C_{j-1}$ for $j=1,\ldots, q$ which satisfy 
the hypotheses of Lemma 7.3. Hence every \inc\ torus in this 
space must bound a compact submanifold of the space. This 
contradiction shows that $V^{\PP}$ must have genus greater 
than one. \end{proof}

%SECTION 12
\section{Every $V^{\PP}$ is an end reduction of $W$.}

Let $\PP=\{p_1,\ldots,p_m\}$. Recall that $V^{\PP}$ is defined 
as the union of the quasi-exhaustion $C^{\PP}$, where for each 
$n\geq 0$ $C^{\PP}_n$ is the union of those $R_{n,j}$ with 
$p_1-1\leq j\leq P_m$, those $L_{n,j}$ with $p_1\leq j\leq p_m$, 
and those $H_{n,p}$ with $p\in \PP$.

%Proposition 12.1
\begin{prop} For each $k\geq 0$ $V^{\PP}$ is an end reduction of 
$W$ at $C^{\PP}_k$. \end{prop}

\begin{proof} We may assume that $k=0$. We will show that $V^{\PP}$ 
is the constructed end reduction of $W$ associated to a certain 
exhaustion $\{K_n\}$ of $W$. For each $n\geq 0$ let $K_n^+$ be the 
union of those $R_{n,j}$ with $p_1-1-n\leq j\leq p_m+n$, and those 
$L_{n,j}$ and $H_{n,j}$ with $p_1-n\leq j\leq p_m+n$. Each $K_n^+$ is 
a cube with $1+2n+p_m-p_1$ handles. Let $K_n$ be $K_n^+$ for $n>0$ and 
$C^{\PP}_0$ for $n=0$. 

Note that for each $n>0$ we have that $C^{\PP}_n$ is obtained by compressing 
$\bd K_n$ in $W-C^{\PP}_0$. More specifically we remove those 1-handles 
$H_{n,j}$ of $K_n$ for which $j\notin \PP$. Note that the compressions are 
confined to $K_{n+1}-\inte C^{\PP}_{n-1}$. It suffices to show that 
each $\bd C^{\PP}_m$, $m>0$,  is \inc\ in $W-\inte C^{\PP}_0$. 

First note that for each $n\geq 0$ we have by Lemmas 7.1, 8.1, 8.2, 8.3, 
8.4, and 8.5 that $C^{\PP}_{n+1}-\inte C^{\PP}_n$ is \birr. It follows 
that for any $m>0$ we have that $C^{\PP}_m-\inte C^{\PP}_0$ is \birr\ and 
that for all $q>0$ we have that $C^{\PP}_{m+q}-\inte C^{\PP}_m$ 
is \birr. Now $K_{m+q}-\inte C^{\PP}_m$ is obtained by attaching 
1-handles to this space, so $\bd C^{\PP}_m$ is \inc\ in this new 
space. Since $K$ is an exhaustion for $W$ it follows that $\bd C^{\PP}_n$ 
is \inc\ in $W-\inte C^{\PP}_m$. Hence we have that $\bd C^{\PP}_m$ 
is \inc\ in $W-\inte C^{\PP}_0$. \end{proof}

%SECTION 13
\section{Every \rirr\ end reduction of $W$ is a $V^{\PP}$.}

%Proposition 13.1
\begin{prop} Every \rirr\ \er\ of $W$ is isotopic to some $V^{\PP}$. 
More precisely, if $J$ is a regular submanifold of $W$, $V$ is an 
\rirr\ \er\ of $W$ at $J$, and $J\sbs\inte C^{\QQ}_n$, then there is a 
good subset $\PP$ of $\QQ$ such that $V$ is isotopic in $W$ to 
$V^{\PP}$. If $V$ is already contained in $V^{\QQ}$ then the isotopy 
can be chosen to lie in $V^{\QQ}$. 
\end{prop}

\begin{proof} By Lemma 2.4 we may assume that $V$ is an end reduction 
of $W$ at a knot $\kappa\sbs \inte J$. By Theorem 5.1 we can isotop 
$V$ rel $\ka$ so that $V$ lies in $V^{\QQ}$ and is an \er\ of $V^{\QQ}$ 
at $\ka$. 

Let $\PP$ be a minimal subset of $\QQ$ such that $V$ can be isotoped 
in $V^{\QQ}$ so that it lies in $V^{\PP}$. 

If $\PP$ is bad, then $V^{\PP}$ contains a non-trivial plane $\Pi$ which 
splits it into components which are isotopic in $V^{\PP}$ 
to $V^{\RRRR}$ and $V^{\SSSS}$, 
where $\RRRR$ and $\SSSS$ form a non-trivial partition of $\PP$. 
By Theorem 6.2 $V$ can be isotoped in $V^{\PP}$ into one of these 
two components. This contradicts the minimality of $\PP$, and so 
$\PP$ must be good. 

Now $\ka$ lies in the interior of some $C_n^{\PP}$.  
We say that $\ka$ is \textit{$\PP$ busting} 
in $C^{\PP}_n$ if for each $p\in \PP$ the annulus 
$H_{n,p}\cap\bd C^{\PP}_n$ is \inc\ in $C^{\PP}_n-\kappa$. 

If $\ka$ is not $\PP$ busting in $C^{\PP}_n$, then for some $p\in\PP$ 
the annulus $H_{n,p}\cap\bd C^{\PP}_n$ is compressible in 
$C^{\PP}_n-\ka$. It follows that $\ka$ can be isotoped in $\inte C^{\PP}_n$ 
so that it lies in $C^{\PP-\{p\}}_n$. This implies by Theorem X that 
$V$ can be isotoped in $V^{\PP}$ so as to lie in $V^{\PP-\{p\}}$, 
contradicting the minimality of $\PP$. Thus $\ka$ is $\PP$ busting 
in $C^{\PP}_n$. 

Whenever $\ka$ lies in the interior of 
a handlebody $K$ we say that $\ka$ is \textit{disk busting} in $K$ if 
$\bd K$ is \inc\ in $K-\ka$. 

%Lemma 13.2
\begin{lem} $\kappa$ is disk busting in $C^{\PP}_{n+2}$. \end{lem}

The proof of this lemma appears below. We note that $\ka$ is not 
disk busting in $C^{\PP}_{n+1}$. By the proof of Proposition 10.1 
we know that $\bd C^{\PP}_{n+2}$ is \inc\ in $V^{\QQ}-\inte C^{\PP}_{n+2}$. 
Thus $\bd C^{\PP}_{n+2}$ is \inc\ in $V^{\QQ}-\kappa$. By Lemma 2.3 we have 
that $V$ is an end reduction of $V^{\QQ}$ at $C^{\PP}_{n+2}$ and therefore 
is isotopic to $V^{\PP}$. This concludes the proof of Proposition 11.1. 
\end{proof}

\begin{proof}[Proof of Lemma 13.2] We first consider the case $m=1$. 
Then $\kappa$ lies in the solid torus $C^p_n$. Any compressing disk 
for $\bd C^p_n$ must be isotopic to a compressing disk for 
$H_{n,p}\cap\bd C^p_n$. So $\bd C^p_n$ is \inc\ in $C^p_n-\kappa$. 
Since $C^p_{n+1}-\inte C^p_n$ and $C^p_{n+2}-\inte C^p_{n+1}$ are 
\birr\ we have that $\bd C^p_{n+2}$ is \inc\ in $C^p_{n+2}-\kappa$. 

We now assume that $m>1$. Isotop $\kappa$ in $C^{\PP}_n$ so that 
it is in minimal general position with respect to the union $\mathcal{D}$ 
of the disks $D_{n+2,i}\cap C^{\PP}_n$, where $2p_1+1\leq i\leq 2p_m$. 
Note that since $\kappa$ is $\PP$ busting in $C^{\PP}_n$ it must 
meet those 3-balls in $C^{\PP}_n$ which are regular neighborhoods of 
the arcs $\beta_{n,p}$ 

Recall that the configuration consisting of $C^{\PP}_{n+1}$ and the graph 
of which $C^{\PP}_n$ is a regular neighborhood is split by certain of 
the disks $D_{n+1,i}$ into a left hitch, junctions, eyebolts, 
and a right hitch as in the 
proof of Proposition 12.1. $C^{\PP}_n$ meets each of these configurations 
in a collection of 3-balls which contain 
proper subarcs of $\kappa$. Our stategy will be to apply Lemma 7.1 
to the case of $Y=X(\kappa, C^{\PP}_{n+2})$ with $S$ equal to the union 
of those $Y\cap D_{n+2,i}$ which arise in the above decomposition, and 
with $Z=\bd C^{\PP}_{n+2}$. The proof of Lemma 13.2 will then follow 
from the following sequence of lemmas. 
\end{proof}

%Lemma 13.3
\begin{lem}Let  $(L_{n+2,p}\cup H_{n+2,p},
\beta_{n+1,p}\cup\varepsilon_{n+1,p})$ be an eyebolt. Let 
$\kappa\p$ be the intersection of $\kappa$ with $L_{n+2,p}\cup H_{n+2,p}$. 
Let $X\p=X(\kappa\p,L_{n+2,p}\cup H_{n+2,p})$ and  
$F_{n+2,i}\p=X\p\cap D_{n+2,i}$, $i=2p,2p+1$. Then 

(1) $X\p$ is \irr, 

(2) $F_{n+2,2p}\p$, $F_{n+2,2p+1}\p$, and $G_{n+2,p}$ are \inc\ in $X\p$, and 

(3) $(X\p,F_{n+2,2p}\p\cup F_{n+2,2p+1}\p,G_{n+2,p})$ has the halfdisk 
property. \end{lem}

\begin{proof} Let $N_{n+1}=C^{\PP}_{n+1}\cap(L_{n+2,p}\cup H_{n+2,p})$. 
It has two components, namely the  regular neighborhoods $N_{n+1}^+$ of 
$\beta_{n+1,p}$ and $N_{n+1}^-$ of $\varepsilon_{n+1,p}$ in 
$L_{n+2,p}\cup H_{n+2,p}$. 
Let $X_{n+2}$ be the exterior in 
$L_{n+2,p}\cup H_{n+2,p}$ of the union of these arcs, so we have 
$L_{n+2,p}\cup H_{n+2,p}=X_{n+2}\cup N_{n+1}$. 
Let $A_{n+1}^{\pm}$ be the annuli $N_{n+1}^{\pm}\cap X_{n+2}$. 
By Lemma 6.1 we have that $X_{n+2}$ is \irr, 
$F_{n+2,2p}\cup F_{n+2,2p+1}$, $A_{n+1}^{\pm}$, and 
$G_{n+2,p}$ are \inc\ in $X_{n+2}$, and $(X_{n+2}, F_{n+2,2p}\cup 
F_{n+2,2p+1}, \bd X_{n+2}-\inte (F_{n+2,2p}\cup F_{n+2,2p+1}))$ 
has the halfdisk property. 

Next let $N_n=C^{\PP}_n\cap (L_{n+2,p}\cup H_{n+2,p})$. 
It has four components, namely the two components of the regular 
neighborhood $N_n^+$ of the union of two arcs 
$\beta_{n,p}\cap N_{n+1}^+$ in $N_{n+1}^+$, 
the regular neighborhood $N_n^-$ of $\beta_{n,p}\cap N_{n+1}^-$ in 
$N_{n+1}^-$, and the regular neighborhood 
$N_n^0$ of $\varepsilon_{n,p}\cap N_{n+1}^-$ in $N_{n+1}^-$. 
Let $X_{n+1}$ be the exterior in $N_{n+1}$ of the union 
of these four arcs. Let $A_n^+$ be the union of two annuli 
$N_n^+\cap X_{n+1}$, $A_n^-$ the annulus $N_n^-\cap X_{n+1}$,  
and $A_n^0$ the annulus $N_n^0\cap X_{n+1}$. 

$X_{n+1}$ has two components. The component 
$X_{n+1}^+$ lying in $N_{n+1}^+$ is the exterior of a Whitehead 
clasp. By \cite{My simple} we have that $X_{n+1}^+$ is \irr, 
$X_{n+1}^+\cap (D_{n+2,2p}\cup D_{n+2,2p+1})$, 
$A_{n+1}^+$ and $A_n^+$ are \inc\ in $X_{n+1}^+$, 
and $(X_{n+1}^+, X_{n+1}\cap (D_{n+2,2p}\cup D_{n+2,2p+1})), 
\bd N_{n+1}^+\cup \bd N_n^+)$ has the halfdisk property.  
The component $X_{n+1}^-$ lying in $N_{n+1}^-$ is the product 
of a disk with two holes and a closed interval. So it is \irr, 
$X_{n+1}^-\cap (D_{n+2,2p}\cup D_{n+2,2p+1})$, $A_{n+1}^-$, 
$A_n^-$, and $A_n^0$ are \inc\ in $X_{n+1}^-$, and 
$(X_{n+1}^-,X_{n+1}\cap(D_{n+2,2p}\cup D_{n+2,2p+1}), 
A_{n+1}^-\cup A_n^-\cup A_n^0)$ has the halfdisk property.

Let $N\p$ be a regular neighborhood of $\kappa\p$ in $N_n$. 
Let $X_n$ be the exterior of $\kappa\p$ in $N_n$.  
Since $\kappa\p$ consists of proper arcs in a collection of 
3-balls it follows that $X_n$ is \irr. Denote the 
components of $X_n$ contained in $N_n^+$, $N_n^-$, and $N_n^0$ by 
$X_n^+$, $X_n^-$, and $X_n^0$, respectively.

We have that $X\p=X_{n+2}\cup X_{n+1}\cup X_n$.  

Since $\kappa$ must join a left hitch in $C^{\PP}_n$ 
to a right hitch and $p\in\PP$, we must have that $\kappa\p$ meets 
each component of 
$N_n^+\cup N_n^-$. 
It follows from the minimality of $\kappa\cap\mathcal{D}$ that  
$A_n^+\cup A_n^-$ must be \inc\ in $X_n$. 
Now $\kappa\p$ may or may not meet $N_n^0$. If it does meet $N_n^0$, 
then as above we have that $A_n^0$  is \inc\ 
in $X_n$. In this case it follows that $X\p$ is \irr\ and that 
$G_{n+2,p}$ is \inc\ in $X\p$. If $\kappa\p$ does not meet $N_n^0$, 
then since $N_n^0$ is a product 3-ball in the product 3-ball $N_{n+1}^-$ 
we have that $A_n^-$ is parallel in $X\p$ to $A_{n+1}^-$. 
Hence we again get that $X\p$ is \irr\ and $G_{n+2}$ is \inc\ in $X\p$. 
Thus we have (1) and a part of (2). 

Suppose $D$ is a compressing disk for, say, $F_{n+2,2p}\p$ in $X\p$. 
Assume first that $\kappa\p$ meets $N_n^0$. 
Put $D$ in minimal general position with respect to 
$A_{n+1}^+\cup A_{n+1}^-\cup A_n^+\cup A_n^-\cup A_n^0$. 
There are no simple closed curve intersections. If $\nu$ is an 
outermost arc on $D$ cutting off an outermost disk $\Delta$, then 
the arc $\nu\p=\bd\Delta-\inte\nu$ is \bpl\ in the annulus 
containing it. Note that this uses the halfdisk property for 
the exterior of the Whitehead clasp $X_{n+1,p}^+$. 
It follows from the incompressibility of $F_{n+2,2p}\p$ 
in $X_{n+2}$, of $X_{n+1}\cap D_{n+2,2p}$ in $X_{n+1}$, and of 
$X_n\cap D_{n+2,2p}$ in $X_n$ that $\nu$ can be removed by an 
isotopy of $D$. Assume next that $\kappa\p$ misses $N_n^0$. 
Put $D$ in minimal general position with respect to 
$A_{n+1}^+\cup A_{n+1}^-\cup A_n^+\cup A_n^-$. 
Since $A_n^-$ and $A_{n+1}^-$ are \pl\ a similar argument 
shows that we can remove all intersections of $D$ with 
these annuli. Thus $D$ must lie in $X_{n+2}$, $X_{n+1}$, or 
$X_n$ and $\bd D$ must lie in the intersection of this manifold 
with $D_{n+2,2p}$, which is \inc. This completes the proof of (2). 

Now suppose that $D$ is a proper disk in $X\p$ with, say, 
$D\cap F_{n+2,2p}\p$ an arc $\theta$ and with $\bd D-\inte\theta$ 
an arc $\theta\p$ in $G_{n+2,p}$. As in the previous paragraph 
we put $D$ in minimal general position with respect to the 
appropriate collection of annuli and use minimality to remove 
all the intersections. $D$ must then lie in $X_{n+2}$, and we 
use Lemma 6.1 to complete the proof of (3). \end{proof}

%Lemma 13.4
\begin{lem} Let $(J_{n+2,p},\beta_{n+2,p}\cup\gamma_{n+2,p}\cup 
\delta_{n+2,p}\cup\varepsilon_{n+2,p})$ be a right hitch. 
Let $\kappa\p=\kappa\cap J_{n+2,p}$. Let $X\p=X(\kappa\p,J_{n+2,p})$. 
Let $F_{n+2,2j}\p=D_{n+2,2j}\cap X\p$. Then 

(1) $X\p$ is \irr, 

(2) $F_{n+2,2j}\p$ and $\bd J_{n+2,j}-\inte D_{n+2,2j}$ are \inc\ 
in $X\p$, and 

(3) $(X\p,F_{n+2,2j}\p,\bd J_{n+2,j}-\inte D_{n+2,2j})$ has the halfdisk 
property. \end{lem}

\begin{proof}Let $N_{n+1}=C^{\PP}_{n+1}\cap J_{n+2,p}$. It is a regular 
neighborhood of the arc $\beta_{n+1,p}\cup\gamma_{n+1,p}\cup
\delta_{n+1,p}$ in $J_{n+2,p}$. Note that it includes $R_{n+1,p}$.  
Let $X_{n+2}$ be the exterior in $J_{n+2,p}$ of this arc, so we 
have $J_{n+2,p}=X_{n+2}\cup N_{n+1}$. Let $A_{n+1}$ be the annulus 
$N_{n+1}\cap X_{n+2}$. By Lemma 8.2 we have that $X_{n+2}$ is \irr, 
$F_{n+2,2j}$, $\bd J_{n+2,p}-\inte D_{n+2,2j}$, and $A_{n+1}$ are 
\inc\ in $X_{n+2}$, and $(X_{n+2},F_{n+2,2j},\bd X_{n+2}-\inte F_{n+2,2j})$ 
has the halfdisk property. 

Next let $N_n=C^{\PP}_n\cap J_{n+2,p}$. It has two components $N_n^+$ and 
$N_N^-$. $N_n^+$ is a regular neighborhood in $N_{n+1}$ of the component 
of $\beta_{n,p}\cap J_{n+2,p}$ which lies entirely within $L_{n+2}\cup 
H_{n+2,p}$. $N_n^-$ is a regular neighborhood in $N_{n+1}$ of the 
other component of 
$(\beta_{n,p}\cup\gamma_{n,p}\cup\delta_{n,p}\cup\varepsilon_{n,j})\cap
J_{n+2,p}$. 
Let $X_{n+1}$ be the exterior of $N_n$ in $N_{n+1}$. Let $A_n^{\pm}$ 
be the annulus $N_n^{\pm}\cap X_{n+1}$. 
$N_n$ sits in $N_{n+1}$ as a regular neighborhood of a 
Whitehead clasp in which a trefoil knot has been tied in one component. 
From the properties of the unknotted Whitehead clasp \cite{My simple} it is 
easily checked that $X_{n+1}$ is \irr, $X_{n+1}\cap D_{n+2,2j}$, 
$A_{n+1}$ and $A_n^{\pm}$ are \inc\ in $X_{n+1}$, and 
$(X_{n+1},X_{n+1}\cap D_{n+2,2j}, A_{n+1}\cup A_n^+\cup A_n^-)$ has 
the halfdisk property. 

Let $N\p$ be a regular neighborhood of $\kappa\p$ in $N_n$. Let $X_n$ 
be the exterior of $\kappa\p$ in $N_n$. Then $X\p=X_{n+2}\cup X_{n+1}
\cup X_n$. Since $p\in\PP$ we must have that $\kappa\p$ meets both 
components of $N_n$. From the minimality of $\kappa\cap\mathcal{D}$ we 
have that $A_n^+\cup A_n^-$ is \inc\ in $X_n$. We now put a compressing 
disk or halfdisk in minimal general position with respect to 
$A_{n+1}\cup A_n^+\cup A_n^-$ and argue as in Lemma 13.3 to complete 
the proof. \end{proof} 

We warn the reader that the proof of the next lemma, which occupies the 
rest of this section, is a very 
lengthy checking of special cases, subcases, and subsubcases. 
It may be advisable to skip it on a first reading. 

%Lemma 13.5
\begin{lem} Let $(R_{n+2,p},\gamma_{n+1,p}\cup\delta_{n+1,p}\cup
\rho_{n+1,p}\cup\alpha_{n+1,p}\cup\zeta_{n+1,p})$ be a junction, 
where $p, p+1\in\PP$.  
Let $\kappa\p=\kappa\cap R_{n+2,p}$. Let $X\p=X(\kappa\p,R_{n+2,p})$. 
Let $F_{n+2,i}\p=D_{n+2,i}\cap X\p$, $i=2p+1,2p+2$. Then 

(1) $X\p$ is \irr, 

(2) $F_{n+2,2p+1}\p$, $F_{n+2,2p+2}\p$, and $O_p$ are \inc\ in $X\p$, and 

(3) $(X\p,F_{n+2,2p+1}\p\cup F_{n+2,2p+2}\p,O_p)$ has the halfdisk property. 
\end{lem}

\begin{proof} Let $N_{n+1}=C^{\PP}_{n+1}\cap R_{n+2,p}$. It is a 3-ball 
which is a regular  
neighborhood in $R_{n+2,p}$ of $\gamma_{n+1,p}\cup\delta_{n+1,p}\cup
\rho_{n+1,p}\cup\alpha_{n+1,p}\cup\zeta_{n+1,p}$. Let $X_{n+2}$ be the 
exterior in $R_{n+2,p}$ of this union of arcs, so we have $R_{n+2,p}=
X_{n+2}\cup N_{n+1}$. By Lemma 8.3 we have that $X_{n+2}$ is \irr, 
$F_{n+2,2p+1}$, $F_{n+2,2p+2}$, $O_p$, and the lateral surface 
$S_{n+1}=S(\gamma_{n+1,p}\cup\delta_{n+1,p}\cup
\rho_{n+1,p}\cup\alpha_{n+1,p}\cup\zeta_{n+1,p},R_{n+2,p}$ are 
\inc\ in $X_{n+2}$, and $(X_{n+2},F_{n+2,2p+1}\cup F_{n+2,2p+2}, 
\bd X_{n+2}-\inte(F_{n+2,2p+1}\cup F_{n+2,2p+2}))$ has the halfdisk 
property. 

Next let $N_n=C^{\PP}_n\cap R_{n+2,p}$. It has three components 
$N_n^-$, $N_n^+$, and $N_n^0$. $N_n^-$ is a regular neighborhood in 
$N_{n+1}$ of that component of $\beta_{n,p}\cap N_{n+1}$ which 
joins one component of $N_{n+1}\cap D_{n+2,2p+1}$ to the other. 
This arc is \bpl\ in $N_{n+1}$ via a disk which misses the other 
two components of $N_n$. Similarly $N_n^+$ is a regular neighborhood in 
$N_{n+1}$ of that component of $\beta_{n,p+1}\cap N_{n+1}$ which 
joins one component of $N_{n+1}\cap D_{n+2,2p+2}$ to the other. 
This arc is also \bpl\ in $N_{n+1}$ via a disk which misses the other 
two components of $N_n$. $N_n^0$ is a 3-ball which sits in $N_{n+1}$ 
in the same fashion that $N_{n+1}$ sits in $R_{n+2,p}$. Let $X_{n+1}$ 
be the exterior of $N_n$ in $N_{n+1}$. It is homeomorphic to the 
space obtained from $X_{n+2}$ by attaching a 1-handle to $F_{n+2,2p+1}$ 
and a 1-handle to $F_{n+2,2p+2}$. It is therefore \irr. We can also 
regard $X_{n+1}$ as being the space obtained by 
drilling out two \bpl\ arcs from a homeomorphic copy $\widehat{X}_{n+1}$ 
of $X_{n+2}$, where each arc has its endpoints in the surface 
corresponding to $F_{n+2,2p+1}$ or $F_{n+2,2p+2}$. Since $O_p$ is \inc\ in 
$X_{n+2}$, we have that $S_{n+1}$ is \inc\ in the smaller space $X_{n+1}$. 
Likewise the surface $S_n=N_n^0\cap X_{n+1}$ is \inc\ in $X_{n+1}$. 
The surface $X_{n+1}\cap (D_{n+2,2p+1}\cup D_{n+2,2p+2})$ has four 
components each of which is a disk with two holes; it is \inc\ in 
$X_{n+1}$ for homological reasons. Let $A_n^{\pm}=N_n^{\pm}\cap X_{n+1}$. 
Each of these annuli is \inc\ in $X_{n+1}$ for homological reasons. 
Suppose $D$ is a halfdisk in $X_{n+1}$ with one boundary arc $\theta$ in, 
say, $X_{n+1}\cap D_{n+2,2p+1}$. Then we may assume that the other 
boundary arc $\theta\p$ misses $A_n^-$. By Lemma 8.3 we have that 
$\bd D=\bd D\p$ for a disk $D\p$ in $\bd \widehat{X}_{n+1}$. Now 
there is a proper arc $\nu$ in $D\p$ which lies in $\bd S_{n+1}\cup \bd S_n$ 
and splits $D\p$ into a disk $\Delta$ with $\bd \Delta=\theta\cup\nu$ 
and a disk $\Delta\p$ with $\bd \Delta\p=\theta\p\cup\nu$. 
$\Delta$ lies in a component of $\widehat{X}_{n+1}\cap D_{n+2,2p+1}$, and 
$\Delta\p$ lies in either $S_n$ or $S_{n+1}$. $\Delta$ must in fact 
lie in $X_{n+1}\cap D_{n+2,2p+1}$, for otherwise we could find a 
compressing disk for $A_n^-$ in $X_{n+1}$. Thus 
$(X_{n+1},X_{n+1}\cap (D_{n+2,2p+1}\cup D_{n+2,2p+2}),S_n\cup S_{n+1}\cup 
A_n^-\cup A_n^+)$ has the halfdisk property. 

Let $N\p$ be a regular neighborhood of $\kappa\p$ in $N_n$. Let $X_n$ be 
the exterior of $\kappa\p$ in $N_n$. Since $\kappa\p$ consists of proper 
arcs in 3-balls we have that $X_n$ is \irr. Denote the components of 
$X_n$ contained in $N_n^-$, $N_n^+$, and $N_n^0$ by $X_n^-$, $X_n^+$, 
and $X_n^0$, respectively. 

Since $p,p+1\in\PP$ we must have that $\kappa\p$ meets all three 
components of $N_n$. It follows from the minimality of $\kappa\cap\mathcal{D}$ 
that $X_n\cap (D_{n+2,2p+1}\cup D_{n+2,2p+2})$ is \inc\ in $X_n$ and 
hence that $A_n^{\pm}$ is \inc\ in $X_n^{\pm}$, and that 
$(X_n^-,X_n^-\cap D_{n+2,2p+1},A_n^-)$ and 
$(X_n^+,X_n^+\cap D_{n+2,2p+2},A_n^+)$ have the halfdisk property. 
Note, however, that $S_n$ may or may not be \inc\ in $X_n^0$ 
and that $(X_n^0,X_n^0\cap(D_{n+2,2p+1}\cup D_{n+2,2p+2}),S_n)$ 
may or may not have the halfdisk property. 

First suppose that both of these properties are satisfied. Then we 
regard $X\p$ as $X_{n+2}\cup X_{n+1}\cup X_n$ and apply the techniques 
of the previous two lemmas to obtain the result. 

Now suppose that at least one of these properties is not satisfied. 
We will show how to reorganize $X\p$ as $X_{n+2}\p\cup X_{n+1}\p\cup 
X_n\p\cup X_n^-\cup X_n^+$ in such a way as to obtain the result. We 
will replace $N_n^0$ by a 3-manifold $N_n\p$ and let $X_n\p$ be the 
exterior of $\kappa\p\cap N_n\p$ in $N_n\p$. We will replace $N_{n+1}$ 
by a 3-manifold $N_{n+1}\p$ and let $X_{n+1}\p$ be the exterior of 
$N_n\p\cup N_n^-\cup N_n^+$ in $N_{n+1}\p$. Then we will let $X_{n+2}\p$ 
be the exterior of $N_{n+1}\p$ in $R_{n+2,p}$. There are several 
possibilities and possibly several steps in constructing these 
manifolds. 

\textit{Case 1:} $S_n$ is compressible in $X_n^0$ via a compressing 
disk $D$ such that $\bd D$ is not \bpl\ in $S_n$. 

The surface $S_n^1$ resulting from the compression consists of two 
annuli. Since $p,p+1\in\PP$ each of these annuli must join \dpo\ 
to \dpt. The compression splits $N_n^0$ into a 3-manifold $N_n^1$ 
which consists of two 3-balls each of which meets both \dpo\ and \dpt. 
Let $X_n^1$ be the exterior of $\kappa\p\cap N_n^1$ in $N_n^1$. 

\textit{Subcase (a):} $S_n^1$ is \inc\ in $X_n^1$. 

By the minimality of $\kappa\cap\mathcal{D}$ we have that 
$X_n^1\cap(\dpo\cup\dpt)$ is \inc\ in $X_n^1$. Since $S_n^1$ 
consists of annuli we then have that $(X_n^1,X_n^2\cap(\dpo\cup\dpt),
S_n^1)$ has the halfdisk property. 

We now have $N_n^1\cup N_n^-\cup N_n^+$ sitting in $N_{n+1}$ with 
exterior $X_{n+1}^1$. Clearly $X_{n+1}^1$ is \irr. 
$X_{n+1}^1\cap(\dpo\cup\dpt)$ is \inc\ in $X_{n+1}^1$ for homological 
reasons. Thus the boundary of any compressing disk for $S_{n+1}$ in 
$X_{n+1}^1$ would have to separate two of the components of $\bd S_{n+1}$ 
from the other two. Hence the compressing disk would have to separate 
two of the disks in $(N_n^1\cup N_n^-\cup N_n^+)\cap(\dpo\cup\dpt)$ 
from the other two. Since the components of this intersection which 
lie in \dpo\ are joined by $N_n^-$ and those which lie in \dpt\ are 
joined by $N_n^+$ we must have that the compressing disk separates 
the pair of disks in \dpo\ from the pair in \dpt. But this is impossible 
since each component of $N_n^1$ joins \dpo\ to \dpt. Thus $S_{n+1}^1$ is 
\inc\ in $X_{n+1}^1$. 

If $(X_{n+1}^1,X_{n+1}^1\cap(\dpo\cup\dpt),S_{n+1}\cup S_n^1\cup A_n^-\cup 
A_n^+)$ does not have the halfdisk property, then there is a proper disk 
$\Delta$ in $X_{n+1}^1$ which meets, say, \dpo\ in an arc $\theta$ and 
$S_{n+1}$ in an arc $\theta\p$ such that $\bd\Delta=\theta\cup\theta\p$, 
and neither arc is \bpl. $\Delta$ then separates $N_n^-$ from a 
component of $N_n^1$. It must therefore separate the two components of 
$N_n^1$ from each other. But this is impossible since they meet disks which 
are joined by $N_n^+$. Thus this triple has the halfdisk property. 

We now let $N_n\p=N_n^1$, $X_n\p=X_n^1$, $N_{n+1}\p=N_{n+1}^1$, and 
$X_{n+2}\p=X_{n+2}$. Then we apply the usual methods to complete the proof. 

\textit{Subcase (b):} $S_n^1$ is compressible in $X_n^1$. 

One component $S_n^2$ of $S_n^1$ must be \inc\ in $X_n^1$. Since 
$p,p+1\in\PP$ it must join a meridian of $\gamma_{n,p}$ to a meridian 
of $\alpha_{n,p}$. We let $N_n^2$ and $X_n^2$ be the components of $N_n^1$ and 
$X_n^1$ which meet $S_n^2$. By the arguments given in (a) we have that 
$X_n^2\cap(\dpo\cup\dpt)$ is \inc\ in $X_n^2$ and that $(X_n^2, X_n^2\cap 
(\dpo\cup\dpt), S_n^2)$ has the halfdisk property. 

We let $X_{n+1}^2$ be the exterior of $N_n^2\cup N_n^-\cup N_n^+$ 
in $N_{n+1}$. As in (a) we have that $X_{n+1}^2$ is \irr\ and that 
$S_{n+1}$ and $X_{n+1}^2\cap(\dpo\cup\dpt)$ are \inc\ in $X_{n+1}^2$. 

Note that $(X_{n+1}^2,X_{n+1}^2\cap(\dpo\cup\dpt),S_{n+1}\cup S_n^2\cup 
A_n^-\cup A_n^+)$ does nt have the halfdisk property. There are disjoint 
proper disks $\Delta_1^2$ and $\Delta_2^2$ in $X_{n+2}^2$ such that, 
for $i=1,2$, $\Delta_i^2$ meets $D_{n+2,2p+i}$ in an arc $\theta_i$ and 
$S_{n+1}$ in an arc $\theta_i\p$ such that 
$\bd \Delta_i^2=\theta_i\cup\theta_i\p$ and neither arc is \bpl. 
$\Delta_1^2\cup\Delta_2^2$ splits $N_{n+1}$ into a \tm\ $N_{n+1}^3$ 
consisting of three 3-balls $B_{n+1}^-$, $B_{n+1}^0$, and $B_{n+1}^+$. 
$B_{n+1}^{\pm}$ contains $N_n^{\pm}$. $B_{n+1}^0$ contains $N_n^2$. 
We let $S_{n+1}^3$ be the result of boundary compressing $S_{n+1}$ along 
$\Delta_1^2\cup\Delta_2^2$. It consists of three annuli. Let $Q_{n+1}^{\pm}$ 
be the exterior of $N_n^{\pm}$ in $B_{n+1}^{\pm}$. Let $Q_{n+1}^0$ be 
the exterior of $N_n^2$ in $B_{n+1}^0$. Let $X_{n+1}^3=Q_{n+1}^-\cup Q_{n+1}^+
\cup Q_{n+1}^0$. 

$X_{n+1}^3$ is \irr. $S_{n+1}^3$, $A_n^-$, $A_n^+$, and $S_n^2$ are \inc\ in 
$X_{n+1}^3$, as is $X_{n+1}^3\cap(\dpo\cup\dpt)$. The triple $(X_{n+1}^3, 
X_{n+1}^3\cap(\dpo\cup\dpt),S_{n+1}^3\cup A_n^-\cup A_n^+\cup S_n^2)$ 
has the halfdisk property. 

We now consider the exterior $X_{n+2}^3$ of $N_{n+1}^3$ in $R_{n+2,p}$. 
It is \irr. $O_p$ and $S_{n+1}^3$ are \inc\ in $X_{n+2}^3$ for 
homological reasons. The boundary of any compressing disk for 
$X_{n+2}\cap(\dpo\cup\dpt)$ in $X_{n+2}^3$ would, for homological reasons, 
bound a disk in, say, \dpo\ which contains $B_{n+1}^-\cap\dpo$ but 
does not contain $B_{n+1}^0\cap\dpo$. The compressing disk could be isotoped 
in the exterior of $B_{n+1}^-\cup B_{n+1}^+$ in $R_{n+2,p}$ to a compressing 
disk for $O_p$ in this space. But this space is \hm\ to the exterior of the 
tangle $\tau_{n+1,p}$ in $R_{n+2,p}$, which is \birr. Thus $X_{n+2}^3\cap
(\dpo\cup\dpt)$ is \inc\ in $X_{n+2}^3$. 

Suppose $(X_{n+2}^3,X_{n+2}^3\cap(\dpo\cup\dpt),O_p\cup S_{n+2}^3)$ 
does not have the halfdisk property. Then there is a proper disk $\Delta$ 
in $X_{n+2}^3$ with, say, $\Delta\cap \dpo$ an arc $\theta$ and 
$\Delta \cap S_{n+1}^3$ an arc $\theta\p$ such that $\theta\cup\theta\p=
\bd \Delta$ and neither arc is \bpl. $\theta\p$ must lie in the component 
$A_{n+1}^3$ of $S_{n+1}^3$ which joins \dpo\ to itself and must join the 
two components of $\bd A_{n+1}^3$. It follows that $\theta\cup\theta\p$ 
must be a trefoil knot in $R_{n+2,p}$ and so cannot bound a disk $\Delta$. 

Thus our triple has the halfdisk property, and so we set $X_i\p=X_i^3$ for 
$n\leq i\leq n+2$ and complete the proof as usual. 

\textit{Case 2:} $S_n$ is compressible in $X_n^0$ via a compressing disk 
$D$ such that $\bd D$ is \bpl\ in $S_n$. 

By Case 1 we may assume that there is no compressing disk for $S_n$ in 
$X_n^0$ whose boundary is not \bpl\ in $S_n$. The surface resulting from 
the compression along $D$ consists of a disk with two holes $S_n^1$ and 
a disk $S\p$. We may assume that $\bd S\p$ lies in \dpo. $D$ splits $N_n^0$ 
into a \tm\ consisting of 3-balls $N_n^1$ and $B\p$ such that $\kappa\p
\cap B\p=\ns$ and $S\p$ is \bpl\ across $B\p$. Since $p,p+1\in\PP$ we have 
that $\bd S\p$ is a meridian of $\zeta_{n,p}$. Let $X_n^1$ be the exterior 
of $\kappa\p\cap N_n^1$ in $N_n^1$. Because $p,p+1\in\PP$ any compressing 
disk for $S_n^1$ in $X_n^1$ would have boundary a meridian of $\zeta_{n,p}$. 
But then one could form a band sum along an arc in $S_n$ of this compressing 
disk and $D$ which would be a compressing disk for $S_n$ in $X_n^0$ whose 
boundary is not \bpl\ in $S_n$. Therefore $S_n^1$ is \inc\ in $X_n^1$. 

\textit{Subcase (a):} $(X_n^1, X_n^1\cap(\dpo\cup\dpt),S_n^1)$ has the 
halfdisk property. 

Let $X_{n+1}^1$ be the exterior of $N_n^1\cup N_n^-\cup N_n^+$ in $N_{n+1}$. 
$X_{n+1}^1\cup(\dpo\cup\dpt)$ and $S_n^1$ are \inc\ in $X_{n+1}^1$ for 
homological reasons. The boundary of a compressing disk for $S_{n+1}$ 
in $X_{n+1}^1$ would have to separate two components of $\bd S_{n+1}$ from 
the other two. Thus in $N_{n+1}$ the compressing disk would have to 
separate two disks of $N_{n+1}\cap(\dpo\cup\dpt)$ from the other two while 
missing $N_n^1\cup N_n^-\cup N_n^+$. This is impossible, and so $S_{n+1}$ 
is \inc\ in $X_{n+1}^1$. 

Note that $(X_{n+1}^1,X_{n+1}^1\cap(\dpo\cup\dpt),S_{n+1}\cup S_n^1
\cup A_n^-\cup A_n^+)$ does not have the halfdisk property. There is a 
proper disk $\Delta^1$ in $X_{n+1}^1$ such that $\Delta^1$ meets \dpo\ in 
an arc $\theta_1$ which lies in the non-annulus component of $X_{n+1}^1
\cap \dpo$. It meets $S_{n+1}$ in an arc $\theta_1\p$ such that 
$\bd \Delta^1=\theta_1\cup\theta_1\p$ and neither arc is \bpl. Moreover 
$\Delta^1$ splits $N_{n+1}$ into 3-bals $B_{n+1}^-$ and $B_{n+1}^+$ 
such that $B_{n+1}^-$ contains $N_n^-$ and $B_{n+1}^+$ contains 
$N_n^+\cup N_n^1$. Let $N_{n+1}^2=B_{n+1}^-\cup B_{n+1}^+$. We let 
$S_{n+1}^2$ be the result of boundary compressing $S_{n+1}$ along 
$\Delta^1$. Let $Q_{n+1}^-$ be the exterior of $N-n^-$ in $B_{n+1}^-$. 
Let $Q_{n+1}^+$ be the exterior of $N_n^+\cup N_n^1$ in $B_{n+1}^+$. 
Let $X_{n+1}^2=Q_{n+1}^-\cup Q_{n+1}^+$. It is \irr. 

$S_{n+1}^2\cap Q_{n+1}^-$ is an annulus $A_{n+1}^2$ which is \pl\ across 
$Q_{n+1}^-$ to $A_n^-$, and so $A_{n+1}^2$, $A_n^-$, and $Q_{n+1}^-
\cap\dpo$ are all \inc\ in $Q_{n+1}^-$, and $(Q_{n+1}^-,Q_{n+1}^-
\cap\dpo, A_n^-\cup A_{n+1}^2)$ has the halfdisk property. 

$S_{n+1}^2\cap Q_{n+1}^+$ is a disk with two holes $\Sigma_{n+1}^2$. 
$S_n^1$, $\Sigma_{n+1}^2$, and $Q_{n+1}^+\cap(\dpo\cup\dpt)$ are \inc\ in 
$Q_{n+1}^+$ for homological reasons. 

Suppose $(Q_{n+1}^+,Q_{n+1}^+\cap(\dpo\cup\dpt),\Sigma_{n+1}^2\cup S_n^1
\cup A_m^+)$ does not have the halfdisk property. Then there is a proper 
disk $\Delta$ in $Q_{n+1}^+$ such that $\Delta$ meets \dpt\ in an arc 
$\theta$ and $\Sigma_{n+1}^2 \cup S_n^1$ in an arc $\theta\p$ such that 
$\bd\Delta=\theta\cup\theta\p$ and neither arc is \bpl. First assume 
that $\theta\p$ lies in $\Sigma_{n+1}^2$. Then $\Delta$ separates 
$B_{n+1}^+\cap\dpo$ from one of the components of $B_{n+1}^+\cap\dpt$ 
in $B_{n+1}^+$ while missing $N_n^1\cap N_n^+$. This is impossible since 
$N_n^1$ meets all three components of $B_{n+1}\cap(\dpo\cup\dpt)$. 
So $\theta\p$ lies in $S_n^1$. $\bd\theta\p$ must lie in a single 
component of $\bd S_n^1$ and separate the other two components from 
each other. $\theta$ must lie in \dpt. Then in the exterior of $N_n^1$ in 
$B_{n+1}^+$ we can isotop $\theta\cup\theta\p$ to a non-trivial simple 
closed curve on $S_n^1$. This contradicts the fact that $S_n^1$ is 
\inc\ in this space for homological reasons. Thus our triple must have 
the halfdisk property. 

We now consider the exterior $X_{n+2}^2$ of $N_{n+1}^2$ in $R_{n+2,p}$. 
It is \irr. $O_p$ and $S_{n+1}^2$ are \inc\ in $X_{n+2}^2$ for 
homological reasons. The boundary of a compressing disk for $X_{n+2}^2
\cap(\dpo\cup\dpt)$ in $X_{n+2}^2$ would, for homological reasons, bound 
a disk in \dpo\ which contains $B_{n+1}^-\cap\dpo$ but does not contain 
$B_{n+1}^+\cap\dpo$. The compressing disk could be isotoped in the 
exterior of $N_n^-\cup N_n^+$ in $R_{n+2,p}$ to give a compressing disk 
for $O_p$ in this space. However this space is \hm\ to the exterior of 
the tangle $\tau_{n+1,p}$ in $R_{n+2,p}$, which is \birr\ \cite{}. Thus 
$X_{n+2}^2\cap(\dpo\cup\dpt)$ is \inc\ in $X_{n+2}^2$. 

Suppose $(X_{n+2}^2,X_{n+2}^2\cap(\dpo\cup\dpt), O_p\cup S_{n+1}^2)$ 
does not have the halfdisk property. Then there is a proper disk $\Delta$ 
in $X_{n+2}^2$ with $\Delta\cap(\dpo\cup\dpt)$ an arc $\theta$ and 
$\Delta\cap S_{n+1}^2$ an arc $\theta\p$ with $\bd\Delta=\theta\cup\theta\p$ 
and neither arc is \bpl. Recall that $S_{n+1}^2=A_{n+1}^2\cup\Sigma_{n+1}^2$, 
where $A_{n+1}^2$ is an annulus with $\bd A_{n+1}^2$ in \dpo\ and 
$\Sigma_{n+1}^2$ is a disk with two holes having one boundary component in 
\dpo\ and the other two in \dpt. 

Assume $\theta$ lies in \dpo. Suppose $\theta\p$ lies in $A_{n+1}^2$. 
Then $\theta\cup\theta\p$ is a trefoil knot in $R_{n+2,p}$ and so cannot 
bound a disk $\Delta$. So $\theta\p$ must lie in $\Sigma_{n+1}^2$ 
and separate the components of $\Sigma_{n+1}^2\cap\dpt$. But this implies 
that a meridian of the arc $\alpha_{n+1,p}$ is compressible in the exterior 
of the arc $\alpha_{n+1,p}\cup\zeta_{n+1,p}$ in $R_{n+2,p}$, which is 
homologically impossible. 

So $\theta$ must lie in \dpt\ and $\theta\p$ in $\Sigma_{n+1}^2$. 
Suppose $\theta\p$ joins different components of $\bd\Sigma_{n+1}^2$. 
then $\theta\cup\theta\p$ is a trefoil knot in $R_{n+2,p}$ and so 
cannot bound a disk $\Delta$. So $\theta\p$ must join a component of 
$\bd \Sigma_{n+1}^2$ to itself and separate the other two boundary 
components. But this implies that either a meridian of the arc 
$\alpha_{n+1,p}$ is compressible in the exterior of the arc 
$\alpha_{n+1,p}\cup\rho_{n+1,p}\cup\delta_{n+1,p}$ in $R_{n+2,p}$ or 
that a meridian of the arc $\zeta_{n+1,p}$ is compressible in the exterior 
of the arc $\zeta_{n+1,p}\cup\rho_{n+1,p}\cup\delta_{n+1,p}$ in $R_{n+2,p}$, 
both of which are homologically impossible. 

Thus our triple has the halfdisk property. Letting $X_i\p=X_i^2$ we complete 
the proof in the usual manner. 

\textit{Subcase (b):} $(X_n^1,X_n^1\cap(\dpo\cup\dpt),S_n^1)$ does not have 
the halfdisk property. 

There is a proper disk $\Delta^1$ in $X_n^1$ such that $\Delta^1$ meets 
$\dpo\cup\dpt$ in an arc $\theta_1$ and $S_n^1$ in an arc $\theta_1\p$. 
$\bd \theta_1\p$ must lie in one component of $\bd S_n^1$ and separate 
the other two components. Let $N_n^2$ be the \tm\ obtained by splitting 
$N_n^1$ along $\Delta^1$. Let $S_n^2$ be the result of performing the 
corresponding boundary compression on $S_n^1$ along $\Delta^1$. 
$S_n^2$ consists of a pair of annuli and $N_n^2$ a pair of 3-balls. 
There are three possibilities for how $N_n^2$ joins the components of 
$N_n^1\cap(\dpo\cup\dpt)$, depending on which component of this set meets 
both components of $N_n^2$. 

Let $X_n^2$ be the exterior of $\kappa\p\cap N_n^2$ in $N_n^2$. 
$X_n^2\cap(\dpo\cup\dpt)$ and $S_n^2$ are \inc\ in $X_n^2$. Since $S_n^2$ 
consists of annuli $(X_n^2,X_n^2\cap(\dpo\cup\dpt),S_n^2)$ has the 
halfdisk property. 

Let $X_{n+1}^2$ be the exterior of $N_n^2\cup N_n^-\cup N_n^+$ in $N_{n+1}$. 
$X_{n+1}^2\cap(\dpo\cup\dpt)$ and $S_n^2$ are \inc\ in $X_{n+1}^2$ for 
homological reasons. The boundary of a compressing disk for $S_{n+1}$ in 
$X_{n+1}^2$ would have to separate two components of $\bd S_{n+1}$ from 
the other two. Thus in $N_{n+1}$ the compressing disk would have to 
separate two disks of $N_{n+1}\cap(\dpo\cup\dpt)$ from the other two 
while missing $N_n^2\cup N_n^-\cup N_n^+$. This is impossible, and so 
$S_{n+1}$ is \inc\ in $X_{n+1}^2$. 

Note that $(X_{n+1}^2,X_{n+1}^2\cap(\dpo\cup\dpt),S_{n+1}\cup S_n^2
\cup A_n^-\cup A_n^+)$ does not have the halfdisk property. There is 
a proper disk $\Delta^2$ in $X_{n+1}^2$ such that $\Delta^2$ meets 
\dpo\ in an arc $\theta_2$ which lies in the non-annulus component of 
$X_{n+1}^2\cap\dpo$. It meets $S_{n+1}$ in an arc $\theta_2\p$ such that 
$\bd \Delta^2=\theta_2\cup\theta_2\p$ and neither arc is \bpl. Moreover 
$\Delta^2$ splits $N_{n+1}$ into 3-balls $B_{n+1}^-$ and $B_{n+1}^+$ 
such that $B_{n+1}^-$ contains $N_n^-$ and $B_{n+1}^+$ contains 
$N_n^=\cup N_n^2$. Let $N_{n+1}^3=B_{n+1}^-\cup B_{n+1}^+$. We let 
$S_{n+1}^3$ be the result of boundary compressing $S_{n+1}$ along $\Delta^2$. 
Let $Q_{n+1}^-$ be the exterior of $N_n^-$ in $B_{n+1}^-$. Let $Q_{n+1}^+$ 
be the exterior of $N_n^+\cup N_n^2$ in $B_{n+1}^+$. Let $X_{n+1}^3=
Q_{n+1}^-\cup Q_{n+1}^+$. It is \irr. 

$S_{n+1}^3\cap Q_{n+1}^-$ is an annulus $A_{n+1}^3$ which is \pl\ across 
$Q_{n+1}^-$ to $A_n^-$, and so $A_{n+1}^3$, $A_n^-$, and $Q_{n+1}^-\cap\dpo$ 
are all \inc\ in $Q_{n+1}^-$, and $(Q_{n+1}^-,Q_{n+1}^-\cap\dpo,
A_n^-\cup A_{n+1}^3))$ has the halfdisk property. 

$S_{n+1}^3\cap Q_{n+1}^+$ is a disk with two holes $\Sigma_{n+1}^3$. 
$S_n^2$, $\Sigma_{n+1}^3$, and $Q_{n+1}^+\cap(\dpo\cup\dpt)$ are \inc\ in 
$Q_{n+1}^+$ for homological reasons. 

Suppose $(Q_{n+1}^+,Q_{n+1}^+\cap(\dpo\cup\dpt),\Sigma_{n+1}^3\cup S_n^2
\cup A_n^+)$ does not have the halfdisk property. Then there is a proper 
disk $\Delta$ in $Q_{n+1}^+$ such that $\Delta$ meets $\Sigma_{n+1}^3$ in an 
arc $\theta\p$ such that $\bd\theta\p$ lies in a single component of $\bd 
\Sigma_{n+1}^3$ and $\theta\p$ separates the other two components. 
$\Delta$ meets $Q_{n+1}^+\cap(\dpo\cup\dpt)$ in an arc $\theta$ such that 
$\bd\Delta=\theta\cup\theta\p$. Neither arc is \bpl. $\Delta$ separates 
the two disks of $B_{n+1}^+\cap(\dpo\cup\dpt)$ which do not contain $\theta$ 
from each other and does not meet $N_n^+\cup N_n^2$. It follows that one 
component $C^+$ of $N_n^2$ joins the two components of $B_{n+1}^+\cap\dpt$ 
while the other component $C$ of $N_n^2$ joins $B_{n+1}^+\cap\dpo$ to a 
component of $B_{n+1}^+\cap\dpt$. Let $\widehat{Q}_{n+1}^+=
Q_n^+\cup N_n^+$. Then $C\cup C^+$ is a regular neighborhood of a tangle 
$\widehat{\tau}=\sigma\cup\sigma^+$ in $B_{n+1}^+$. 

There are two possibilities for $\widehat{\tau}$. 

If $\theta\p$ lies in the component of $B_{n+1}^+\cap\dpt$ which meets 
$\alpha_{n,p}$, then $\sigma$ is isotopic to $\gamma_{n,p}\cup\rho_{n,p}
\cup\alpha_{n,p}$ and $\sigma^+$ is isotopic to a parallel copy of 
$\alpha_{n,p}\cup\zeta_{n,p}$. In this case one can isotop $\sigma^+$ by 
moving one of its endpoints along $\sigma$ so as to obtain the graph 
$\gamma_{n,p}\cup\rho_{n,p}\cup\alpha_{n,p}\cup\zeta_{n,p}$. By isotoping 
$\rho_{n,p}$ so that $\rho_{n,p}\cap\zeta_{n,p}$ is moved along $\zeta_{n,p}$ 
to $\bd B_{n+1}^+$ and then off $\zeta_{n,p}$ we obtain a tangle in 
$B_{n+1}^+$ which is equivalent to the tangle $\tau_{n,p}$ in $R_{n+1,p}$. 
Thus $\widehat{Q}_{n+1}^+$ is \hm\ to the exterior of $\tau_{n,p}$ in 
$R_{n+1,p}$ and so is \birr. Thus $\bd\Delta=\bd\Delta\p$ for a disk 
$\Delta\p$ in $\bd \widehat{Q}_{n+1}^+$. If $\Delta\p$ does not lie in 
$\bd Q_{n+1}^+$, then for homological reasons it must contain both 
components of $N_n^+\cap \dpt$. But this is impossible since then $\Delta\p$ 
would also contain the component of $B_{n+1}^+\cap\dpt$ which meets 
$\zeta_{n,p}$, and so $\Delta\p$ could not lie in $\bd\widehat{Q}_{n+1}^+$. 
Thus $\Delta\p$ lies in $\bd Q_{n+1}^+$. 

If $\theta\p$ lies in the component of $B_{n+1}^+\cap\dpt$ which meets 
$\zeta_{n,p}$, then $\sigma$ is isotopic to $\gamma_{n,p}\cup\rho_{n,p}\cup
\zeta_{n,p}$ and $\sigma^+$ is isotopic to a parallel copy of 
$\alpha_{n,p}\cup\zeta_{n,p}$. Then $\widehat{\tau}$ is equivalent to the 
tangle $\tau_{n,p}$ in $R_{n+1,p}$, and so the exterior of $\widehat{\tau}$ 
in $B_{n+1}^+$ is \birr. Thus $\bd \Delta=\bd \Delta\p$ for a disk $\Delta\p$ 
in $\bd \widehat{Q}_{n+1}^+$. If $\Delta\p$ does not lie in $\bd Q_{n+1}^+$, 
then for homological reasons it must contain both components of $N_n^+\cap 
\dpt$. This is impossible since then $\Delta\p$ would also contain the 
component of $B_{n+1}^+\cap \dpt$ which meets $\alpha_{n,p}$, and so 
$\Delta\p$ could not lie in $\bd Q_{n+1}^+$. Thus $\Delta\p$ lies in 
$\bd Q_{n+1}^+$. 

We have thus shown that our triple has the halfdisk property. 

We now consider the exterior $X_{n+2}^3$ of $N_{n+1}^3$ in $R_{n+2,p}$. 
This is the same as the exterior $X_{n+2}^2$ of $N_{n+1}^2$ in $R_{n+2,p}$ 
that we had in subcase (a). Thus $X_{n+2}^3$ is \irr, $O_p$, $S_{n+1}^3$, and 
$X_{n+2}^3\cap(\dpo\cup\dpt)$ are \inc\ in $X_{n+2}^3$, and $(X_{n+2}^3, 
X_{n+2}^3\cap(\dpo\cup\dpt),O_p\cup S_{n+1}^3)$ has the halfdisk property. 
We let $X_i\p=X_i^3$ and complete the proof as usual. 

\textit{Case 3:} $S_n$ is \inc\ in $X_n^0$, but 
$(X_n^0,X_n^0\cap(\dpo\cup\dpt),S_n)$ does not have the halfdisk property. 

There is a proper disk $\Delta^1$ in $X_n^0$ such that, say, $\Delta^1
\cap\dpo$ is an arc $\theta_1$ and $\Delta^1\cap S_n$ is an arc $\theta_1\p$ 
such that $\bd \Delta^1=\theta_1\cup\theta_1\p$ and neither arc is \bpl. 
$\Delta^1$ splits $N_n^0$ inot a \tm\ consisting of two 3-balls $B_n^-$ and 
$B_n$. Let $S_n^1$ be the surface obtained by boundary compressing $S_n$ along 
$\Delta^1$. It has two components $\Sigma_n^-$ and $\Sigma_n$, where 
$\Sigma_n^-$ is an annulus with $\bd \Sigma_n^-$ in \dpo\ and $\Sigma_n$ is 
a disk with two holes having one boundary component in \dpo\ and the other 
two in \dpt. Let $X_{n+1}^1$ be the exterior of $N_n^1\cup N_n^-\cup N_n^+$ 
in $N_{n+1}$. We choose the notation so that $\Sigma_n^-=B_n^-\cap X_{n+1}^1$ 
and $\Sigma_n=B_n\cap X_{n+1}^1$. Let $Q_n^-$ be the exterior of 
$\kappa\p\cap F_n^-$ in $B_n^-$, and let $Q_n$ be the exterior of 
$\kappa\p\cap B_n$ in $B_n$. Let $X_n^1=Q_n^-\cup Q_n$. 

We have that $\Sigma_n^-$ and $Q_n^-\cap \dpo$ are \inc\ in $Q_n^-$ and 
that $(Q_n^-,Q_n^-\cap(\dpo\cup\dpt,\Sigma_n^-)$ has the halfdisk property. 
Also $\Sigma_n$ and $Q_n\cap(\dpo\cup\dpt)$ are \inc\ in $Q_n$. 

\textit{Subcase (a):} $(Q_n,Q_n\cap(\dpo\cup\dpt),\Sigma_n)$ has the 
halfdisk property. 

Consider $X_{n+1}^1$. $A_n^-\cup A_n^+\cup \Sigma_n^-\cup \Sigma_n$ is 
\inc\ in $X_{n+1}^1$ for homological reasons, as is $X_{n+1}^1\cap
(\dpo\cup\dpt)$. The boundary of any compressing disk for $S_{n+1}$ in 
$X_{n+1}^1$ thus has to separate two of the components of $\bd S_{n+1}$ 
from the other two. So in $N_{n+1}$ the compressing disk must separate 
two of the components of $N_{n+1}\cap(\dpo\cup\dpt)$ from the other two 
components. But this is impossible since each pair of these disks is 
joined by some component of $N_n^-\cup N_n^+\cup N_n^1$. So $S_{n+1}$ 
is \inc\ in $X_{n+1}^1$. 

Suppose $(X_{n+1}^1,X_{n+1}^1\cap(\dpo\cup\dpt),S_{n+1}\cup A_n^-\cup A_n^+
\cup\Sigma_n^-\cup\Sigma_n)$ does not have the halfdisk property. Then there 
is a proper disk $\Delta$ in $X_{n+1}^1$ with $\Delta\cap(\dpo\cup\dpt)$ 
an arc $\theta$ and $\Delta\cap(S_{n+1}\cup\Sigma_n)$ an arc $\theta\p$ 
such that $\theta\cup\theta\p=\bd\Sigma$ and neither arc is \bpl. 

Assume $\theta\p$ lies in $\Sigma_n$. Then $\bd\theta\p$ lies in one 
component of $\bd \Sigma_n$ and $\theta\p$ separates the other two 
components. Then in the \tm\ $\widehat{X}_{n+1}^1$ obtained by adjoining 
$N_n^-\cup N_n^+\cup B_n^-$ to $X_{n+1}^1$ we have that $\bd\Delta$ is 
isotopic to one of these two components of $\bd\Sigma_n$. But this is 
impossible since $\Sigma$ is \inc\ in $\widehat{X}_{n+1}^1$ for 
homological reasons. 

Thus $\theta\p$ lies in $S_{n+1}$. $\theta$ lies in the unique component 
of $N_{n+1}\cap\dpo$ which meets both $B_n^-$ and $B_n$. $\Delta$ separates 
$N_n^-\cup B_n^-$ from $N_n^+\cup B_n$ in $N_{n+1}$. But $B_n^-$ contains 
the arc $\gamma_{n,p}\cup\delta_{n,p}$, and $B_n$ contains the arc 
$\alpha_{n,p}\cup\zeta_{n,p}$. Thus the tangle $\tau_{n,p}$ in $N_{n+1}$ is 
splittable, contradicting the fact that its exterior is \birr. 

Hence our triple has the halfdisk property. 

Let $X_{n+2}^1$ be the exterior of $N_{n+1}$ in $R_{n+2,p}$. We already 
know that it is \irr, $X_{n+2}^1\cap(\dpo\cup\dpt)$ and $S_{n+1}$ 
are \inc\ in $X_{n+2}^1$ and $(X_{n+2}^1,X_{n+2}^1\cap(\dpo\cup\dpt),
O_p\cup S_{n+1})$ has the halfdisk property. We let $X_i\p=X_i^1$ and 
complete the proof as usual.

\textit{Subcase (b):} $(Q_n,Q_n\cap(\dpo\cup\dpt),\Sigma_n)$ does not 
have the halfdisk property. 

There is a proper disk $\Delta^2$ in $Q_n$ with $\Delta^2\cap(\dpo\cup\dpt)$ 
an arc $\theta_2$ and $\Delta^2\cap\Sigma_n$ an arc $\theta_2\p$ such 
that $\bd\Delta^2=\theta_2\cup\theta_2\p$ and neither arc is \bpl. 
$\Delta$ splits $B_n$ into a pair of 3-balls $B_n^0$ and $B_n^+$. 
Let $N_n^2=B_n^-\cup B_n^0\cup B_n^+$. Let $X_{n+1}^2$ be the exterior 
of $N_n^-\cup N-N^+\cup N_n^2$ in $N_{n+1}$. Let $S_n^2$ be the 
surface obtained by boundary compressing $S_n^1$ along $\Delta^2$. 
It consists of three annuli $\Sigma_n^-$, $\Sigma_n^0$, and $\Sigma_n^+$ 
with $\Sigma_n^{\pm}=B_N^{\pm}\cap X_{n+1}^2$ and $\Sigma_n^0=B_n^0\cap 
X_{n+1}^2$. Let $Q_n^{\pm}$ and $Q_n^0$ be the exteriors of $\kappa\p\cap
B_n^{\pm}$ and $\kappa\p\cap B_n^0$, respectively. Let $X_n^2=
Q_n^-\cup Q_n^0\cup Q_n^+$. 

$X_n^2$ is \irr. $X_n^2\cap(\dpo\cup\dpt)$ and $S_n^2$ are \inc\ in $X_n^2$. 
$(X_n^2,X_n^2\cap(\dpo\cup\dpt),S_n^2)$ has the halfdisk property. 

Consider $X_{n+1}^2$. It is \irr. $A_n^-$, $A_n^+$, $\Sigma_n^-$, 
$\Sigma_n^0$, and $\Sigma_n^+$ are \inc\ in $X_{n+1}^2$ for homological 
reasons, as is $X_{n+1}^2\cap(\dpo\cup\dpt)$. Thus the boundary of any 
compressing disk for $S_{n+1}$ in $X_{n+1}^2$ must separate two of the 
components of $\bd S_{n+1}$ from the other two. Hence the compressing 
disk must separate two of the components of $N_{n+1}\cap(\dpo\cup\dpt)$ 
from the other two in $N_{n+1}$. But this is impossible since the union of 
these four disks with $N_n^2$ is connected. Thus $S_{n+1}$ is \inc\ in 
$X_{n+1}^2$. 

Suppose $(X_{n+1}^2,X_{n+1}^2\cap(\dpo\cup\dpt),S_{n+1}\cup A_n^-\cup 
A_n^+\cup\Sigma_n^-\cup\Sigma_n^0\cup\Sigma_n^+)$ does not have the 
halfdisk property. Then there is a proper disk $\Delta$ in $X_{n+1}^2$ 
which meets $\dpo\cup\dpt$ in an arc $\theta$ and $S_{n+1}$ in an arc 
$\theta\p$ such that $\theta\cup\theta\p=\bd\Delta$ and neither arc 
is \bpl. Let $N_{n+1}^{\pm}$ be the two 3-balls into which $N_{n+1}$ is 
split by $\Delta$. 

There are six patterns in which $B_n^-$, $B_n^0$, and $B_n^+$ 
can connect the four components of $N_{n+1}\cap(\dpo\cup\dpt)$. 
These patterns can be described as follows. 
Note that $C_n^{\PP}\cap N_{n+1}$ is a regular neighborhood of 
the union of a set of arcs $\alpha$, $\gamma$, $\delta$, $\zeta$, $\rho$, 
$\beta^-$, and $\beta^+$, where $(N_{n+1},\gamma\cup\delta\cup\rho\cup\alpha
\cup\zeta)$ is equivalent to $(R_{n+1,p},\gamma_{n,p}\cup\delta_{n,p}\cup
\rho_{n,p}\cup\alpha_{n,p}\cup\zeta_{n,p})$ and $\beta^-$ and $\beta^+$ 
are proper \bpl\ arcs in the complement of this graph. 
The 3-balls $B_n^-$, $B_n^0$, and $B_n^+$ are 
regular neighborhoods of arcs which are certain unions of the 
$\alpha$, $\gamma$, $\delta$, $\zeta$, and $\rho$ in $N_{n+1}$, 
where two such arcs are pushed apart slightly to make the 3-balls 
disjoint. 

\textit{(i)} $B_n^-$, $B_n^0$, and $B_n^+$ are regular neighborhoods of 
$\gamma\cup\delta$, $\gamma\cup\rho\cup\alpha$, and $\gamma\cup\rho\cup\zeta$, 
respectively. 

Then $N_n^-\cup B_n^-$ lies in $N_{n+1}^-$, and $N_n^+\cup B_n^0\cup B_n^+$ 
lies in $N_{n+1}^+$. Thus $\Delta$ separates $\gamma\cup\rho\cup\alpha$ 
and a copy of $\gamma\cup\delta$ in $N_{n+1}$. The exterior of this tangle 
is \hm\ to the exterior of the graph $\gamma\cup\delta\cup\rho\cup\alpha$ in 
$N_{n+1}$ which in turn is \hm\ to the exterior of the tangle 
$\tau=(\gamma\cup\delta)\cup(\alpha\cup\zeta)$ in $N_{n+1}$ and is 
therefore \birr. So $\bd \Delta=\bd\Delta\p$ for a disk $\Delta\p$ in the 
boundary of the exterior of $B_n^-\cup B_n^0$ in $N_{n+1}$. But this is 
impossible since $\bd \Delta$ splits the boundary of this space into a 
pair of punctured tori. So this case cannot occur.

\textit{(ii)} $B_n^-$, $B_n^0$, and $B_n^+$ are regular neighborhoods of 
$\gamma\cup\delta$, $\delta\cup\rho\cup\zeta$, and $\delta\cup\rho\cup
\alpha$, respectively. 

Then $N_n^-\cup B_n^-$ lies in $N_{n+1}^-$, and $N_n^+\cup B_n^0\cup B_n^+$ 
lies in $N_{n+1}^+$. Thus $\Delta$ separates $\delta\cup\rho\cup\alpha$ and 
a copy of $\gamma\cup\delta$ in $N_{n+1}$. The exterior of this tangle 
is \hm\ to the exterior of the tangle $\tau=(\gamma\cup\delta)\cup
(\alpha\cup\zeta)$ in $N_{n+1}$ and so is \birr. So this case cannot occur. 

\textit{(iii)} $B_n^-$, $B_n^0$, and $B_n^+$ are regular neighborhoods of 
$\gamma\cup\delta$, $\gamma\cup\rho\cup\alpha$, and $\alpha\cup\zeta$, 
respectively. 

There are then two possibilities. The first is that $N_n^-\cup B_n^-$ 
lies in $N_{n+1}^-$ and $N_n^+\cup B_n^0\cup B_n^+$ lies in $N_{n+1}^+$. 
Thus $\Delta$ separates $\gamma\cup\delta$ and $\alpha\cup\zeta$ in $N_{n+1}$. 
This is impossible since the exterior of this tangle is \birr. The second 
is that $N_n^-\cup B_n^-\cup B_n^0$ lies in $N_{n+1}^-$ and $N_n^+\cup B_n^+$ 
lies in $N_{n+1}^+$. Again $\Delta$ separates $\gamma\cup\delta$ and 
$\alpha\cup\zeta$ in $N_{n+1}$, so this case cannot occur. 

\textit{(iv)} $B_n^-$, $B_n^0$, and $B_n^+$ are regular neighborhoods of 
$\gamma\cup\delta$, $\delta\cup\rho\cup\zeta$, and $\alpha\cup\zeta$, 
respectively. 

The first possibility is that $N_n^-\cup B_n^-$ lies in $N_{n+1}^-$ and 
$N_n^+\cup B_n^+\cup B_n^0$ lies in $N_{n+1}^+$. The second is that 
$N_n^-\cup B_n^-\cup B_n^0$ lies in $N_{n+1}^-$ and $N_n^+\cup B_n^+$ 
lies in $N_{n+1}^+$. As in (iii) either of these implies that 
$\Delta$ separates $\gamma\cup\delta$ and $\alpha\cup\zeta$ in $N_{n+1}$, 
so this case cannot occur. 

\textit{(v)} $B_n^-$, $B_n^0$, and $B_n^+$ are regular neighborhoods of 
$\gamma\cup\delta$, $\gamma\cup\rho\cup\zeta$, and $\alpha\cup\zeta$, respectively. 

The first possibility is that $N_n^-\cup B_n^-$ lies in $N_{n+1}^-$ and 
$N_n^+\cup B_n^0\cup B_n^+$ lies in $N_{n+1}^+$. The second is that 
$N_n^-\cup B_n^-\cup B_n^0$ lies in $N_{n+1}^-$ and $N_n^+\cup B_n^+$ 
lies in $N_{n+1}^+$. Thus $\Delta$ separates $\gamma\cup\delta$ and 
$\alpha\cup\zeta$ in $N_{n+1}$, so this case cannot occur. 

\textit{(vi)} $B_n^-$, $B_n^0$, and $B_n^+$ are regular neighborhoods 
of $\gamma\cup\delta$, $\delta\cup\rho\cup\gamma$, and $\alpha\cup\zeta$, 
respectively. 

The first possibility is that $N_n^-\cup B_n^-$ lies in $N_{n+1}^-$ and 
$N_n^+\cup B_n^0\cup B_n^+$ lies in $N_{n+1}$. The second is that 
$N_n^-\cup B_n^-\cup B_n^0$ lies in $N_{n+1}^-$ and $N_n^+\cup B_n^+$ lies 
in $N_{n+1}^+$. Again $\Delta$ separates $\gamma\cup\delta$ and $\alpha\cup
\zeta$ in $N_{n+1}$, so this case cannot occur. 

So our triple has the halfdisk property. 

Let $X_{n+2}^2$ be the exterior of $N_{n+1}$ in $R_{n+2,p}$. We already 
know that it is \irr, $X_{n+2}^2\cap(\dpo\cup\dpt)$ and $S_{n+1}$ are 
\inc\ in $X_{n+2}^2$ and $(X_{n+2}^2,X_{n+2}^2\cap(\dpo\cup\dpt),O_p
\cup S_{n+1})$ has the halfdisk property. We let $X_i\p=X_i^2$ and 
complete the proof as usual. \end{proof}
  
%Section 14
\section{Isotopy classification of the $V^{\PP}$}

%Proposition 14.1
\begin{prop} Let $\PP$ and $\QQ$ be finite non-empty sets of integers. 
$V^{\PP}$ and $V^{\QQ}$ are isotopic if and only if $\PP=\QQ$. \end{prop}

\begin{proof}Let $\PP=\{p_1,\ldots, p_m\}$ and $\QQ=\{q_1,\ldots, q_n\}$ 
with their natural orderings. Suppose $h_t:V^{\PP}\ra W$ is an isotopy 
with $h_0$ the inclusion map of $V^{\PP}$ into $W$ and $h_1(V^{\PP})
=V^{\QQ}$. Suppose that there is a $p\in\PP$ such that $p\notin\QQ$. 
We may assume that $V^p$ is an \er\ of $W$ at a knot $\ka$. Then $h_1(V^p)$ 
is an \er\ of $W$ at $h_1(\ka)$, and $h_1(\ka)\sbs h_1(V^p)\sbs V^{\QQ}$.

First suppose that $q_1\leq p\leq q_n$. Then as in the proof of Theorem 10.7 
there is a plane $\Pi$ in $V^{\QQ}$ which is proper and non-trivial in 
$V^{\QQ}$ such that $V^{\QQ}-\Pi$ has components isotopic to $V^{\RRRR}$ 
and $V^{\SSSS}$, where $\RRRR=\QQ\cap[q_1,p]$ and $\SSSS=\QQ\cap[p,q_n]$. 
Since $h_1(V^p)$ is \ei\  rel $h_1(\ka)$ in $W$ we have that $h_1(V^p)$ is 
\ei\ rel $h_1(\ka)$ in the smaller set $V^{\QQ}$. Thus by Theorem 6.2 
$h_1(V^p)$ can be isotoped to lie in $V^{\RRRR}$ or $V^{\SSSS}$. Hence we 
may assume that $\QQ$ has the property that either $q<p$ for all $q\in\QQ$ 
or $q>p$ for all $q\in\QQ$. 

There is an $i\geq0$ such that $h_1(\ka)\sbs \inte C^{\QQ}_i$. 
Since $V^p$ has genus one and is therefore \rirr\ we have by Theorem 13.1 
that $h_1(V^p)$ is isotopic to $V^{\RRRR}$ for some good subset of $\QQ$. 
Since $h_1(V^p)$ has genus one we have by Proposition 11.1 that $\RRRR$ has a 
single element. Thus it suffices to prove the following result. \end{proof}

\begin{lem}If $p\neq q$, then $V^p$ is not isotopic to $V^q$. \end{lem}

\begin{proof} We may assume that $p<q$ and that $V^p$ is an \er\ of $W$ 
at a knot \ka\ in $V^p$. Suppose $h_t:V^p\ra W$ is an isotopy with $h_0$ the 
inclusion map of $V^p$ into $W$ and $h_1(V^p)=V^q$. 

Let $T=\cup_{t\in[0,1]}h_t(\ka)$ be the track of \ka\ under this isotopy. 
There exist integers $n$ and $r\leq s$ such that $T\sbs\inte C_n^{[r,s]}$. 
Note that $r\leq p$ and $q\leq s$. By the Covering Isotopy Theorem 
\cite{Ce, EK}  
there is an ambient isotopy $g_t:W\ra W$ such that $g_0$ is the identity 
of $W$, $g_t(x)=h_t(x)$ for all $x\in \ka$ and $t\in[0,1]$, and $g_t(x)=x$ 
for all $x\in W-\inte C_n^{[r,s]}$ and $t\in[0,1]$. 

Recall that $C_n^{[r,s]}$ is the union of all those $R_{n,j}$ with 
$r-1\leq j\leq s$ and all those $L_{n,j}$ and $H_{n,j}$ with $r\leq j\leq s$. 
Now $C_n^p=R_{n,p-1}\cup L_{n,p}\cup H_{n,p}\cup R_{n,p}$,  and 
$C^q_n=R_{n,q-1}\cup L_{n,q}\cup H_{n,q}\cup R_{n,q}$. We have that 
$\ka\sbs\inte C^p_n$, and the annulus $\bd H_{n,p}-\inte(H_{n,p}\cap L_{n,p})$ 
is \inc\ in $C^p_n-\ka$. Since $C^p_n$ meets the rest of $C_n^{[r,s]}$ in 
disjoint disks we have that this annulus is \inc\ in $C_n^{[r,s]}-\ka$. 
Since $g_1$ is the identity on $\bd C_n^{[r,s]}$, we have that the annulus 
is \inc\ in $C_n^{[r,s]}-g_1(\ka)$. But $g_1(\ka)=h_1(\ka)\sbs V^q$, and 
so $g_1(\ka)$ does not meet $H_{n,p}$, from which it follows that the annulus 
must be compressible in $C_n^{[r,s]}-g_1(\ka)$. This contradiction completes 
the proof. \end{proof}

%Section 15
\section{Homeomorphism classification of the $V^p$}

Until now the sense of the Whitehead clasp in the 1-handle $H_{n,p}$ in 
our construction of $W$ has been immaterial. We will now modify our 
construction by choosing different senses for the clasp depending on $n$ 
and $p$. This will be used to construct uncountably many pairwise 
non-\hm\ $W$ which cover 3-manifolds with infinite cyclic fundamental 
group such that the only 3-manifolds non-trivially covered by them must 
have infinite cyclic fundamental group. This modification will also be 
used to construct uncountably many pairwise non-\hm\ $W$ which cannot 
non-trivially cover any \tm. 

Let $C$ be a solid torus with oriented meridian $m$ and longitude $\ell$ on 
$\bd C$ as shown in Figure 9. It is regarded as embedded in 
$\RRR\sbs S^3$ in the manner shown. We will blithely confuse an oriented 
simple closed curve on $\bd C$ and its homology class in $H_1(\bd C)$. 
Let $C^+$ and $C^-$ be the solid tori in the interior of $C$ as indicated in 
Figures 9(a) and 9(b), respectively. 

%INSERT FIGURE 9 HERE

\begin{figure}
\epsfig{file=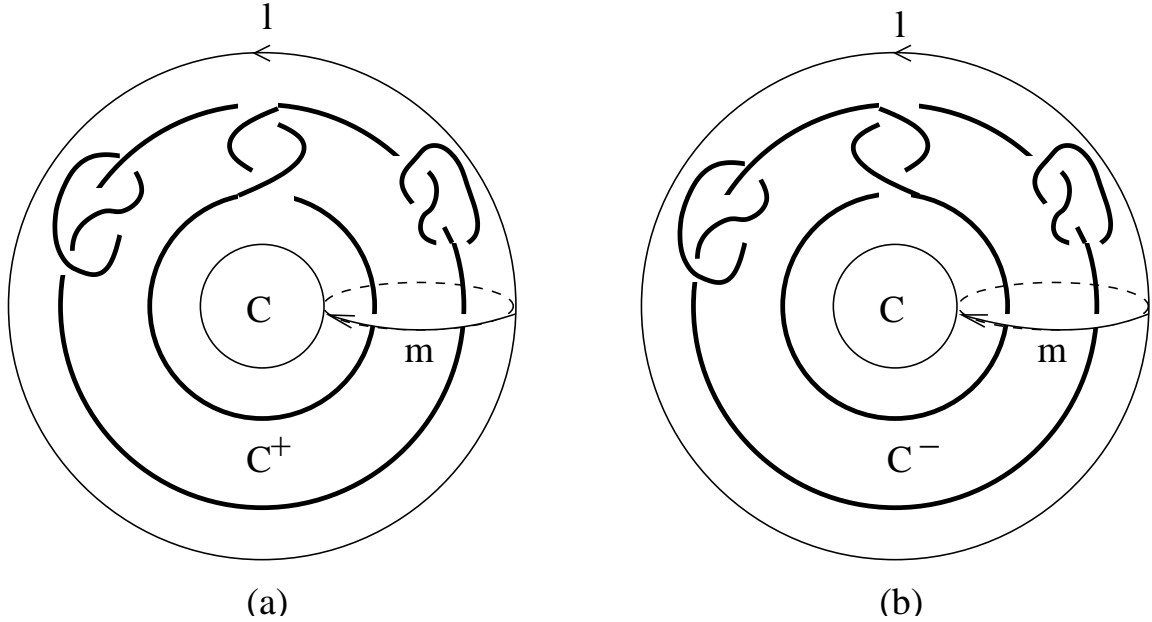, width=6in}
\caption{Granny Whitehead solid tori $C^+$ and $C^-$}
\end{figure}

%LEMMA 15.1
\begin{lem} There is no homeomorphism $h:C\ra C$ such that $h(C^+)=C^-$ 
and $h(\ell)=\pm\ell$. \end{lem}

\begin{proof} Since $h(\ell)=\pm\ell$ we have that $h$ extends from $C$ 
to $S^3$. It must be orientation preserving since otherwise the granny knot is 
invariant under an orientation reversing homeomorphism of $S^3$. This 
cannot occur because the granny knot has signature $\pm 4$ and is therefore 
non-amphicheiral. (See Rolfsen \cite{Ro}.) Thus $h(m)=\pm m$, with the same sign 
as $h(\ell)=\pm\ell$. 

Let $t:C\ra C$ be the homeomorphism obtained by cutting $C$ along a 
meridinal disk, twisting, and regluing so that $t(m)=m$ and $t(\ell)=\ell+m$. 
The results of applying $t$ to $C^+$ and $C^-$ are shown in Figure 10. 

%INSERT FIGURE 10 HERE

\begin{figure}
\epsfig{file=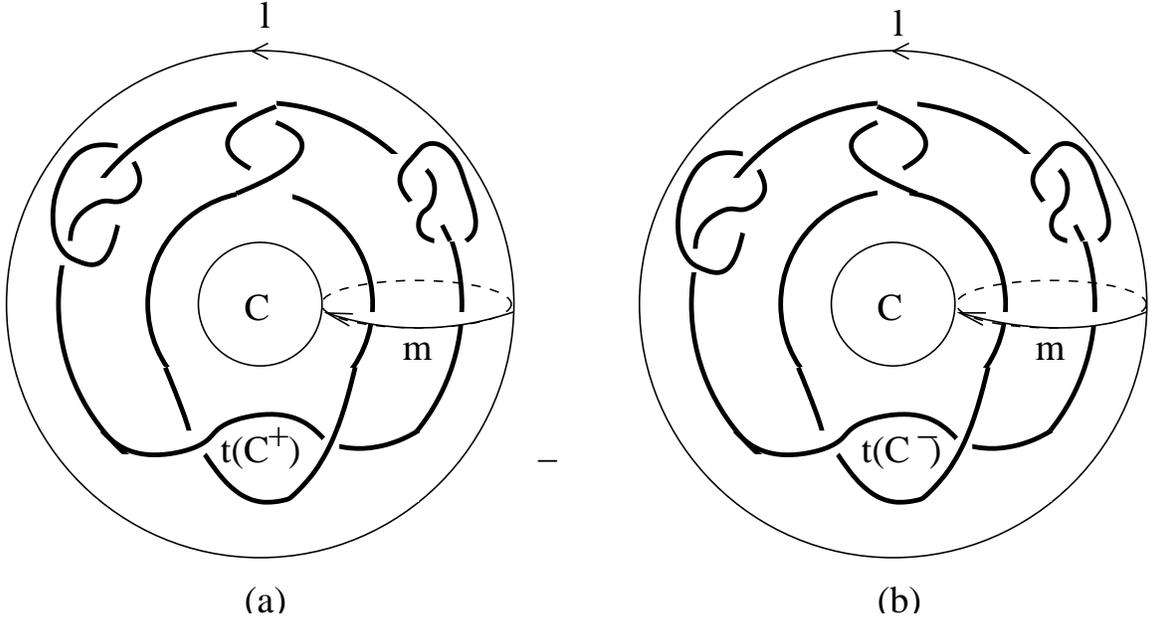, width=6in}
\caption{Twisted Granny Whitehead solid tori $t(C^+)$ and $t(C^-)$}
\end{figure}

Now $t\n(\ell)=\ell-m$, and so $tht\n(\ell)=th(\ell-m)=t(\pm\ell\mp m)=
\pm(\ell+m)\mp m=\pm\ell\pm m \mp m = \pm \ell$. Thus $tht\n$ extends 
to a homeomorphism of $S^3$ carrying $t(C^+)$ to $t(C^-)$. But this 
is impossible sinice the core of $t(C^+)$ is the sum of a granny knot and 
a figure eight knot while the core of $t(C^-)$ is the sum of a granny 
knot and a trefoil knot. These knots can be distinguished by, for example, 
their Alexander polynomials. \end{proof}

We remark that with slightly more work the condition $h(\ell)=\pm\ell$ 
can be dropped from the hypotheses of this lemma. However, we will not 
need this stronger result. 

Now let $\varphi:\mathbf{N}\ra\{\pm1\}$ be a function. We construct a 
contractible open \tm\ $V[\varphi]$ by specifying an exhaustion 
$\{C_n\}_{n\geq0}$ where the model for the pair $(C_{n+1},C_n)$ is 
$(C,C^{\pm})$ for $\varphi(n)=\pm1$. Note that $V[\varphi]$ is a 
modified version of $V^p$. 

%Proposition 15.2
\begin{prop} $V[\varphi]$ and $V[\varphi\p]$ are \hm\ if and only if 
there are integers $m_0$, $n_0\geq0$ such that $\varphi(m_0+i)=
\varphi\p(n_0+i)$ for all $i\geq0$. \end{prop}

\begin{proof} Denote $V[\varphi]$ and $V[\varphi\p]$ by $V$ and $V\p$, 
respectively. Distinguish their defining exhaustions by $\{C_n\}$ 
and $\{C_n\p\}$, respectively. Clearly the condition is sufficient, so 
assume that $h:V\ra V\p$ is a homeomorphism. 

Choose $p\geq0$ such that $h(C_0)\sbs\inte C\p_p$. Then choose $q\geq1$ such 
that $C\p_p\sbs\inte h(C_q)$. Now $h(\bd C_q)$ is \inc\ in $V\p-h(C_0)$ and 
so is \inc\ in the smaller space $V\p-C_p\p$. Put $h(\bd C_q)$ in \mgp\ 
with respect to $\cup_{n>p} \bd C\p_n$. Then no component of the 
intersection bounds a disk on either surface. If the intersection is 
non-empty, then let $n$ be the smallest $n>p$ such that $h(\bd C_q)\cap 
\bd C\p_n\neq\ns$. There is then a component $A$ of $h(\bd C_q)\cap C\p_n$ 
which is an annulus lying in $Y\p=C_n\p-\inte C_{n-1}\p$. Recall that 
$Y\p$ is the union of a Whitehead link exterior $Y_0\p$ and two disjoint 
trefoil knot exteriors $Y\p_1$ and $Y\p_2$ with $Y_0\p\cap Y_j\p=A\p_j$ 
an annulus, $j=1$, $2$. As in the proof of Lemma 9.6 we have that $A$ is 
\pl\ in $Y\p$ to an annulus in $\bd C\p_n$, contradicting minimality. 
Thus we have that the intersection is empty. 

So $h(\bd C_q)$ lies in $Y\p=C\p_n-\inte C\p_{n-1}$ for some $n>p$. 
Put it in \mgp\ with respect to $A\p_1\cup A\p_2$. As usual there are 
no trivial intersection curves. 

\textit{Case 1:} The intersection is non-empty. Then there is an annulus 
component $A$ of, say, $h(\bd C_q)\cap Y_1\p$. Since it has meridian 
boundary components on $\bd Y_1\p$ it must be \pl\ in $Y\p_1$ to an 
annulus $A\p$ in $\bd Y_1\p$. By minimality $A\p$ cannot lie in $A_1\p$, so 
it must contain $Y_1\p\cap\bd C\p_{n-1}$. By minimality we may assume that no 
component of $h(\bd C_q)\cap Y_1\p$ lies between $A$ and 
$Y_1\p\cap \bd C\p_{n-1}$ in the product $I$ bundle joining $A$ and 
$Y_1\p\cap \bd C_{n-1}\p$. 

Consider the annuli $A_0=h(\bd C_q)-\inte A$ and 
$A_0\p=(Y_1\p\cap \bd C_{n-1}\p)\cup(A\p\cap A_1\p)$. 
$A\p_0$ is a proper \inc\ annulus in $h(C_q-\inte C_0)$ with 
$\bd A_0\p$ in $h(\bd C_q)$. It follows from the anannularity of 
the Whitehead link exterior \cite{My simple} that $A\p_0$ is \bpl\ in this space. 
It cannot be \pl\ to $A$ since $A\cup A_0\p$ bounds a solid torus 
containing $C_{n-1}\p$ and hence $h(C_0)$. Therefore $A_0\p$ is 
\pl\ to $A_0$. It follows that $h(\bd C_q)$ is \pl\ in $Y\p$ to 
$\bd C\p_{n-1}$. We then isotop $h$ so that $h(C_q)=C_{n-1}\p$. 
In this case we let $m_0=q$ and $n_0=n-1$. 

\textit{Case 2:} $h(\bd C_q)\cap(A_1\p\cup A_2\p)=\ns$. Then 
$h(\bd C_q)$ lies in $Y_0\p$. Since $Y_0\p$ is atoroidal $h(\bd C_q)$ is 
\bpl\ in $Y_0\p$. We isotop $h$ so that $h(\bd C_q)$ is a component of 
$\bd Y_0\p$. 

Suppose $h(\bd C_q)=\bd Y_0\p-\bd C_m\p$. Then $A_1\p$ is a proper 
\inc\ annulus in $h(C_q-\inte C_0)$ with $\bd A_1\p$ in $h(\bd C_q)$, so it 
must be \bpl\ in this space. But this is impossible since $A_1\p$ splits 
this space into the trefoil knot space $Y_1\p$ and a space with two 
boundary components. So this situation cannot occur. 

Thus $h(\bd C_q)=\bd C_n\p$. In this case we let $m_0=q$ and $n_0=n$. 

We now have that $h(C_{m_0})=C_{n_0}\p$. Using the methods just 
employed we isotop $h$ rel $C_{m_0}$ so that $h(C_{m_0+1})=C_{n_0+1}\p$. 
By Lemma 15.1 we have that $\varphi(m_0)=\varphi\p(n_0)$. We then isotop 
$h$ rel $C_{m_0+1}$ so that $h(C_{m_0+2})=C_{n_0+2}\p$ and conclude 
that $\varphi(m_0+1)=\varphi\p(n_0+1)$. We then continue this process 
to get that $\varphi(m_0+i)=\varphi\p(n_0+i)$ for all $i\geq0$. 
\end{proof}

%Section 16
\section{The complex of end reductions}

Let $W$ be a Whitehead manifold. An end reduction $V$ of $W$ is 
\textit{minimal} if whenever $U$ is an \er\ of $W$ 
which is contained in $V$ we have that $U$ is isotopic to $V$. 

%Theorem 16.1
\begin{thm} Genus one end reductions of $W$ are minimal and \rirr. \end{thm}

\begin{proof} Suppose $V$ is a genus one end reduction of $W$. 
Then $V$ has an exhaustion $\{C_n\}_{n\geq0}$ such that each 
$\bd C_n$ is a torus. $\bd C_n$ is \inc\ in $W-\inte C_0$. By Lemma ? 
$V$ is a Whitehead manifold, and so $\bd C_n$ must be compressible in $C_n$. 
Since $V$ is \irr\ and $V-\inte C_n$ is non-compact $C_n$ must be a 
solid torus. 

Let $J$ be a regular submanifold of $W$ such that $J\sbs V$. Then 
$J\sbs \inte C_n$ for some $n$. If $\bd C_n$ is compressible in 
$C_n-J$, then $J$ lies in a 3-ball, contradicting the fact that it is 
regular. Therefore $\bd C_n$ is \inc\ in $C_n-J$ and hence in 
$W-J$. By Lemma X we then have that any end reduction of $W$ at $J$ is 
isotopic to $V$. 

The \RR-irreducibility of genus one Whitehead manifolds was proven 
by Kinsoshita \cite{Kn}. \end{proof}

We remark that \rirr\ end reductions need not be minimal (any $V^{\PP}$ 
with $\PP$ a good set having more than one element is an example), and 
minimal end reductions need not be \rirr\ (the double of the Tucker 
manifold \cite{Tu} can be shown to be an example). 

We now define the \textit{simplicial complex of minimal, \rirr\ 
end reductions of $W$}, denoted $\mathcal{S}(W)$. 

The vertices of $\SSSS(W)$ are the isotopy classes $[V]$ of minimal, 
\rirr\ end reductions of $W$. 

Two distinct vertices $[V_0]$ and $[V_1]$ are joined by an edge if there 
is an \rirr\ end reduction $E_{0,1}$ of $W$ having the following properties: 
(1) $E_{0,1}$ contains representatives of $[V_0]$ and $[V_1]$. 
(2) Every \rirr\ \er\ of $W$ which is contained in $E_{0,1}$ is 
isotopic to $V_0$, $V_1$, or $E_{0,1}$. (3) Among \rirr\ \er s of $W$ 
one has that $E_{0,1}$ is unique up to isotopy with respect to 
properties (1) and (2). 

Three distinct vertices $[V_0]$, $[V_1]$, and $[V_2]$ span a 2-simplex if 
each pair of vertices is joined by an edge and there is an \rirr\ 
end reduction $T_{0,1,2}$ of $W$ having the following properties: 
(1) $T_{0,1,2}$ contains representatives of $[V_0]$, $[V_1]$, $[V_2]$, 
$[E_{0,1}]$, $[E_{1,2}]$, and $[E_{2,0}]$. (2) Every \rirr\ \er\ of $W$ 
which is contained in $T_{0,1,2}$ is isotopic to one of these six \er s 
or $T_{0,1,2}$. (3) Among \rirr\ \er s $T_{0,1,2}$ is unique with respect 
to properties (1) and (2). 

There is an obvious generalization of these definitions which inductively 
defines simplices of higher dimension.  

%Theorem 16.2
\begin{thm} Let $W$ be a member of the family $\FF$ of Whitehead manifolds 
constructed in section 9. 
Then $\SSSS(W)$ is isomorphic to a triangulation of $\R$. \end{thm}

\begin{proof} Each $V^p$ has genus one and so by Theorem 16.1 is minimal 
and \rirr. By Theorem 14.1 $V^p$ and $V^q$ are isotopic if and only if 
$p=q$. 

Suppose $V$ is \rirr\ and minimal. By Theorem 13.1 $V$ is isotopic to $V^{\PP}$ 
for some good set $\PP$. If $\PP$ has more than one element, then for 
$p\in\PP$ we have by Theorem 11.1 that $V^p$ is not isotopic to $V^{\PP}$, 
and so $V$ is not minimal. This contradiction implies that $V$ must be 
isotopic to some $V^p$. We thus can associate the vertices of $\SSSS(W)$ 
with the integers. 

Given $p$, let $\RRRR=\{p,p+1\}$. Then $V^p$ and $V^{p+1}$ are contained 
in $V^{\RRRR}$. By Corollary 10.5 $V^{\RRRR}$ is \rirr. Let $V$ be an 
\rirr\ \er\ of $W$ which is contained in $V^{\RRRR}$. By Theorem 13.1 $V$ is 
isotopic to $V^{\QQ}$ for some subset $\QQ$ of $\RRRR$. Hence $V$ is isotopic  
to $V^p$, $V^{p+1}$, or $V^{\RRRR}$. We will see in the next paragraph that 
$V^{\RRRR}$ is unique up to isotopy with respect to these properties, and 
thus the vertices corresponding to $p$ and $p+1$ are joined by an edge 
in $\SSSS(W)$. 

Consider integers $p<q$. 
Suppose $E$ is an \rirr\ \er\ which contains representatives of 
$[V^p]$ and $[V^q]$ and has the property that every end reduction of 
$W$ which is contained in $E$ is isotopic to $V^p$, $V^q$, or $E$. 
Theorem 13.1 says that $E$ must be isotopic to $V^{\RRRR}$ for some 
good set $\RRRR$. Thus $V^p$ and $V^q$ must be isotopic to 
subsets of $V^{\RRRR}$. Since they are \rirr\ by Theorem 13.1 they 
are isotopic to $V^{\PP}$ and $V^{\QQ}$, respectively, where \PP\ and \QQ\ 
are good subsets of \RRRR. Since they have genus one these sets must each 
have one element. By Corollary 10.5 we then have that $\PP=\{p\}$ and 
$\QQ=\{q\}$. Thus $p,\,q\in\RRRR$. If $\RRRR\neq\{p,q\}$ then choose 
a third element $r\in\RRRR$. By Lemma 14.2 $V^r$ is not isotopic to 
$V^p$ or $V^q$. Since $V^r$ has genus one and $V^{\RRRR}$ does not 
we have that $V^r$ is not isotopic to $V^{\RRRR}$. This contradicts 
a property of $E$. Thus $\RRRR=\{p,q\}$. If $q\neq p+1$, then by 
Theorem 10.7 $V^{\RRRR}$ is not \rirr. This contradicts another property 
of $E$. Thus $q=p+1$, $\RRRR=\{p,p+1\}$ and $E$ is isotopic to $V^{\RRRR}$. 

Thus two vertices are joined by an edge if and only if they correspond 
to consecutive integers. It follows that there are no higher dimensional 
simplices, and so $\SSSS(W)$ is isomorphic to a triangulation of $\R$. 
\end{proof} 

%Section 17
\section{Covering translations}

%Theorem 17.1
\begin{thm} Let $W$ be a Whitehead manifold, and let $\SSSS(W)$ be 
the simplicial complex of minimal, \rirr\ \er s of $W$. Suppose 
$W$ is a covering space of a \tm\ $M$. Then $\pi_1(M)$ is 
isomorphic to a fixed point free, torsion free group of simplicial 
automorphisms of $\SSSS(W)$. \end{thm}

\begin{proof} Since $W$ is contractible $M$ is a finite dimensional 
$K(\pi,1)$, implying that $\pi_1(M)$ has finite cohomological dimension. 
Since groups with torsion have infinite cohomological dimension \cite{He} 
have that $\pi_1(M)$ is torsion free. 

$\pi_1(M)$ is isomorphic to the group of covering translations of 
the covering map $p:W\rightarrow M$. From the definition of $\SSSS(W)$ 
it is immediate that any group of homeomorphisms of $W$ induces a group 
of simplicial automorphisms. If some point of $\SSSS(W)$ were fixed by 
the automorphism $\ga$ associated to a non-trivial covering translation $g$, 
then the simplex of minimum dimension containing that point would be carried 
to itself by $\ga$. In particular the vertices of the simplex would be  
permuted by $\ga$. So some power $\ga^k$, $k>0$,  would fix a vertex 
$[V]$. This would imply that $g^k(V)$ is isotopic to $V$. By the Main 
Theorem X, $g^k$ would be the identity. Since $\pi_1(M)$ is torsion 
free this would imply that $g$ is trivial, a contradiction. \end{proof}

%Theorem 17.2
\begin{thm} There are uncountably many pairwise non-homeomorphic 
\RR-ir\-re\-du\-ci\-ble Whitehead manifolds $W$ which are covering spaces 
of 3-manifolds $W^{\#}$ with $\pi_1(W^{\#})$ infinite cyclic and have 
the property that whenever $W$ non-trivially covers a 3-manifold $M$ 
one has that $\pi_1(M)$ is infinite cyclic. \end{thm}

\begin{proof} Choose a function $\varphi:\mathbf{N}\rightarrow\{\pm1\}$. 
We let $W^{\#}[\varphi]$  be the 3-manifold with $\pi_1(W^{\#}[\varphi])$ 
infinite cyclic which was 
constructed in Section 9, where the sense of the clasp in the 1-handle 
$H_n$ is 
determined by $\varphi(n)$ as explained in Section 15. We let $W[\varphi]$ 
be its universal covering space. All of the $V^p$ in $W[\varphi]$ 
are homeomorphic to the genus one Whitehead manifold $V[\varphi]$ defined 
in Section 15. If $W[\varphi]$ and $W[\varphi\p]$ are \hm, then the 
minimal \rirr\ end reductions of one must be \hm\ to those of the 
other, and so we have that $V[\varphi]$ and $V[\varphi\p]$ are 
homeomorphic. By Theorem 15.2 we have that there exist $m_0$ and $n_0$ such that 
for all $i\geq0$ one has that $\varphi(m_0+i)=\varphi\p(n_0+i)$. 
There are uncountably many equivalence classes of functions 
$\varphi$ under this relation and so uncountably many pairwise 
non-\hm\ $W[\varphi]$. 

(More concretely, given any sequence $\{s_k\}$ 
of $1$'s and $-1$'s construct the sequence $\{t_n\}$ for which 
$t_0=s_0$, $t_1$ and $t_2$ equal $s_1$, $t_3$ through $t_6$ equal $s_2$, 
$t_7$ through $t_{14}$ equal $s_3$, etc. For each $k$ one has a block 
of $2^k$ copies of $s_k$. Then let $\varphi(n)=t_n$. Two sequences 
$\{s_k\}$ and $\{s_k\p\}$ are equal if and only if the corresponding 
sequences $\{t_n\}$ and $\{t_n\p\}$ are equal and thus are equal if and 
only if the functions $\varphi$ and $\varphi\p$ are equal. But 
$\varphi(m_0+i)=\varphi\p(n_0+i)$ if and only if $\varphi=\varphi\p$.)

Now suppose that one of these $W$ non-trivially covers a 3-manifold $M$.      
By the previous theorem $\pi_1(M)$ is isomorphic to a subgroup of the 
simplicial automorphism group of $\SSSS(W)$. Since $\SSSS(W)$ is 
isomorphic to a triangulation of $\R$, its automorphism group is 
isomorphic to the infinite dihedral group, whose only non-trivial 
torsion free subgroups are infinite cyclic. \end{proof}

\begin{thm} There are uncountably many pairwise non-\hm\ 
\RR-ir\-re\-du\-ci\-ble, 
non-eventually end-irreducible Whitehead manifolds $W$ which cannot 
non-trivially cover any 3-manifold. \end{thm}

\begin{proof} We modify the construction of $W[\varphi]$ as follows. 
Choose a function $\psi:\mathbf{N}\rightarrow{\pm1}$ such that 
$V[\varphi]$ is not \hm\ to $V[\psi]$. Replace the original manifold $V^0$ 
which was \hm\ to $V[\varphi]$ by a new $V^0$ \hm\ to $V[\psi]$. 
This can be done by 
changing the clasps in the 1-handles $H_{n,0}$ in the expression of 
$W[\varphi]$ as a monotone union of infinite genus handlebodies. Call 
the new manifold $W$. 

Any homeomorphism $g$ of $W$ must induce an automorphism $\ga$ of  
$\SSSS(W)$ which fixes $[V^0]$. Thus $g(V^0)$ is isotopic to $V^0$. 
If $g$ were a covering translation it would therefore have to be 
the identity. \end{proof}

\end{document}